\documentclass[12pt,leqno]{article}
\usepackage{amsthm,amssymb,latexsym,amsxtra}
\input diagrams.tex
\diagramstyle[small,labelstyle=\scriptstyle,midshaft]

\theoremstyle{plain}
\newtheorem{thm}{Theorem}[section]
\newtheorem{cor}[thm]{Corollary}
\newtheorem{lem}[thm]{Lemma}
\newtheorem{pro}[thm]{Proposition}

\theoremstyle{definition}
\newtheorem{dfn}[thm]{Definition}
\newtheorem{exa}[thm]{Example}

\newtheorem{que}[thm]{Question}
\newtheorem{rem}[thm]{Remark}
\newtheorem{block}[thm]{}

\newcommand{\lto}{\longrightarrow}
\newcommand{\eto}{\hookrightarrow}

\DeclareMathOperator{\End}{End}
\DeclareMathOperator{\Sing}{Sing}

\begin{document}
\title{\bf Hurwitz spaces of triple coverings of elliptic curves
 and moduli spaces of abelian threefolds}
\date{}
 
\author{Vassil Kanev}

\maketitle

\begin{flushright}
\emph{Dedicated to the memory of Fabio Bardelli}
\end{flushright}

\def\thefootnote{}
\footnote{Research supported by the italian MIUR under the  program
"Geometria sulle variet\`{a} algebriche".
 }


\begin{abstract}
\noindent
We prove that the moduli spaces $\mathcal{A}_3(D)$ of polarized abelian 
threefolds  with polarizations of types $D=(1,1,2), (1,2,2), (1,1,3)$ or 
$(1,3,3)$ are unirational. The result is based on the study of 
families
of simple coverings of elliptic curves of degree 2 or 3 and on the 
study of the corresponding period 
mappings associated with holomorphic differentials  with trace 0. In 
particular we prove  the unirationality of the Hurwitz space 
$\mathcal{H}_{3,A}(Y)$ which parameterizes simply branched triple 
coverings of an elliptic curve $Y$ with  determinants of the 
Tschirnhausen modules  isomorphic to  $A^{-1}$. 
\end{abstract}

\section*{Introduction}
The problem of calculating the Kodaira dimension of the 
moduli spaces $\mathcal{A}_{g}(d_{1},d_2,\ldots,d_g)$ which parameterize abelian varieties of dimension $g$ with  polarizations of type 
$D=\linebreak
(d_{1},d_2,\ldots,d_g),\; d_{i}|d_{i+1}$ has been a topic of intensive study in the last 20 years. 
The case of principal polarizations is almost settled.
Due to the work of Clemens, Donagi, Tai, Freitag and Mumford it is known that   the moduli space $\mathcal{A}_{g}
= \mathcal{A}_{g}(1,\ldots,1)$  is unirational if $g\leq 5$ and is of general type if $g\geq 7$. The Kodaira dimension of $\mathcal{A}_{6}$ is unknown. Much work has been devoted to the modular varieties $\mathcal{A}_{2}(1,d)$. 
Due to the work of Birkenhake, Lange, van Straten, Horrocks and Mumford, Manolache and  Schreyer, O'Grady, Gross and  Popescu 
it is known that $\mathcal{A}_{2}(1,d)$ is 
 unirational if $2\leq d\leq 11$. The interested reader may find a detailed discussion on these and related results in \cite{GP2}.
Gritsenko proved in \cite{Gri1} that $\mathcal{A}_{2}(1,d)$ is 
not 
uniruled if $d\geq 13$ and $d\neq 14,15,16,18,20,24,30, 36$.
 Sankaran proved in  \cite{Sa} that $\mathcal{A}_{2}(1,d)$ is of general 
type if $d$ is a prime number $\geq 173$. 

Much less is known about the modular varieties of 
non-principally 
polarized abelian varieties of dimension $\geq 3$. Tai proved in \cite{Ta2} that 
$\mathcal{A}_{g}(D)$ is of general type if $g\geq 16$ and $D$ is an 
arbitrary polarization type or if $g\geq 8$ for certain $D$.  The only result about unirationality of such modular varieties known to the author is due to Bardelli, Ciliberto 
and Verra who proved in \cite{BCV} that $\mathcal{A}_{4}(1,2,2,2)$ is unirational.
In the recent paper \cite{BL3} Birkenhake and Lange 
proved that $\mathcal{A}_{g}(D)\cong \mathcal{A}_{g}(\hat{D})$ where 
$D=(d_{1},d_2,\ldots,d_g)$, 
$\hat{D}=(\hat{d}_{1},\hat{d}_2,\ldots,\hat{d}_g)$ with 
$\hat{d}_{i}=\frac{d_{1}d_{g}}{d_{g-i+1}}$. As a consequence 
$\mathcal{A}_{4}(1,1,1,2)$ is unirational. Whether 
$\mathcal{A}_{4}(1,1,2,2)$ is unirational seems to be unknown.

In the present paper we consider threedimensional abelian varieties 
with polarizations of exponent 2 or 3. We prove that the following 
modular varieties are unirational: $\mathcal{A}_{3}(1,1,2)$,\; 
$\mathcal{A}_{3}(1,2,2)$,\; $\mathcal{A}_{3}(1,1,3)$ and 
$\mathcal{A}_{3}(1,3,3)$. By the result of Birkenhake and Lange 
cited above it suffices to verify the unirationality only of $
\mathcal{A}_{3}(1,1,2)
$ and $\mathcal{A}_{3}(1,1,3)$. The idea of the proof is the 
following. We observe that given a covering $\pi :X\to Y$ of degree 
$d=2$ or 3, with $g(X)=4,\; g(Y)=1$, the associated Prym variety 
$P=Ker(Nm_{\pi}:J(X)\to J(Y) )$ 
has polarization of type $(1,1,d)$. We consider a smooth elliptic 
fibration $\mathcal{Y}\to Z\subset \mathbb{P}^{1}$ with a section 
obtained from a general pencil of cubic curves in $\mathbb{P}^{2}$. We 
construct an open subset $T\subset \mathbb{P}(\mathbb{H})$ where 
$\mathbb{H}$ is a certain vector bundle over $Z$ and a family of 
simple degree $d$ coverings $ p:\mathcal{X}\to \mathcal{Y}_{T} =
\mathcal{Y}\times_{Z}T$. With this family one associates the Prym 
mapping $\varPhi:T\to \mathcal{A}_{3}(1,1,d)$. We prove that the 
differential of the Prym mapping $d\varPhi$ is generically surjective 
and this fact implies the unirationality of $\mathcal{A}_{3}(1,1,d)$. 
Both the 
construction of $T$ and the generic surjectivity of 
$d\varPhi$ are easier when $d=2$. 
The case $d=3$ requires the use of Prym varieties in a generalized 
sense, associated with triple coverings. Here we use 
 Miranda's result which relates triple coverings with rank 2 
vector bundles on the base \cite{Mi} and Atiyah's results about vector 
bundles over elliptic curves \cite{At}. Defining the Prym mapping and calculating its differential is a necessary work that we do in Section~\ref{s4} and Section~\ref{s2} for families of coverings of arbitrary degree. Using the results of the present paper and following the same pattern we prove in \cite{K2} the unirationality of $\mathcal{A}_{3}(1,1,4)$. The modular variety $\mathcal{A}_{3}(1,1,5)$ may be studied in a similar manner which we intend to address elsewhere.

Here is an outline of the content of the paper by sections. Section~\ref{s1} 
contains two lemmas which connect coverings of elliptic curves 
with abelian varieties with polarization of type $(1,\ldots,1,d)$. 
In Section~\ref{s3} we study the Hurwitz space 
$\mathcal{H}_{d,n}(Y)$ which parameterizes degree $d$ coverings $\pi 
:X\to Y$ of an elliptic curve $Y$ simply branched in $n\geq 2$ points. 
The direct summand $E^{\vee}$ in the decomposition 
$\pi_{*}\mathcal{O}_{X}\cong \mathcal{O}_{Y}\oplus E^{\vee}$ is called 
the Tschirnhausen module of the covering. Given $A\in Pic^{n/2}Y$ we 
denote by $\mathcal{H}_{d,A}(Y)$ the subset of $\mathcal{H}_{d,n}(Y)$ 
parameterizing coverings with $\det E\cong A$.
We define and give some simple 
properties of $\mathcal{H}_{d,n}(Y)$ and $\mathcal{H}_{d,A}(Y)$ in (\ref{s3.41}) -- (\ref{s3.41d}).
In (\ref{s3.41e}) we focus on triple coverings of elliptic 
curves and make a dimension count of the number of parameters on which 
depend triple 
coverings with a Tschirnhausen module of a given type. In order to 
make a conclusion about the type of the Tschirnhausen module of a general 
triple covering in Proposition~\ref{s3.49} we need the technical result 
of Lemma~\ref{s3.45bis}. We then prove in 
Theorem~\ref{s3.53} that $\mathcal{H}_{3,A}(Y)$ is rational if $\deg 
A$ is odd and is unirational if $\deg A$ is even. This result is of 
independent interest and is analogous to the classically known
unirationality of the
moduli space of trigonal curves.    The way we prove the unirationality of 
$\mathcal{H}_{3,n}(Y)$ suggests the construction in Proposition~\ref{s3.57}
of the family of coverings 
$p:\mathcal{X}\to \mathcal{Y}_{T}$ over a rational base $T$ of 
dimension 6 which we discussed in the previous paragraph.

Large part of Section~\ref{s4} is an overview of polarized Hodge structures 
of weight one and their variations. Some simple facts 
we need are usually included as particular cases of more general 
theorems and also usually only unimodular polarizations are considered. It 
seems to us appropriate to include in the paper some of the material we 
use. With every covering \linebreak
$\pi :X\to Y$ of smooth projective curves one 
can associate two dual abelian varieties, the Prym variety 
$P=Ker(Nm_{\pi}:J(X)\to 
J(Y))^{0}$, and its dual  $\Hat{P} = Pic^{0}X/\pi^{*}Pic^{0}Y$ polarized 
naturally by dual 
polarizations. Given a family of coverings over a smooth base $T$
one obtains respectively two morphisms $\varPhi :T\to \mathcal{A}(D)$ 
and $\hat{\varPhi}:T\to \mathcal{A}(\hat{D})$ into moduli spaces of 
abelian varieties. The morphisms $\varPhi$ and $\hat{\varPhi}$ are 
constructed by means of variations of Hodge structures of weight one. 
In Proposition~\ref{s4.75f} we give a multiplicative formula for the 
differential of $\varPhi$. 

 Section~\ref{s2} is devoted to the local study of the Prym mapping. 
In the first part of the section we consider an arbitrary simple branched covering of 
smooth, projective curves $\pi :X\to Y$ with $g(Y)\geq 1$ and its 
minimal versal 
deformation over the base $N\times H$. Here $N$ is the base of a 
minimal versal deformation of $Y$ and $H$ is a product of small disks 
centered at the branch points of $\pi$. Considering the spaces of 
holomorphic differentials with trace 0 one obtains a polarized 
variation of Hodge structures of weight 1 and a corresponding period 
mapping $\tilde{\varPhi}:N\times H\to \mathbb{D}$ where $\mathbb{D}$ 
is a period domain biholomorphically equivalent to a Siegel upper 
half space. In (\ref{s2.0}) -- (\ref{s2.9}) we work out various details 
necessary for obtaining in Proposition~\ref{s2.9} a formula for the 
differential $d\tilde{\varPhi}$. 
In the remaining part of the section 
we  restrict ourselves to the 
case of elliptic $Y$. In (\ref{s2.11}) -- (\ref{s2.18}) we obtain a 
geometric criterion in terms of the cover $X$ in order  that 
the kernel of the differential $d\tilde{\varPhi}(s_0)$, 
evaluated at the reference point corresponding to $\pi :X\to Y$, has 
dimension 1 (the minimal possible).
The check of this criterion for 
double coverings of elliptic curves is easy and is done in 
Proposition~\ref{s2.18a}. The remaining part of the section is devoted 
to the proof that the criterion is valid for general triple covers 
of genus 4, a fact needed for the proof of the unirationality of 
$\mathcal{A}_{3}(1,1,3)$. We first verify it for a union of 
two elliptic curves which intersect transversally in 3 points, and then 
obtain the result for general triple covers of genus 4 by smoothing.

In Section~\ref{s5} we prove our main results. In Theorem~\ref{s4.76} 
we give an alternative proof of a result of Birkenhake and Lange that 
$\mathcal{A}_{2}(1,2)$ and $\mathcal{A}_{2}(1,3)$ are unirational 
\cite{BL1}. Using the same argument and applying the results of 
Gritsenko and Sankaran cited above we prove in Theorem~\ref{s4.77a} 
that if  $d\geq 13$ and $d\neq 14, 15, 16, 18, 20, 24, 30, 36$ and if $A\in 
Pic^{2}Y$ is a fixed 
invertible sheaf then every connected component of the Hurwitz space 
$\mathcal{H}_{d,A}(Y)$ corresponding to coverings for which $Ker\; 
(\pi^*:J(Y)\to J(X))\; =\; 0$ is 
not uniruled. Moreover every connected component of $\mathcal{H}_{d,A}(Y)$ is of general type  if $d$ is prime and $d\geq 173$. 
Finally 
in Theorem~\ref{s4.77} we prove that $\mathcal{A}_{3}(1,1,2)$ and 
$\mathcal{A}_{3}(1,1,3)$ are unirational.

An argument of Section~\ref{s2} uses that degenerating  coverings of 
smooth projective curves into a covering of reduced curves the limit 
of the trace mapping of holomorphic differentials equals the trace 
mapping of regular differentials. As discussed in \cite{Li}~p.7
such a statement does not seem to follow  from relative 
duality. We give a proof of the 
statement we need in  Appendix A.

In Appendix B we prove the openness of the stability, the 
semistability and the regular polystability conditions for families of vector bundles over a family 
of elliptic curves. These are well-known facts that we use but for which we could not find a reference.

\par
\medskip
\noindent {\bf Notation and conventions.} 
We use the term morphism only in 
the category of schemes. When working with complex analytic spaces we use the term holomorphic mapping. Unless otherwise specified
 we make distinction between locally free sheaves and 
vector bundles and we denote differently their projectivizations. If $E$ is a 
locally free sheaf of $Y$ and if $\mathbb{E}$ is the corresponding vector 
bundle, i.e. $E\cong \mathcal{O}_Y(\mathbb{E})$, then $\mathbf{P}(E):=
\mathbf{Proj}(S(E)) \cong \mathbb{P}(\mathbb{E}^{\vee})$. A morphism (or 
holomorphic mapping) $\pi : X\to Y$ is called covering if it is finite, 
surjective and flat. If $X$ and $Y$ are smooth, then finiteness and 
surjectivity imply flatness (see e.g. \cite{Mat}~p.179 and \cite{Fi}~p.158). 
A covering of irreducible projective curves $\pi : X\to Y$ of degree 
$d$ is called simple if $X$ and $Y$ are   
smooth and for each $y\in Y$ one has $d-1\leq \#\; \pi^{-1}(y) \leq d$. 
All schemes are assumed separated of finite type over the algebraically closed base field. 
Unless otherwise specified, or clear from the context, curve means 
integral scheme of dimension one. Unless 
otherwise specified we assume the base field $k=\mathbb{C}$.

\par
\medskip
\noindent {\bf Acknowledgments.}
The author is grateful to R. Donagi, D. Markushevich and Y. Zarkhin for 
useful discussions and to M. Cornalba for sending a preliminary 
version of the book \cite{ACGH}~Vol.II.

\bigskip
\noindent
\noindent {\bf Contents.}
1. Prym varieties of coverings of elliptic curves,\;
2. Hurwitz spaces of triple coverings of elliptic curves,\;
3. Families of coverings and variations of Hodge structures,\;
4. Local study of the Prym mapping,\;
5. Unirationality results,\;
Appendix A. Traces of differential forms, \;
Appendix B. Openness conditions for families of vector bundles over families of elliptic curves.

\section{
Prym varieties of coverings of elliptic curves}\label{s1}

Let $\pi :X\to Y$ be a  covering of smooth, projective curves of 
degree $d\geq 2$, suppose $g(Y)\geq 1$. Let 
$P= Ker(Nm_{\pi}:J(X)\to J(Y))^{0} $ be the Prym variety of the covering. 
Let $\Theta$ be the canonical polarization of $J(X)$ and let 
$\Theta_{P}$ be its restriction on $P$. We give a proof of the 
following fact stated in \cite{BNR}~Remark 2.7.

\begin{lem}\label{s1.3}
The following three conditions are equivalent: 
$\pi_{*}:H_1(X,\mathbb{Z})\to H_1(Y,\mathbb{Z})$ is surjective; 
$\pi^{*}:J(Y)\to J(X)$ is injective; 
$Ker(Nm_{\pi })$ is connected. Suppose these conditions hold and let 
$P=Ker(Nm_{\pi })$. Then the polarization $\Theta_{P}$ is of type 
$(1,\ldots,1,d,\ldots,d)$ where the $d$'s are repeated $g(Y)$ times. 
\end{lem}
{\bf Proof.}
To prove the first statement we may look at these abelian varieties as real 
tori:\linebreak $J(X)\cong H_{1}(X,\mathbb{R})/H_1(X,\mathbb{Z}), \quad 
Pic^{0}X\cong 
H^{1}(X,\mathbb{R})/H^1(X,\mathbb{Z})$ and similar isomorphisms hold for 
$Y$. Then $Ker(Nm_{\pi 
})/Ker(Nm_{\pi })^{0} \cong H_1(Y,\mathbb{Z})/\pi_{*}H_1(X,\mathbb{Z})$ 
The dual finite group $\subset H^1(Y,\mathbb{R})/H^1(Y,\mathbb{Z})$ 
equals $Ker(\pi ^{*}:Pic^{0}Y\to Pic^{0}X)$. Let us prove now the second 
statement. Let 
$J(X)=V/\Lambda$. Let $B=\pi^{*}J(Y)$. We have $J(X)=B+P$. Let 
$V=V_{B}\oplus V_{P}$ be the respective decomposition. Let 
$\Lambda_{B}=V_{B}\cap \Lambda, \Lambda_{P}=V_{P}\cap \Lambda$. If $E$ 
is the Riemann form of the polarization $\Theta$, 
$E|_{\Lambda}=(\cdot,\cdot)_{X}$, let $E_{B}$ and $E_{P}$ be the 
restrictions on $V_{B}$ and $V_{P}$ respectively. The sublattices 
$\Lambda_{B}$ and $\Lambda_{P}$ are primitive and it is a standard 
fact that $\Lambda^{*}_{B}/\Lambda_{B}\cong 
\Lambda^{*}_{P}/\Lambda_{P}$ where the duals are taken with respect to 
$E_{B}, E_{P}$. By hypothesis $\pi^{*}:J(Y)\to B$ is an isomorphism 
and by the projection formula $(\pi^{*}\lambda,\pi^{*}\mu)_{X} = 
d(\lambda,\mu)_{X}$. Thus $E_{B}|_{\Lambda_{B}}$ is $d$ times a 
unimodular form, so 
$\Lambda^{*}_{B}/\Lambda_{B}\cong 
\Lambda^{*}_{P}/\Lambda_{P}\cong (\mathbb{Z}/d\mathbb{Z})^{2g(Y)}$.
\hfill $\Box$

\begin{lem}\label{s1.1} \label{s1.3a}
Let $\pi :X\to Y$ be a covering of smooth, projective curves of degree $d$, 
let $g(Y)=1$ and let
$g(X)\geq 2$. Let $P= Ker(Nm_{\pi}:J(X)\to J(Y))^{0} $ and $\Theta 
_{P}=\Theta |_{P}$ be as above. Let 
$d_{2}=|H_1(Y,\mathbb{Z}):\pi_{*}H_1(X,\mathbb{Z})| = |Ker\; \pi^*:J(Y)\to 
J(X)|$.  
Then $\pi$ may 
be decomposed as 
$X\overset{\pi_{1}}{\to}\tilde{Y}\overset{\pi_{2}}{\to}Y$ where 
$\pi_{2}$ is an isogeny of degree $d_{2}$ and $\deg \pi_{1} = d_{1} =
\frac{d}{d_{2}}$.  
The type of the polarization $\Theta_{P}$ is   $(1,\ldots ,1,d_{1})$. In 
particular if $d$ is prime then $Ker(Nm_{\pi}:J(X)\to J(Y))$ is connected and 
the type of $\Theta_{P}$ is   $(1,\ldots ,1,d)$.
\end{lem}
{\bf Proof.}
The decomposition $X\overset{\pi_{1}}{\to}\tilde{Y}\overset{\pi_{2}}{\to}Y$ is 
clear and furthermore $\pi_{1*}:H_1(X,\mathbb{Z})\to 
H_1(\tilde{Y},\mathbb{Z})$ is surjective. The statements of the lemma follow 
thus from Lemma~\ref{s1.3}
\hfill $\Box$

\bigskip
\noindent
Let $Y$ be an elliptic curve. Let $X$ be a smooth, irreducible, 
projective curve  which is  a cover $\pi : X\to Y$ of degree $d$ 
simply branched at $B\subset Y$. By Hurwitz' formula 
$g=g(X)=\frac{\# B}{2}+1$. Let $Nm_{\pi}:J(X)\to J(Y)$ be the induced 
map of the Jacobians. Then $P=Ker(Nm_{\pi})^{0}$ is an abelian variety of 
dimension $g-1=\frac{\# B}{2}$.A composition of $\pi$ with an automorphism of 
$Y$ does not change 
$P\subset J(X)$. Counting parameters we see that there are two cases 
in which one might obtain a generic abelian $(g-1)$-fold with 
polarization of type $(1,\ldots ,1,d)$ by this construction.

\medskip
\noindent
{\bf Case A.}\quad $\dim P =2$. Here $\dim \mathcal{A}_{2}(1,d)=3$. 
One obtains the same number of moduli by fixing $Y$, varying 4 branch 
points of  simple $d$-sheeted coverings and subtracting 1 for the 
action of $Aut(Y)$.

\medskip
\noindent
{\bf Case B.}\quad $\dim P = 3$. Here $\dim \mathcal{A}_{3}(1,1,d)=6$. 
This number of moduli is obtained by varying $Y$, varying the 6 branch 
points of  
simple $d$-sheeted coverings and subtracting 1 for the 
action of $Aut(Y)$.

\bigskip
\noindent
The reader is referred to Corollary~\ref{s4.80} for results about the 
representation of general abelian surfaces  and general abelian threefolds 
(with appropriate polarizations) as Prym varieties of 
coverings of elliptic curves.

\section{
Hurwitz spaces of triple coverings of elliptic curves}\label{s3}

\bigskip
\noindent
Let $\pi :X\to Y$ be 
a  covering of an elliptic curve $Y$ of degree $d$. We say that 
$\pi' :X'\to Y$ is equivalent to $\pi $ if there is an isomorphism 
$f:X\to X'$ such that $\pi = \pi' \circ f$. We are mainly interested 
in simple ramified coverings. Let $R\subset X$ be the ramification 
locus and 
$B\subset Y$ be the discriminant locus bijective to $R$. We have $\# 
B=n=2e$.

\begin{lem}\label{s3.41} Suppose $d\geq2$. Let $n=2e\geq 2$ and 
let $B\subset Y$ be a subset of $n$ points. Then there exists a simple 
covering $\pi :X\to Y$ of degree $d$ branched in $B$ with monodromy group 
$S_d$.
\end{lem}
{\bf Proof.}
Let $B=\{b_{1},\ldots,b_{n}\}$ and let $y_{0}\in Y-B$. 
The fundamental group of $Y-B$ is isomorphic to 
\[
\pi_{1}(Y-B,y_{0})\; \cong \; \langle \alpha _{1},\ldots ,\alpha _{n}, 
\gamma, \delta \: |\: \alpha _{1}\cdots \alpha _{n}\gamma \delta 
\gamma^{-1}\delta^{-1} = 1\rangle
\]
where $\alpha _{i}$ is a simple loop around $b_{i}$. Consider the 
homomorphism $m:\pi_{1}(Y-B,y_{0})\to S_{d}$ defined by:
\[
m(\alpha _{1})=(12), \ldots , m(\alpha _{n})=(12),\; 
m(\gamma) = (12\ldots d),\; m(\delta) = (12\ldots d)^{-1}.
\]
Since $(12)$ and $(12\ldots d)$ generate $S_{d}$ 
 by Riemann's existence theorem\; $m$\; is the 
monodromy homomorphism of a connected simple covering $\pi 
:X\to Y$ of degree $d$ branched at $B$.
\hfill $\Box$

\begin{block}\label{s3.41b}
The set of equivalence classes of simple coverings of $Y$ of degree 
$d$ branched at $n\geq 2$ points is parameterized by the Hurwitz space 
$\mathcal{H}_{d,n}(Y)$ (see e.g. \cite{Mo}). The Hurwitz space is an 
\'{e}tale cover of $Y^{(n)}-\Delta$, where $\Delta $ is the 
codimension one subvariety consisting of nonsimple divisors of degree 
$n$. So it is smooth equidimensional of dimension $n$. We denote by 
$p:\mathcal{X}\to Y\times \mathcal{H}_{d,n}(Y)$ the 
universal family of simple $d$-sheeted coverings branched in $n$ 
points. It has the following properties:

\begin{itemize}
\item $\mathcal{X}$ is smooth;
\item $p:\mathcal{X}\to Y\times \mathcal{H}_{d,n}(Y)$ is a covering of 
degree $d$ and the composition $\pi_{2}\circ p:\mathcal{X}\to 
 \mathcal{H}_{d,n}(Y)$ is a smooth, proper morphism with connected 
fibers;
\item for every simple $d$-sheeted covering $\pi :X\to Y$ branched in 
$n$ points there is a unique $s\in \mathcal{H}_{d,n}(Y)$ such that 
$p_{s}:\mathcal{X}_{s}\to Y\times \{s\}$ is equivalent to $\pi :X\to 
Y$.
\end{itemize}
\end{block}

\noindent
The connection between triple coverings and rank 2 bundles is due to 
R. Miranda \cite{Mi}.  We recall some basic facts taken from \cite{CE}. Let 
$\pi :X\to Y$ be a covering of degree $d$ of smooth, projective 
curves. The Tschirnhausen module of the covering is the quotient sheaf 
$E^{\vee}$ defined by the exact sequence
\[
0\lto 
\mathcal{O}_{Y}\overset{\pi^{\#}}{\lto}\pi_{*}\mathcal{O}_{X}\lto 
E^{\vee}\lto 0.
\]
One has $E^{\vee}\cong Ker(Tr_{\pi}:\pi_{*}\mathcal{O}_{X}\to 
\mathcal{O}_{Y})$ and this is a locally free sheaf of rank $d-1$. There 
is a canonical embedding $i:X\to \mathbf{P}(E)$ and 
$i^{*}\mathcal{O}_{\mathbf{P}(E)}(1)\cong \omega_{X/Y}\cong 
\omega_{X}\otimes (\pi^{*}\omega_{Y})^{-1}$ (see \cite{CE}~p.448).

\begin{lem}\label{s3.41c}
Let $R$ be the ramification divisor of $\pi :X\to Y$. Then $\deg R 
=2\deg E$.
\end{lem}
{\bf Proof.}
One has $\chi (\mathcal{O}_{X})=\chi 
(\pi_{*}\mathcal{O}_{X})=\chi(\mathcal{O}_{Y})\oplus \chi(E^{\vee})$. 
By Riemann-Roch $\chi(E^{\vee})=\deg E^{\vee}+(d-1)(1-g(Y))$. Therefore 
\[
1-g(X)\ =\ \deg E^{\vee} + d(1-g(Y)).
\]
The equality $2\deg E = \deg R$ follows thus from Hurwitz' formula.
\hfill $\Box$

\begin{block}\label{s3.51}
Let $\mathcal{X}\to Y\times \mathcal{H}_{d,n}(Y)$ be the universal 
family of simple coverings of degree $d\geq 2$ branched in $n=2e$ 
points. Let 
$\mathcal{E}^{\vee}$ be the corresponding Tschirnhausen module. Let 
$\mathcal{A}=\det \mathcal{E}$. By Lemma~\ref{s3.41c} for every $z\in 
\mathcal{H}_{d,n}(Y)$ one has 
$\deg(\mathcal{A}_{z}) = e$. We obtain a morphism 
$h:\mathcal{H}_{d,n}(Y)\to Pic^{e}Y,\; h(z)=
\mathcal{A}_{z}$. 
 Let $A\in Pic^{e}Y$. We let 
$\mathcal{H}_{d,A}(Y)=h^{-1}(A)$.
\end{block}

\begin{lem}\label{s3.51a} 
Let $d$ and $e$ be integers such that $d\geq 2,\; e\geq 1 $. Let $n=2e$ and 
let $A\in 
Pic^{e}Y$. 
 The following properties hold.
\renewcommand{\theenumi}{\roman{enumi}}
\begin{enumerate}
\item 
The morphism $h:\mathcal{H}_{d,n}(Y)\to Pic^eY$ is surjective.
\item
If $A'\in Pic^{e}Y$ then $\mathcal{H}_{d,A}(Y)\cong 
\mathcal{H}_{d,A'}(Y)$.
\item
$\mathcal{H}_{d,A}(Y)$ is smooth, equidimensional of dimension $n-1$.
\end{enumerate}
\end{lem}
{\bf Proof.}
(ii).\quad Let $\alpha \in Pic^{0}Y$ and let $t_{\alpha }:Y\to Y$ be 
the corresponding translation. If $E^{\vee}$ is the Tschirnhausen 
module of $\pi :X\to Y$ then $t_{\alpha }^{*}E^{\vee}$ is the 
Tschirnhausen module of $t_{-\alpha }\circ \pi :X\to Y$. So $A=\det E$ and if 
$A'=t_{\alpha }^{*}A$ then the mapping $[\pi :X\to Y]\mapsto [
t_{-\alpha }\circ \pi :
X\to Y]$ yields an isomorphism between $\mathcal{H}_{d,A}(Y)$ and 
$\mathcal{H}_{d,A'}(Y)$.

(i) and (iii).\quad Given an element $[\pi:X\to Y]\in \mathcal{H}_{d,n}(Y)$ 
the 
translation $[t_{-\alpha}\circ \pi :X\to Y]$ with $\alpha \in Pic^{0}Y$ 
belongs to the same connected component of $\mathcal{H}_{d,n}(Y)$. Thus 
the morphism $h:\mathcal{H}_{d,n}(Y)\to Pic^{e}Y$ has the property 
that its restriction on every connected component of 
$\mathcal{H}_{d,n}(Y)$ is surjective, so Part (i) holds. Hence every 
sufficiently  
general fiber $\mathcal{H}_{d,A}(Y)$ is equidimensional of dimension 
$n-1$ and is furthermore smooth \cite{Ha}~III.10.7. Using (ii) we 
conclude that every fiber $\mathcal{H}_{d,A}(Y)$ has these properties.
\hfill $\Box$

\begin{block}\label{s3.41d}
If $\pi :X\to Y$ is a triple covering of smooth projective curves  
with Tschirnhausen module isomorphic to $E^{\vee}$ then
$\mathbf{P}(E)$ is a 
ruled surface, $\varphi :\mathbf{P}(E)\to Y$, and $i(X)$ is a divisor 
of the linear system $|\mathcal{O}_{\mathbf{P}(E)}(3)\otimes 
\varphi^{*}(\det E)^{-1}|$. Hence $\pi :X\to Y$ is uniquely determined 
by $E$ and by an element $\langle \eta \rangle \in 
\mathbb{P}H^{0}(Y,S^{3}E\otimes 
(\det E)^{-1})$. Two equivalent coverings $\pi :X\to Y$ and $\pi' :X'\to 
Y$ yield isomorphic Tschirnhausen modules such that $\langle\eta\rangle$ is 
transformed into $\langle\eta'\rangle$ under the isomorphism. Conversely, 
given 
a rank 2 locally free sheaf $E$ on $Y$ with associated ruled surface 
$\varphi :\mathbf{P}(E)\to Y$ any nonsingular, irreducible divisor 
$X\in 
|\mathcal{O}_{\mathbf{P}(E)}(3)\otimes 
\varphi^{*}(\det E)^{-1}|$ determines a triple covering $\pi :X\to Y$ 
with Tschirnhausen module isomorphic to $E^{\vee}$. The group 
$PGL_{Y}(E)$ acts faithfully on the set of reduced divisors of 
$|\mathcal{O}_{\mathbf{P}(E)}(3)\otimes 
\varphi^{*}(\det E)^{-1}|$ and the orbit $
PGL_{Y}(E)\cdot X$
 corresponds to the equivalence class $[X\to Y]$.
\end{block}

\begin{block}\label{s3.41e}
We now consider simple triple branched coverings $\pi :X\to Y$ where $Y$ is 
elliptic. The 
number of branch points $n=\# B=\# R=2\deg E, \: n\geq 2$. Let 
$\mathcal{H}_{3,n}(Y)$  be the corresponding Hurwitz space. We wish to 
determine the types of the Tschirnhausen modules of general coverings. 
We need to bound above the number of moduli of $[X\to Y]\in 
\mathcal{H}_{3,n}(Y)$ with certain types of Tschirnhausen modules 
(Cases 1 -- 5 considered below). We do not need to address the 
question whether such coverings with smooth, irreducible $X$ actually 
exist.
Let us first consider the case of decomposable Tschirnhausen modules. 
Let $\deg E=e,\: n=2e,\: E=Lu\oplus Mv$ where $u\in \Gamma(Y,E\otimes 
L^{-1}),\: 
v\in \Gamma(Y,E\otimes M^{-1})$. Let $a=\deg L,\: b=\deg M, \: e=a+b$. 
we may assume $a\leq b$. We have
\begin{equation}\label{es3.42}
S^3E\otimes(\det E)^{-1} = 
L^{2}M^{-1}u^{3}+Lu^{2}v+Muv^{2}+M^{2}L^{-1}v^{3}.
\end{equation}
Hence a global section $\eta \in H^0(Y,S^3E\otimes(\det E)^{-1})$ may 
be decomposed as $\eta = \alpha u^{3}+\beta u^{2}v+\gamma 
uv^{2}+\delta v^{3}$. In order that $E^{\vee}$ is the Tschirnhausen 
module of the irreducible covering $\pi :X\to Y$ determined by $\eta$ it 
is necessary that $\alpha \ne 0$ and $\delta \ne 0$ 
(see \cite{Mi}~p.1145 or \cite{CDC}~p.266). Hence $L^{2}\geq M$ and $M^{2}\geq 
L$. If $a\leq 0$ then $b\leq 2a\leq 0$ which is absurd. Hence $a\geq 
1,\, b\geq 1,\, 2a\geq b,\,$ and $2b>a$ since $b\geq a$. We see that 
if $n=2$ and if $X$ is irreducible then the Tschirnhausen module is 
indecomposable. The splitting \eqref{es3.42} yields
\[
h^0(Y,S^3E\otimes(\det E)^{-1})\ =\ (2a-b+\epsilon) +a+b+2b-a\ =\ 2e + 
\epsilon
\]
where $\epsilon =1$ if $L^2\cong M$ and
$\epsilon =0$ otherwise. We have
\[
\dim Aut_{Y}(E) = h^0(Y,End(E)) = 2+h^0(Y,LM^{-1})+h^0(Y,ML^{-1})
\]
{\bf Case 1.} \quad $a<b, \, 2a>b$. Here $L$ and $M$ may vary 
independently, so
\[
\# \, moduli \leq 2 +h^0(Y,S^3E\otimes(\det E)^{-1})-h^0(Y,End(E)) = 2 
+ 2e -2 - (b-a) = n - (b-a)
\]
{\bf Case 2.}\quad $a=b, \, L\ncong M$. This case is possible only if 
$n \equiv 0(mod\: 4)$ since $n=2(a+b)$. Here $\# \, moduli \leq 
2+2e-2=n$.

\medskip
\noindent
{\bf Case 3.}\quad $a=b,\, L\cong M$. Here
$\# \, moduli \leq 1+2e-4=n-3$.

\medskip
\noindent
{\bf Case 4.}\quad $2a=b$. Here one has two subcases: $L^2\ncong M$ and 
$L^2\cong M$. In both subcases
$\# \, moduli \leq 
2+2e-2-(b-a)=n-(b-a)$.

\smallskip
\noindent
Recall from \cite{At}~p.432 that on an elliptic curve for every 
$r\geq 1$ up to isomorphism there is a unique indecomposable locally 
free sheaf $F_r$ of rank $r$ and degree 0 with $h^0(Y,F_r)\ne 0$. 
Furthermore any indecomposable locally free sheaf $E$ of degree $e$ and 
rank $r$ such that $r|e$ is isomorphic to $L\otimes F$ where $L$ is an
invertible sheaf of degree $\frac{e}{r}$.

\smallskip
\noindent
{\bf Case 5.}\emph{The Tschirnhausen module of $\pi :X\to Y$ is 
indecomposable of even degree}. Let $n=2e$, $e$ is even. Let us fix an 
indecomposable $E$ with $\deg E=e$. Let $E\cong L\otimes F_2$. We have 
by \cite{At}~p.438 that 
\[
S^3E\otimes(\det E)^{-1} \cong L^{3}\otimes S^{3}F_{2}\otimes 
L^{-2}\cong L\otimes F_{4}.
\]
Hence by \cite{At}~p.430
\[
h^0(S^3E\otimes(\det E)^{-1})=h^0(L\otimes F_{4})=\deg(L\otimes 
F_{4})=4\frac{e}{2}=n.
\]
One has $End_{Y}(E)\cong F_{2}\otimes F_{2}^{\vee}\cong F_{2}\otimes 
F_{2}\cong F_{1}\oplus F_{3}$ by \cite{At}~pp.433,437, hence 
$h^0(End_{Y}(E))=2$. Varying $E$ (i.e. $L\in Pic^{e/2}Y$) we obtain 
\(
\# moduli\ \leq \ 1 + n - 2 = n-1.
\)
\end{block}

\medskip
\noindent
We need a technical result related to \cite{CE}~Theorem~3.6. Although 
we need it only in the case of families of elliptic curves we state 
and prove it for arbitrary dimensions and arbitrary algebraically 
closed base field $k$ of characteristic 0. Let $Y$ be a smooth integral 
scheme over $k$. We 
recall from \cite{CE}~Definition~3.3 that given a rank 2 locally free 
sheaf $E$ on $Y$ an element $\eta \in H^0(Y,S^3E\otimes(\det E)^{-1})$ 
is called of right codimension in every $y\in Y$ if 
$\eta(y)\in(S^3E\otimes(\det E)^{-1})\otimes _{\mathcal{O}_{Y}}k(y)$ 
is nonzero for every $y\in Y$. Every such $\eta$ determines a 
Gorenstein triple covering $X_{\eta}\to Y,\: X_{\eta}\subset 
\mathbf{P}(E)$ with Tschirnhausen module isomorphic to $E^{\vee}$ 
\cite{CE}~p.449.

\begin{lem}\label{s3.45bis} 
Let $q:\mathcal{Y}\to Z$ be a smooth proper morphism with connected 
fibers, where $Z$ is smooth.  Let $\mathcal{E}$ be a 
locally free sheaf of rank 2 on $\mathcal{Y}$ such that 
$h^0(\mathcal{Y}_{z},S^3\mathcal{E}_{z}\otimes(\det 
\mathcal{E}_{z})^{-1})$ is independent of $z\in Z$ and is $\ne 0$. 
Consider the locally free sheaf $\mathcal{H}=q_{*}(S^3\mathcal{E}\otimes(\det 
\mathcal{E})^{-1})$ on $Z$. Let $f:\mathbb{H}\to Z$ be the associated 
vector bundle with fibers 
$\mathbb{H}_z=H^0(\mathcal{Y}_{z},S^3\mathcal{E}_{z}\otimes(\det 
\mathcal{E}_{z})^{-1})$. Then the subset $\mathbb{H}_{0}\subset \mathbb{H}$ 
consisting of $\eta$ which satisfy the following three conditions is 
Zariski open in $\mathbb{H}$.
\renewcommand{\theenumi}{\alph{enumi}}
\begin{enumerate}
\item
If $f(\eta)=z$ then $\eta$ is of right codimension for every $y\in 
\mathcal{Y}_{z}$.
\item
Assuming (a), if $\pi_{\eta} :X_{\eta}\to \mathcal{Y}_{z}, \: X_{\eta}\subset 
\mathbf{P}(\mathcal{E}_{z})$ is the triple covering determined by 
$\eta$, then $X_{\eta }$ is smooth and irreducible.
\item
Assuming (a) and (b) the discriminant scheme of $\pi_{\eta}:X_{\eta}\to 
\mathcal{Y}_{z}$ is a smooth subscheme of $\mathcal{Y}_{z}$.
\end{enumerate}
Suppose $\mathbb{H}_{0}\ne \emptyset$. Consider the base change 
$\mathcal{Y}'=\mathcal{Y}\times _{Z}\mathbb{H}_{0}$ and let 
$\mathcal{E}'=\pi_{1}^{*}\mathcal{E}$. Then there is a smooth 
subscheme $\mathcal{X}\subset \mathbf{P}(\mathcal{E}')$ and a triple 
covering $p:\mathcal{X}\to \mathcal{Y}'$ such that for every $\eta \in 
\mathbb{H}_{0}$ with $f(\eta)=z$ the fiber
$\mathcal{X}_{\eta}\to \mathcal{Y}'_{\eta}$ is equivalent to 
$X_{\eta}\to \mathcal{Y}_{z}$.
\end{lem}
{\bf Proof.}
The statement is local with respect to $Z$ so we may assume $Z$ is irreducible.
If $\mathbb{H}_{0}$ is empty there is nothing to prove. Suppose 
$\mathbb{H}_{0}\ne \emptyset$.

\emph{Step 1.}\quad Let $\mathbb{H}'$ be the set of $\eta \in \mathbb{H}$ 
for which (a) holds. We claim $\mathbb{H}'$ is Zariski open in 
$\mathbb{H}$. Consider the incidence correspondence $\Gamma \subset 
\mathbf{P}(\mathcal{E})\times _{Z}\mathbb{H}$ defined as follows. 
\[
\Gamma = \{(x,\eta)|\eta(x)=0\; \text{where}\; x\in 
\mathbf{P}(\mathcal{E})_{y}, \eta \in \mathbb{H}_{z}, y\in 
\mathcal{Y}_{z}\}.
\]
Consider the projection $\varepsilon :\Gamma \to \mathcal{Y}\times_Z
\mathbb{H},\; \varepsilon (x,\eta)=(y,\eta)
$. An element $\eta \in \mathbb{H}_{z}$ fails to be of right codimension 
in $y\in \mathcal{Y}_{z}$ if and only if $\eta(x)=0$ for $\forall x\in 
\mathbf{P}(\mathcal{E})_{y}$. Equivalently $(y,\eta)\in \Sigma$ where 
$\Sigma \subset \mathcal{Y}\times _{Z}\mathbb{H}$ is the subset of 
points for which $\dim \varepsilon ^{-1}(y,\eta)\geq 1$. Hence $\Sigma$ is 
closed in $\mathcal{Y}\times _{Z}\mathbb{H}$. Since properness is 
preserved under base change the projection of $\Sigma$ in $\mathbb{H}$ 
is closed. This proves $\mathbb{H}'$ is open in $\mathbb{H}$.

\emph{Step 2.}\quad  Let $\mathbb{H}''\subset \mathbb{H}'$ be the set 
of $\eta\in \mathbb{H}$ for which both (a) and (b) hold. We claim 
$\mathbb{H}''$ is open in $\mathbb{H}'$. There is a commutative 
diagram
\[
\begin{diagram}
          &                         & S^3\mathbb{E}\otimes(\det \mathbb{E})^{-1}\\
          & \ruTo^{N'} & \dTo               \\
\mathcal{Y}\times _{Z}\mathbb{H} & \rTo^{\pi_{1}} &\mathcal{Y}\\
\end{diagram}
\]
where $N'(y,\eta)=\eta(y)$. Let 
$\mathbb{E}_{\mathbb{H}}=\pi_{1}^{*}\mathbb{E}$. Since 
$\pi_{1}^{*}(S^3\mathbb{E}\otimes(\det \mathbb{E})^{-1})\cong 
S^3\mathbb{E}_{\mathbb{H}}\otimes(\det \mathbb{E}_{\mathbb{H}})^{-1}$
we obtain a tautological section $N \in H^0(
\mathcal{Y}\times _{Z}\mathbb{H}
,S^3\mathbb{E}_{\mathbb{H}}\otimes(\det 
\mathbb{E}_{\mathbb{H}})^{-1})$. Restricting to 
$\mathcal{Y}\times _{Z}\mathbb{H}'$
 we obtain a section which is of right codimension for every 
$(y,\eta)\in \mathcal{Y}\times _{Z}\mathbb{H}'$. By 
\cite{CE}~Theorem~3.4 one obtains a closed subscheme $\mathcal{X}'\subset 
\mathbf{P}(\mathcal{E}_{\mathbb{H}'})\quad 
(=\mathbb{P}(\mathbb{E}_{\mathbb{H}'}^{\vee}))$ and a Gorenstein 
triple covering $p:\mathcal{X}'\to \mathcal{Y}\times _{Z}\mathbb{H}'$ 
whose Tschirnhausen module is isomorphic to 
$\mathcal{E}_{\mathbb{H}'}^{\vee}$. Let $f:\mathcal{X}'\to \mathbb{H}'$ 
be the projection morphism. Let $f(x)=\eta$. The point $x$ is 
nonsingular in the fiber $\mathcal{X}'_{\eta}$ if and only if 

(i)\; 
$x$ is a nonsingular point of $\mathcal{X}'$,\qquad
(ii)
$T_{x}f:T_{x}\mathcal{X}'\to T_{\eta}\mathbb{H}'$ is surjective

\noindent
(see e.g. \cite{AK}~pp.131,145). Since $f$ is proper $f(\Sing 
\mathcal{X}')$ is closed in $\mathbb{H}'$. Let 
$\mathbb{H}'_{1}=\mathbb{H}'-f(\Sing \mathcal{X}')$. The set of $x\in 
f^{-1}(\mathbb{H}'_{1})$ where the rank of $T_{x}f$ is not maximal is 
closed in $f^{-1}(\mathbb{H}'_{1})$, so again by properness its image 
in $\mathbb{H}'_{1}$ is closed. The complement of the latter in 
$\mathbb{H}'_{1}$ is  the set $U$ consisting of $\eta\in \mathbb{H}'$ for 
which $\mathcal{X}'_{\eta}$ is smooth. Considering the Stein factorization 
 we see that all fibers of $f^{-1}(U)\to U$  have the same number of irreducible components. Since by assumption 
$\mathbb{H}_{0}\ne \emptyset$ and $Z$ is irreducible we conclude that all fibers $\mathcal{X}'_{\eta},\quad \eta \in U$ are irreducible. So the set 
$\mathbb{H}''$ consisting of $\eta \in \mathbb{H}$ for which 
both (a)  and  (b) hold equals $U$ and is therefore open in $\mathbb{H}'$. 

\emph{Step 3.}\quad We claim $\mathbb{H}_{0}$ is open in 
$\mathbb{H}''$. Let $\mathcal{X}''\to \mathcal{Y}\times 
_{Z}\mathbb{H}''$ be the restriction of $
\mathcal{X}'\to \mathcal{Y}\times 
_{Z}\mathbb{H}'
$. Let $\mathcal{B}\subset 
\mathcal{Y}\times 
_{Z}\mathbb{H}''
$ be the discriminant subscheme (see e.g. \cite{AK}~pp.123-124). We 
apply to the projection $\mathcal{B}\to \mathbb{H}''$ the same 
argument as that in Step 2 to conclude that $\mathbb{H}_{0}$ is open 
in $\mathbb{H}''$. 

This proves $\mathbb{H}_{0}$ is open in $\mathbb{H}$. If 
$\mathbb{H}_{0}\ne \emptyset$ let $\mathcal{X}\to \mathcal{Y}'=
\mathcal{Y}\times 
_{Z}\mathbb{H}_{0}
$ be the restriction of $
\mathcal{X}'\to \mathcal{Y}\times 
_{Z}\mathbb{H}'
$ from Step 2. The last claim of the lemma follows from the 
functoriality of Miranda's construction.
\hfill $\Box$

\begin{pro}\label{s3.49}
Let $n=2e\geq 2$. There is a Zariski open dense subset of the Hurwitz 
space $\mathcal{U}\subset \mathcal{H}_{3,n}(Y)$ such that for every 
$[X\to Y]\in \mathcal{U}$ one has 
\begin{enumerate}
\item
if $e\equiv 1(mod\; 2)$ the Tschirnhausen module $E^{\vee}$ is 
indecomposable of degree $-e$.
\item
if $e\equiv 0(mod\; 2)$ the Tschirnhausen module $E^{\vee}$ is 
isomorphic to $L^{-1}\oplus M^{-1}$ where $\deg L=\deg M = 
\frac{e}{2}$ and $L\ncong M$.
\end{enumerate}
\end{pro}
{\bf Proof.}
We need to prove that in each of the cases 1, 3 and 5 from 
(\ref{s3.41e}) either there is no triple covering $\pi :X\to Y$ with 
Tschirnhausen module of that type or if it exists the set of 
equivalence classes of such coverings is contained in a closed subscheme of 
$\mathcal{H}_{3,n}(Y)$ of codimension $\geq 1$. We apply 
Lemma~\ref{s3.45bis} with $\mathcal{Y}=Y\times Z,\; q=\pi_{2}$ where $Z$ and 
$\mathcal{E}$ 
are constructed in the 
various cases as follows. In Case 1 we let $Z=Pic^{a}Y\times Pic^{b}Y$ 
and if $\mathcal{L}\to Y\times Pic^{a}Y$ and $\mathcal{M}\to Y\times 
Pic^{b}Y$ are the Poincar\'{e} invertible sheaves we let 
$\mathcal{E}=\pi_{12}^{*}\mathcal{L}\oplus \pi_{13}^{*}\mathcal{M}$.
Here $\pi_{12}$ and $\pi_{13}$ are the corresponding projections of 
$Y\times Z = Y\times Pic^{a}Y\times Pic^{b}Y$.
In Case 3 we let $Z=Pic^{a}Y$ and $\mathcal{E}=\mathcal{L}\oplus 
\mathcal{L}$.  Case 4 splits into two subcases.
When $L^2\ncong M$ we repeat the construction of Case 1 letting 
$Z\subset Pic^aY \times Pic^b Y$ be the open subset 
$Z=\{(u,v)|\mathcal{L}_u^2\ncong \mathcal{M}_v\}$. In the subcase $L^2\cong M$
 we let $Z=Pic^{a}Y$ and $
\mathcal{E}=\mathcal{L}\oplus 
\mathcal{L}^{2}
$. In Case 5 we let $Z=Pic^{e/2}Y$ and let 
$\mathcal{E}=(\pi_{1}^{*}F_2)\otimes \mathcal{L}$. Suppose there is a 
triple covering of one of the considered types. Let $\mathbb{H}_{0}$ 
and $\mathcal{X}\to Y\times \mathbb{H}_{0}$ be as in 
Lemma~\ref{s3.45bis}. By the universal property of the Hurwitz space 
there is an associated morphism $f:\mathbb{H}_{0}\to 
\mathcal{H}_{3,n}(Y)$. The closure $\overline{f(\mathbb{H}_{0})}$ is a 
closed subscheme of codimension $\geq 1$ according to the calculations 
of the number of moduli in the various cases from (\ref{s3.41e}).
\hfill $\Box$

\begin{thm}\label{s3.53}
Let $Y$ be an elliptic curve, let $n$ be  a pair integer $n=2e\geq 
2$. Let $A\in Pic^{e}Y$. If $d=2$ or 3 the Hurwitz spaces 
$\mathcal{H}_{d,n}(Y)$ and $\mathcal{H}_{d,A}(Y)$ are irreducible. 
The variety 
$\mathcal{H}_{d,A}(Y)$ is rational
if 
$d=2$ or if $d=3$ and $e\equiv 1(mod\; 2)$ and is unirational 
 if $d=3$ and $e\equiv 0(mod\; 2)$.
\end{thm}
{\bf Proof.}
According to Lemma~\ref{s3.51a}\;
$h:\mathcal{H}_{d,n}(Y)\to Pic^{e}Y$\; is surjective. It has fibers 
$\mathcal{H}_{d,A}(Y),\; 
A\in 
Pic^{e}Y$
 which are isomorphic to each other. Hence it suffices to prove the 
statements for $\mathcal{H}_{d,A}(Y)$. 

If $d=2$ what is claimed is obvious, since a simple double covering 
$\pi :X\to Y$ such that $\pi_{*}\mathcal{O}_{X}\cong 
\mathcal{O}_{Y}\oplus A^{-1}$ is uniquely determined by an element 
$\langle \eta \rangle \in \mathbb{P}H^0(Y,A^{2})$ such that 
$div(\eta)$ is a simple divisor (see e.g. \cite{Wa}). If $\Delta \subset 
\mathbb{P}H^0(Y,A^{2})
$ is the closed subset of nonsimple divisors one obtains an 
isomorphism $f:\mathbb{P}H^0(Y,A^{2})-\Delta 
\overset{~}{\lto}\mathcal{H}_{2,A}(Y)$.

If $d=3$ we have two cases.

Case 1.  $e\equiv 1(mod\; 2)$. According to Atiyah's results 
\cite{At} up to isomorphism there is a unique indecomposable rank 2 
locally free sheaf $E$ on $Y$ with $\deg E\cong A$. By 
Proposition~\ref{s3.49} and Lemma~\ref{s3.51a} there exist simple 
triple coverings with Tschirnhausen module isomorphic to $E^{\vee}$. 
Applying Lemma~\ref{s3.45bis} with $Z=\{*\},\mathcal{E}=E$ one obtains a 
covering $\mathcal{X}\to Y\times \mathbb{H}_{0}$, where $\mathbb{H}_0$ 
is Zariski open nonempty subset of $\mathbb{H}=H^0(Y,S^3E\otimes(\det 
E)^{-1})$ which moreover is invariant with respect to  the action of 
$\mathbb{C}^{*}$. Using the universal property of the Hurwitz space 
$\mathcal{H}_{3,n}(Y)$ one obtains a morphism 
$f:\mathbb{P}\mathbb{H}_{0}\to \mathcal{H}_{3,A}(Y)$. This morphism is 
dominant by Proposition~\ref{s3.49}. It is injective since 
$h^0(Y,End(E))=1$ (cf. \cite{At}~Lemma~22 and (\ref{s3.41d})). Hence 
$\mathcal{H}_{3,A}(Y)$ is a rational variety.

Case 2.  $e\equiv 0(mod\; 2)$. Let $e=2a$. Let $\sigma : 
Pic^{a}Y\to Pic^{a}Y$ be the involution $L\mapsto A\otimes L^{-1}$ and 
let $\mu :Pic^{a}Y\to \mathbb{P}^{1}$ be the quotient map. Let $Z\subset 
\mathbb{P}^{1}$ be the complement of the branch locus of $\mu$. 
Consider the double covering $1\times \mu :Y\times Pic^{a}Y \to Y\times 
\mathbb{P}^{1}$. Let 
$\mathcal{L}$ be the Poincar\'{e} invertible sheaf on $Y\times 
Pic^{a}Y$ and let $\mathcal{E}=(1\times \mu)_{*}\mathcal{L}|_{Y\times 
Z}$. By construction for every $z\in Z$ one has $\mathcal{E}_{z}\cong 
\mathcal{L}_{z}\oplus (A\otimes \mathcal{L}_{z}^{-1})$. We apply 
Lemma~\ref{s3.45bis} with $\mathcal{Y}=Y\times Z,\; q=\pi_{2}$
and using Proposition~\ref{s3.49} and 
Lemma~\ref{s3.51a} we conclude that $\mathbb{H}_{0}$ is nonempty. Here 
$\mathbb{H}$ is a vector bundle over $Z$ and its Zariski open subset 
$\mathbb{H}_{0}$ is $\mathbb{C}^{*}$-invariant. The family of triple 
coverings $\mathcal{X}\to Y\times \mathbb{H}_{0}$ yields a morphism 
$f:\mathbb{P}\mathbb{H}_{0}\to \mathcal{H}_{3,A}(Y)$ which is dominant 
by Proposition~\ref{s3.49}. Hence $\mathcal{H}_{3,A}(Y)$ is 
irreducible and unirational.
\hfill $\Box$

\begin{rem}
 Graber, Harris and Starr proved in  \cite{HGS} the 
irreducibility of the space $\mathcal{H}_{d,n}^{^{S_d}}(Y)$ parameterizing 
simple coverings with monodromy group $ S_d$ for any $Y$ of positive genus 
when $n\geq 2d$. This result implies the irreducibility of $H_{3,n}(Y)$ stated in the above theorem when $n\geq 6$.
\end{rem}
The proof of the theorem ($d=3$, Case 1) shows the following result.
\begin{cor}\label{s3.55a}
Let $A\in Pic^{e}Y,\; e\equiv 1(mod\; 2),\; e\geq 1$, and let $E$ be an 
indecomposable rank 2 locally free sheaf on $Y$ with $\det E\cong A$. 
Then a Zariski open nonempty subset of the Hurwitz space 
$\mathcal{H}_{3,A}(Y)$ consists of equivalence classes of coverings 
$[X\to Y]$ with Tschirnhausen module isomorphic to $E^{\vee}$.
\end{cor}

\begin{block}\label{s3.55}
So far in this section we considered  families of triple coverings of 
a fixed elliptic curve. We now want to vary also the elliptic curve. 
We need only the case of triple coverings simply branched at 6 points and 
this is the case we work out. We consider a sufficiently general 
pencil of cubic curves in $\mathbb{P}^{2}$. Blowing up the nine base 
points and discarding the singular fibers we obtain a smooth family 
$q:\mathcal{Y}\to Z,\quad Z\subset \mathbb{P}^{1}$ of elliptic curves 
with 9 sections. Let us choose one of the sections and call it $D$. We 
construct a rank 2 locally free sheaf on $\mathcal{Y}$ as in 
\cite{Ha}~Ch.V~Ex.2.11.6. Namely, the extensions
\begin{equation}\label{es3.55}
0\lto \mathcal{O}_{\mathcal{Y}}\lto \mathcal{F}\lto 
\mathcal{O}_{\mathcal{Y}}(D)\lto 0
\end{equation}
are parameterized by $H^1(\mathcal{Y},\mathcal{O}_{\mathcal{Y}}(-D)) = 
H^0(Z,R^{1}q_{*}\mathcal{O}_{\mathcal{Y}}(-D))$ by Leray's spectral 
sequence. The sheaf $R^{1}q_{*}\mathcal{O}_{\mathcal{Y}}(-D)$ is 
locally free of rank one by Grauert's theorem. Replacing $Z$ by a 
smaller affine 
set we may assume $R^{1}q_{*}\mathcal{O}_{\mathcal{Y}}(-D)$ is 
trivial. A trivializing section yields an extension \eqref{es3.55} 
with the property that $\mathcal{F}_{z}$ is indecomposable over 
$\mathcal{Y}_{z}$ for every $z\in Z$. Let 
$\mathcal{E}=\mathcal{F}\otimes \mathcal{O}_{\mathcal{Y}}(D)$. 
By semicontinuity, replacing $Z$ by a smaller open set, we may assume 
$h^0(\mathcal{Y}_{z},S^3\mathcal{E}_{z}\otimes(\det \mathcal{E}_{z})^{-1})$ is 
independent of $z$. We may now apply Lemma~\ref{s3.45bis}. If 
$\mathcal{H}=q_{*}(S^3\mathcal{E}\otimes(\det \mathcal{E})^{-1})$, if 
$f:\mathbb{H}\to Z$ is the corresponding vector bundle and 
$\mathbb{H}_{0}\subset \mathbb{H}$ the open subset that satisfies the 
three conditions of Lemma~\ref{s3.45bis} then $\mathbb{H}_{0}\ne 
\emptyset $ by Corollary~\ref{s3.55a}. According to that lemma letting 
$\mathbb{P}(\mathbb{H}_{0})=T$
we obtain a family of triple coverings
\begin{equation}\label{es3.56a}
\begin{diagram}
\mathcal{X}& \rTo & \mathcal{Y} \\
\dTo       &      &   \dTo\\
T          & \rTo &  Z      \\
\end{diagram}
\end{equation}
Letting $\mathcal{Y}_{T}=\mathcal{Y}\times 
_{Z}T$ one obtains a family of coverings over $T$:
\begin{equation}\label{es3.56}
\begin{diagram}
\mathcal{X}&    &\rTo^p &       &\mathcal{Y}_{T} \\
           &\rdTo &     & \ldTo  &               \\
           &      &  T  &        &               \\
\end{diagram}
\end{equation}
In the case of double coverings one has an analogous and simpler 
construction. Here one lets 
$\mathcal{A}=\mathcal{O}_{\mathcal{Y}}(3D),\; 
\mathcal{H}=q_{*}\mathcal{A}^{2}$ and $T=\mathbb{P}(\mathbb{H}_{0})$. 
One obtains a family as above with $\deg p=2$.
\end{block}

\begin{pro}\label{s3.57}
Let $d=2$ or 3. The constructed family of coverings of degree $d$ has 
the following properties.
\renewcommand{\theenumi}{\alph{enumi}}
\begin{enumerate}
\item
Every sufficiently general elliptic curve is isomorphic to a fiber of 
$\mathcal{Y}\to Z$.
\item
Let $z\in Z$. The fibers $\mathcal{X}_{\eta}\to \mathcal{Y}_{z}$ with 
$\eta \in T,\: f(\eta )=z$ correspond to a Zariski open nonempty subset 
of the Hurwitz space $\mathcal{H}_{d,A}(\mathcal{Y}_{z})$ with 
$A=\mathcal{O}_{\mathcal{Y}_{z}}(3D_{z})$.
\item T is a rational variety of dimension 6.
\end{enumerate}
\end{pro}
{\bf Proof.}
Part (b) follows from Corollary~\ref{s3.55a} and Lemma~\ref{s3.45bis}. 
The other parts are clear.
\hfill $\Box$

\section{Families of coverings and variations of Hodge structures}\label{s4}
We start by recalling some well known facts and fixing notation 
(cf. \cite{Gr},\cite{GS},\cite{Ke}).
\begin{block}\label{s4.58}
A polarized Hodge structure of weight one is given by the following 
data.
\renewcommand{\theenumi}{\roman{enumi}}
\renewcommand{\theenumii}{\alph{enumi}}
\begin{enumerate}
\item
A free abelian group $M$ of rank $2g$.
\item
A complex subspace $U\subset M_{\mathbb{C}}=M\otimes
_{\mathbb{Z}}\mathbb{C}$ such that 
$M_{\mathbb{C}}=U\oplus \overline{U}$.
\item
An integer valued, skew-symmetric, nondegenerate form $Q:M\times M\to 
\mathbb{Z}$ whose $\mathbb{C}$-bilinear extension $Q_{\mathbb{C}}$ 
satisfies the Riemann relations
\begin{enumerate}
\item $Q_{\mathbb{C}}(U,U) = 0$,
\item $iQ_{\mathbb{C}}(u,\overline{u})>0$ for $\forall u\in U, u\ne 0$
\end{enumerate}
\end{enumerate}
One defines the real Weyl operator $C$ so that it has eigenvalue $i$ 
on $U$ and $-i$ on $\overline{U}$. Given such data one may define a 
polarized abelian variety $P(U)=\overline{U}/\pi^{0,1}(M)$. The 
polarization on $P(U)$ is defined by the Hermitian form $H'$ on 
$\overline{U}$ such that $Im(H')=E':=-Q_{\mathbb{R}}$. Let us check a 
part of the last statement. If $J=-C$ is the complex structure on 
$M_{\mathbb{R}}=M\otimes _{\mathbb{Z}}\mathbb{R}$ defined by 
$J(\varphi +\overline{\varphi})=-i\varphi +i\overline{\varphi}$, then 
$\pi^{0,1}:(M_{\mathbb{R}},J)\to \overline{U}$ is a $\mathbb{C}$-linear 
isomorphism. We want to verify that $E'(J\omega ,\omega) > 0$ for 
every nonzero $\omega\in M_{\mathbb{R}}$. Indeed, if 
$\omega=\varphi+\overline{\varphi}$ then 
\[
-Q_{\mathbb{R}}(J\omega,\omega) = 
-Q_{\mathbb{C}}(-i\varphi+i\overline{\varphi}, 
\varphi+\overline{\varphi}) = 
-(-iQ_{\mathbb{C}}(\varphi,\overline{\varphi})+iQ_{\mathbb{C}}(
\overline{\varphi},\varphi)) = 
2iQ_{\mathbb{C}}(\varphi,\overline{\varphi}) > 0
\]
\end{block}
\begin{exa}\label{s4.59}
Let $X$ be a compact Riemann surface of genus $g$. Let 
$M=H^1(X,\mathbb{Z}), U = H^{1,0}(X), M_{\mathbb{C}} = H^1(X,\mathbb{C}) = 
H^{1,0}(X) 
\oplus H^{0,1}(X) = U \oplus \overline{U}$ be the Hodge decomposition, 
$Q(\varphi, \psi) = \int_{X}\varphi\wedge \psi$. Then 
$P(U)=H^{0,1}(X)/\pi^{0,1}H^1(X,\mathbb{Z}) \cong Pic^{0}(X)$ is polarized by 
an 
hermitian form $H'$ with $E'=Im(H')$ where $E'(\varphi, \psi) = -
\int_{X}\varphi\wedge \psi$.
\end{exa}
\begin{block}\label{s4.59bis}
Consider a polarized Hodge structure as in (\ref{s4.58}). Let 
$\Lambda=M^{*}=Hom_{\mathbb{Z}}(M,\mathbb{Z})$. Let 
$V=(\overline{U})^{\bot}\cong U^{*}, \overline{V}=U^{\bot}\cong 
\overline{U}^{*}$. The transposed Weyl operator $^{t}C$ has eigenvalues 
$i$, $-i$ on  $V$ and $\overline{V}$ respectively. Consider the 
corresponding splitting 
$\Lambda_{\mathbb{C}}=V\oplus \overline{V}$.
 The complex torus $A(U)=V/\pi_V(\Lambda)$ is dual to 
$P(U)$. Indeed, define a complex structure on $\Lambda_{\mathbb{R}}$ 
by $I= {^{t}C}$, i.e. $I(v+\overline{v})=iv-i\overline{v}$. Then the 
$\mathbb{C}$-linear isomorphism $\pi_{V}:(\Lambda_{\mathbb{R}},I) \to V$ 
induces an isomorphism of complex tori
\(
\Lambda_{\mathbb{R}}/\Lambda \overset{~}{\lto} V/\pi_{V}(\Lambda) = A(U)
\)
As we saw above we have $(M_{\mathbb{R}},J)/M \cong P(U)$.  The 
$\mathbb{R}$-extension of the perfect pairing $\langle \, ,\, \rangle : 
\Lambda 
\times M \to 
\mathbb{Z}$ satisfies $\langle Iv,J\varphi\rangle = \langle v,\varphi
 \rangle$. Thus $A(U)$ is dual to $P(U)$. Let 
$\omega_{1},\ldots,\omega_{g}$ be a $\mathbb{C}$-basis of $U$ and let 
$\gamma_{1},\ldots,\gamma_{2g}$ a $\mathbb{Z}$-basis of $\Lambda$. Let 
$\{\gamma_{\rho}^{*}\}$ be the dual basis of $M$ and let 
$\Pi=(\pi_{\alpha \rho})$ be the $g\times 2g$ period matrix with 
entries $\pi_{\alpha \rho} = \langle \gamma_{\rho},\omega_{\alpha 
}\rangle$. Then $\omega_{\alpha }=\sum_{\rho}\pi_{\alpha 
\rho}\gamma_{\rho}^{*}$. In matrix form this can be written as 
\(
^{t}(\omega_{1},\ldots ,\omega_{g}) = \Pi \, {^{t}(\gamma_{1}^{*},\ldots 
,\gamma_{2g}^{*})}.
\)
The matrix $\Pi$ is the period matrix of the torus $A(U)$. In 
Example~\ref{s4.59} we have $\Lambda=H_{1}(X,\mathbb{Z}), V\cong 
H^{1,0}(X)^{*}, A(U)=V/\pi_{V}(\Lambda)$ is the Jacobi variety $J(X)$ of 
$X$.
\end{block}

\begin{block}\label{s4.60}
The material of this paragraph is related to \cite{BL2},\cite{BL3}. 
Let $\Lambda$ be a lattice of rank $2g$. Let $E$ be a nondegenerate, integer 
valued, skew-symmetric form on $\Lambda$. Let $\phi: \Lambda_{\mathbb{R}} \to 
\Lambda_{\mathbb{R}}^*$ be the associated isomorphism $v\mapsto E(v,-)$. We 
define a skew-symmetric form 
 $E^{*}$ on $\Lambda_{\mathbb{R}}^{*}$ by the 
equality $E=\phi^{*}E^{*}$, i.e. $E(v,w) = 
E^{*}(\phi(v),\phi(w))$. We need to multiply $E^*$ by an integer in order to 
make it integer-valued on $\Lambda^*$. Let $A=(a_{\rho \sigma})$ be 
the matrix of 
$E$ in a basis $\gamma_{1},\ldots ,\gamma_{2g}$ of $\Lambda$. Let 
$\{\gamma_{\sigma}^{*}\}$ be the dual basis of $\Lambda^{*}$, $\langle 
\gamma_{\rho},\gamma_{\sigma}^{*}\rangle=\delta_{\rho \sigma}$. Since 
$\phi(\gamma_{\rho})= 
\sum_{\sigma}a_{\rho\sigma}\gamma_{\sigma}^{*}$ one has that the 
matrix of $E^{*}$ in the basis $\{\gamma_{\sigma}^{*}\}$ satisfies 
$A=A\, A^{*}\, {^{t}A}$. Hence $A^{*}= -A^{-1}$. Let $\lambda_{1},\ldots ,
\lambda_{2g}$ be a simplectic basis of $\Lambda$, thus 
\(
A=\left(
\begin{smallmatrix}
0 & D \\
-D & 0
\end{smallmatrix}
\right)
\)
where $D=diag(d_{1},\ldots,d_{g})$ with $d_{i}|d_{i+1}$ and let $e=d_{g}$.  
Then 
\(
A^{*}=\left(
\begin{smallmatrix}
0 & D^{-1} \\
-D^{-1} & 0
\end{smallmatrix}
\right)
\)
and hence $eA^{*}$ is with integer entries. We see that $\hat{E}=eE^{*}$ 
is an integer valued skew-symmetric form on $\Lambda^*$ with 
elementary divisors 
$(\hat{d_{1}},\ldots,\hat{d_{g}})$,\; 
$\hat{d_{i}}=\frac{d_g}{d_{g-i+1}}$. 
If 
$A=\left(a_{\rho\sigma}\right)$ is the matrix of $E$ in the basis 
$\{\gamma_{\rho}\}$ of $\Lambda$ and if $A^{-1}=(b_{\sigma\tau})$, then 
$-eA^{-1}=(-eb_{\sigma\tau})$ is the matrix of $\hat{E}$ in the dual 
basis $\gamma_{\sigma}^{*}$ of $\Lambda^{*}$. Let 
$\hat{\phi}:\Lambda_{\mathbb{R}}^* \to \Lambda_{\mathbb{R}}$ be the 
isomorphism associated 
with $\hat{E}$, $\omega \mapsto \hat{E}(\omega,-)$. One has 
$\hat{\phi}(\gamma_{\sigma}^{*})=\sum_{\tau}-eb_{\sigma\tau}\gamma
_{\tau}$, so 
\begin{equation}\label{es4.62}
\hat{\phi}\circ \phi = -e\, \mathbf{1}, \quad
\phi \circ \hat{\phi} = -e\, \mathbf{1}.
\end{equation}
The equality $\phi^{*}\hat{E}=e\, E$ and \eqref{es4.62} imply 
$\hat{\phi}^{*}E=e\, \hat{E}$, so $\hat{E}^{*}=\frac{1}{e}E$ and since 
the maximal elementary divisor of $\hat{E}$ is $\frac{e}{d_{1}}$ one 
obtains
 $(\hat{E})\sphat=\frac{1}{d_1}E$. Thus if 
$d_1=1$ one has $(\hat{E})\sphat=E,\; d_g = e = \hat{d_g}$ and the following 
relations hold
\begin{equation}\label{es4.62a}
\phi^{*}\hat{E}=e\, E, \quad \hat{\phi}^{*}E=e\, \hat{E}
\end{equation}
Let 
$P=\mathbb{C}^{g}/\Lambda$ be an abelian variety with polarization 
$E=Im\, H$. It is convenient to consider $\mathbb{C}^{g}$ as the real 
space $\Lambda_{\mathbb{R}}$ and the multiplication by $i$ as a linear 
operator $I:\Lambda_{\mathbb{R}}\to \Lambda_{\mathbb{R}}$ with 
$I^{2}=-\mathbf{1}$. The Riemann conditions for $E$ are:

(i) $E(\Lambda,\Lambda)\subset \mathbb{Z}$,\; (ii) $E(Iv,Iw) = E(v,w)$,\;
(iii) $E(Iv,v)>0$ for $\forall v\ne 0$.

\noindent
Let $\hat{P}=\Lambda_{\mathbb{R}}^{*}/\Lambda^{*}$ be the dual complex 
torus. Here the complex structure on $\Lambda_{\mathbb{R}}^{*}$ is 
defined by $J=-\, ^{t}I$ so that $\langle Iv,J\omega\rangle = \langle 
v,\omega\rangle$. The mapping $v\mapsto E(v,-)$ yields a 
$\mathbb{C}$-isomorphism 
$\phi: (\Lambda_{\mathbb{R}},I) \overset{\sim}{\lto} 
(\Lambda_{\mathbb{R}}^{*},J)$ and an isogeny $\varphi:P\to \hat{P}$ such that 
$\varphi_*=\phi$. If 
the polarization is defined by an invertible sheaf $L$ on $P$ then 
$\varphi=\varphi_{L}$ where $\varphi(x)=T_{x}^{*}L\otimes L^{-1}\in 
Pic^{0}P = \hat{P}$ (see \cite{Ke}~p.7).  Now the above construction 
yields a 
polarization $\hat{E}$ on $\hat{P}$ such that the polarization mapping 
$\hat{\varphi}:\hat{P}\to P$ satisfies
\(
\hat{\varphi}\circ \varphi = -e\, \mathbf{1}_{P}, \quad
\varphi \circ \hat{\varphi} = -e\, \mathbf{1}_{\hat{P}}.
\)

With every simplectic basis $\{\lambda_{i}\}$ of $\Lambda,\; 
E(\lambda_{i},\lambda_{g+j})=d_{i}\delta_{ij}$, one may associate a 
normalized period matrix $(Z,D)$ of $P$. Here 
$D=diag(d_{1},\ldots,d_{g}),\; Z$ belongs to the Siegel upper half 
space $\mathfrak{H}_{g}$ and is defined by the equalities \cite{LB}~p.213
\begin{equation}\label{es4.63a}
\lambda_{j} = \sum_{i=1}^{g}Z_{ij}\frac{1}{d_{i}}\lambda_{g+i},\quad 
j=1,\ldots,g
\end{equation}
One has $\phi(\lambda_{j}) = d_{j}\lambda_{g+j}^{*},\: 
\phi(\lambda_{g+j}) = -d_{i}\lambda_{i}^{*}$, so applying $\phi$ 
to \eqref{es4.63a} and dividing by $d_{j}$ one obtains
\[
\lambda_{g+j}^{*} = 
\sum_{i=1}^{g}\frac{1}{d_{j}}Z_{ij}(-\lambda_{i}^{*}) = 
\sum_{i=1}^{g}e\frac{1}{d_{i}}Z_{ij}\frac{1}{d_{j}}
\left(\frac{e}{d_{i}}\right)^{-1}(-\lambda_{i}^{*}) = 
\sum_{i=1}^{g}Z_{ij}'
\left(\frac{e}{d_{i}}\right)^{-1}(-\lambda_{i}^{*})
\]
where $Z'=eD^{-1}ZD^{-1}$. Now, let 
$\hat{\lambda}_{i} = \lambda_{2g-i+1}^{*},\: 
\hat{\lambda}_{g+i} = -\lambda_{g-i+1}^{*}, \: \hat{d}_{i} = 
\frac{e}{d_{g-i+1}},\: \hat{Z}_{ij} = Z'_{g-i+1,g-j+1}$. Then 
$\{\hat{\lambda}_{i}\}$ is a simplectic basis of $\Lambda^{*}$ for 
$\hat{E}$ with elementary divisors 
$(\hat{d}_{1}=1,\hat{d}_{2},\ldots,\hat{d}_{g})$ and the corresponding 
normalized period matrix of $\hat{P}$ is $(\hat{Z},\hat{D})$ where 
$\hat{D}=diag(\hat{d}_{1},\ldots,\hat{d}_{g})$,  and

\begin{equation}\label{es4.63}
\hat{Z} = S(eD^{-1}ZD^{-1})S \qquad \text{with} \qquad
 S =
\begin{pmatrix}
0&\dots &0&1\\
0&\dots &1&0\\
\hdotsfor{4}\\
1&\dots &0&0
\end{pmatrix}
\end{equation}
 \end{block}

\begin{block}\label{s4.63}
Let 
$M,\: M_{\mathbb{C}}=U\oplus \overline{U},\: Q:M\times M\to 
\mathbb{Z}$ be a polarized Hodge structure of weight one. Let $\Lambda 
= M^{*}$, let $\psi :M_{\mathbb{R}}\to \Lambda_{\mathbb{R}}$ be the 
linear isomorphism $\psi(m)=Q(m,-)$ and let $\hat{Q}:\Lambda \times 
\Lambda \to \mathbb{Z}$ be the dual form as defined in (\ref{s4.60}). 
Since 
$Q_{\mathbb{R}}(C\omega_{1},C\omega_{2})=Q_{\mathbb{R}}(\omega_{1},
\omega_{2})$ one obtains a $\mathbb{C}$-linear isomorphism $\psi : 
(M_{\mathbb{R}},C)\overset{\sim}{\lto}(\Lambda_{\mathbb{R}},-{^{t}C})
$. With respect to $-{^{t}C}$ one has the splitting 
$\Lambda_{\mathbb{C}}=W\oplus \overline{W}$ with eigenvalues $i$ and 
$-i$ respectively, where $W=\overline{V}=U^{\bot}$ and 
$\overline{W}=V=(\overline{U})^{\bot}$ (cf. (\ref{s4.59bis})). 
Furthermore $(\Lambda,\Lambda_{\mathbb{C}}=W\oplus 
\overline{W},\hat{Q})$ is a polarized Hodge structure of weight one. 
We call it the \emph{dual polarized Hodge structure} of 
$(M,M_{\mathbb{C}}=U\oplus \overline{U},Q)$. One has $P(W)=A(U),\; 
A(W) = P(U)$.
\end{block}

\noindent
We summarize the arguments of (\ref{s4.58}) -- (\ref{s4.63}) in the 
following statement.
\begin{pro}\label{s4.64}
Let $M,\: M_{\mathbb{C}}=U\oplus \overline{U},\: Q:M\times M\to 
\mathbb{Z}$ be a polarized Hodge structure of weight one. Then one may 
associate to it a 
 pair of 
dual abelian varieties $A(U) = U^{*}/\pi_{_{U^{*}}}(\Lambda),\; P(U) = 
\overline{U}/\pi^{0,1}(M)$ with polarizations $E$, resp. $\hat{E}$ 
which have types $(d_{1}=1,d_{2},\ldots,d_{g})$, resp. 
$(\hat{d}_{1}=1,\hat{d}_{2},\ldots,\hat{d}_{g})$ where 
$\hat{d}_{i} = d_{g}/d_{g-i+1}$. Furthermore $\hat{E} = 
-\frac{1}{c}Q$ where $c$ is the first elementary divisor of $Q$ and $E 
= -\hat{Q}$. In 
appropriate simplectic bases the normalized period matrices of $A(U)$ 
and $P(U)$ are respectively $\Pi = (Z,D),\; \hat{\Pi} = 
(\hat{Z},\hat{D})$ where $\hat{Z} = S(d_{g}D^{-1}ZD^{-1})S$, $D = 
diag(d_{1},\ldots,d_{g}),\; \hat{D} = diag(\hat{d}_{1},\ldots,\hat{d}_g)$ and 
$S$ is the matrix 
defined in \eqref{es4.63}. The dual polarized Hodge structure is 
defined by $\Lambda = M^{*}, \Lambda_{\mathbb{C}} = W\oplus 
\overline{W}, \hat{Q}$, where $W = U^{\bot}, \overline{W} = 
(\overline{U})^{\bot}$ and one has $A(U) = P(W),\; P(U) = A(W)$.
\end{pro}

\begin{exa}\label{s4.65}
Consider $A(U)\cong J(X)$ from (\ref{s4.59}) and (\ref{s4.59bis}). Let 
$E:\Lambda \times \Lambda \to \mathbb{Z}$ be the skew-symmetric form 
$E(\gamma,\delta)= -(
\gamma,\delta
)_{X}$. Then $\phi : \Lambda \to \Lambda^{*} = M$ equals $\phi = 
-D_{X}:H_1(X,\mathbb{Z})\to H^{1}(X,\mathbb{Z})$, where $D_{X}$ is the 
Poincar\'{e} isomorphism. From the property $(\gamma,\delta)_X = \linebreak
\int_{X}D_{X}(\gamma)\wedge D_{X}(\delta)$ (cf. \cite{GH}~Ch.0) one 
concludes $\hat{E}=E'$ where $E'$ was defined in Example~\ref{s4.59}. 
Hence $-(\;,\;)_X = E=(E')\sphat$\; is the canonical principal polarization of 
$J(X)$. The period $g\times g$ 
matrix $Z$ as defined in (\ref{s4.60}) is the same as the classical one. 
In fact if $\{A_{i}, B_{j}\}_{i,j=1}^g$ is a standard system of 
cycles on $X,\; (A_{i},B_{j})_{X}=\delta_{ij}$, then 
$\{\lambda_{i}=B_{i}, \lambda_{g+j}=A_{j}\}$ is a simplectic basis for 
$E$. Thus if $\omega_{1},\ldots,\omega_{g}$ is a normalized basis of 
differentials, $\int_{A_{j}}\omega_{i}=\delta_{ij}$ then 
$Z_{ij}=\int_{B_{j}}\omega_{i}$. The polarized Hodge structure dual to
the one of Example~\ref{s4.59} is $H_1(X,\mathbb{Z}),H_1(X,\mathbb{C}) 
= W \oplus \overline{W}, \hat{Q}$, where $W=H^{1,0}(X)^{\bot},
\overline{W} = H^{0,1}(X)^{\bot},\; 
\hat{Q}(\gamma,\delta) = (\gamma,\delta)_{X}$.
\end{exa}

\begin{block}\label{s4.66}
Let $\pi :X\to Y$ be a covering of smooth projective curves, $g(Y)\geq 
1$. Then 
$\pi^*: H^1(Y,\mathbb{Z})\to H^1(X,\mathbb{Z})$ and $ 
^{t}\pi^{*} = \pi_{*}: H_1(X,\mathbb{Z}) \to H_1(Y,\mathbb{Z})$ induce 
morphisms of the corresponding Hodge structures. We let $M = 
H^1(X,\mathbb{Z})/H^1(X,\mathbb{Z})\cap 
\pi^{*}H^1(Y,\mathbb{R})$ with a dual lattice $\Lambda = 
Ker(\pi_{*}:H_1(X,\mathbb{Z})\to H_1(Y,\mathbb{Z}))$. The Hodge 
structures $M_{\mathbb{C}} = U\oplus \overline{U}$ and 
$\Lambda_{\mathbb{C}} = W\oplus \overline{W}$ are defined respectively 
by $U = H^{1,0}(X)/\pi^{*}H^{1,0}(Y)$, $W = H^{1,0}(X)^{\bot}\cap 
\Lambda_{\mathbb{C}}$. The corresponding pair of dual abelian 
varieties is 
\[
A(U) = Ker(Nm_{\pi}:J(X)\to J(Y))^{0}, \qquad 
P(U) = Pic^{0}(X)/\pi^{*}Pic^{0}(Y).
\]
Consider the restriction of $-(\; ,\; )_{X}$ on $\Lambda$. It is a 
nondegenerate form as evident from the orthogonal decomposition 
$H_1(X,\mathbb{R})=Ker\, \pi_{*}\oplus \pi^{*}H_1(Y,\mathbb{R})$. 
Dividing it by its smallest elementary divisor we obtain a polarization 
$E:\Lambda\times \Lambda\to \mathbb{Z}$ on $A(U)$. The dual form 
$\hat{E}:M\times M\to \mathbb{Z}$ is a polarization on $P(U)$, both 
$E$ and $\hat{E}$ have first elementary divisor 1, have the same exponent and 
their 
types are 
related as in Proposition~\ref{s4.64}. The Hodge structures 
$M_{\mathbb{C}}=U\oplus \overline{U}$ and $\Lambda_{\mathbb{C}}=W\oplus 
\overline{W}$ are polarized respectively by $Q = -\hat{E}$ and 
$\hat{Q} = -E$
\end{block}
\begin{exa}\label{s4.67a}

(i) Let $\pi :X\to Y$ be a covering of smooth, projective curves of 
prime degree $d$, let $g(X)\geq 3,\; g(Y)=1$. Then by Lemma~\ref{s1.1}
$P=Ker(Nm_{\pi })$ is connected  and  the 
polarization $E$ induced from $J(X)$  has type $(1,\ldots ,1,d)$. 
hence the dual abelian variety $\hat{P} = 
Pic^{0}(X)/\pi^{*}Pic^{0}(Y)$ has dual polarization $\hat{E}$ of type 
$(1,d,\ldots ,d)$.

(ii) Let $\pi :X\to Y$ be a double covering of smooth, projective curves where $g(X)=7,\; g(Y)=3$ as in 
\cite{BCV}. Then $Ker\,Nm_{\pi }$ is connected and equals the Prym variety
 $P$ with 
induced polarization $E$ of type $(1,2,2,2)$. Thus the dual variety 
$\hat{P} = Pic^{0}X/\pi^{*}Pic^{0}Y$ has dual polarization $\hat{E}$ 
of type $(1,1,1,2)$.
\end{exa}
\begin{block}\label{s4.67}
We want to adapt the arguments in \cite{Gr}~pp.576,577 to the case 
of arbitrary polarizations. Let $M,\; M_{\mathbb{C}}=U\oplus 
\overline{U},\; Q:M\times M\to \mathbb{Z}$ be a polarized Hodge 
structure. Let $\Lambda=M^{*},\: \Lambda_{\mathbb{C}}=W\oplus 
\overline{W},\; \hat{Q}$ be the dual Hodge structure and let 
$-\hat{Q}=E:\Lambda\times \Lambda\to \mathbb{Z}$ and $\hat{E}:M\times 
M\to \mathbb{Z}$ be the associated skew-symmetric forms as defined in 
Proposition~\ref{s4.64}.
Let   $\{\omega_{\alpha }\}$ be a basis of $U$, let 
$\{\lambda_{\rho}\}$ be a simplectic basis of $\Lambda$ 
with respect to $E$
and let $\Pi$ 
be the corresponding period matrix 
\begin{equation}\label{es4.67}
^{t}(\omega_{1},\ldots ,\omega_{g}) = \Pi\, {^{t}(\lambda_{1}^{*},\ldots 
,\lambda_{2g}^{*})}.
\end{equation}
 The matrix of 
$\hat{E}$ in $\{\lambda_{\rho}^*\}$ is 
\(
(\hat{E}(\lambda_{\rho}^{*},\lambda_{\sigma}^{*}))=\left(
\begin{smallmatrix}
0 & eD^{-1} \\
-eD^{-1} & 0
\end{smallmatrix}
\right)
\)
where $D=diag(d_{1},\ldots,d_{g})$. Since $Q=-c\hat{E},\; c\in 
\mathbb{N}$, the Riemann relations of (\ref{s4.58}) may be written as 
\(
(a)\; 
(\hat{E}(\omega_{\alpha },\omega_{\beta})) = 0,
(b)\; 
i(\hat{E}(\omega_{\alpha },\overline{\omega}_{\beta })) < 0
\)
or equivalently in matrix form for $\Pi = (\Pi_{1}\: \Pi_{2})$
\[
\begin{pmatrix}
\Pi_{1}&\Pi_{2}
\end{pmatrix}
\begin{pmatrix}
0 & D^{-1}\\
-D^{-1}& 0
\end{pmatrix}
\begin{pmatrix}
^{t}\Pi_{1}\\
^{t}\Pi_{2}
\end{pmatrix}
= 0, \qquad
i
\begin{pmatrix}
\Pi_{1}&\Pi_{2}
\end{pmatrix}
\begin{pmatrix}
0 & D^{-1}\\
-D^{-1}& 0
\end{pmatrix}
\begin{pmatrix}
^{t}\overline{\Pi}_{1}\\
^{t}\overline{\Pi}_{2}
\end{pmatrix}
< 0.
\]
In this form the relations are the same as in \cite{LB}~p.77. 
Changing the basis of $U$ by $^{t}(\phi_{1},\ldots,\phi_{g})=A\, 
^{t}(\omega_{1},\ldots,\omega_{g})$ and the simplectic basis of 
$\Lambda$ by 
$(\mu_{1},\ldots,\mu_{2g})=(\lambda_{1},\ldots,\lambda_{2g})R$ one 
obtains an equivalent period matrix $\Pi'\sim \Pi,\; \Pi' = A\Pi R$. 
The same argument as in \cite{Gr}~p.577 shows that $\det \Pi_{1}\ne 
0,\; \det \Pi_{2}\ne 0$, so every $\Pi$ is equivalent to 
$D\Pi_{2}^{-1}(\Pi_{1}\: \Pi_{2}) = (Z\: D)$ with $Z\in 
\mathfrak{H}_{g}$ as follows from Riemann's relations. The matrix 
$R\in M_{2g}(\mathbb{Z})$ satisfies the equality 
\(
^{t}R\left(
\begin{smallmatrix}
0 & D\\
-D & 0
\end{smallmatrix}
\right)
R = \left(
\begin{smallmatrix}
0 & D\\
-D & 0
\end{smallmatrix}
\right).
\)
If we let 
\(
R = \:
{^{t}\left(
\begin{smallmatrix}
a&b\\
c&d
\end{smallmatrix}
\right)}
\), then $\left(
\begin{smallmatrix}
a&b\\
c&d
\end{smallmatrix}
\right)\in \Gamma_{D}=Sp^{D}_{2g}(\mathbb{Z})$, the group defined in 
\cite{LB}~p.219. Multiplying on the right a normalized period matrix 
$(Z\: D)$ by $R= {^{t}\left(
\begin{smallmatrix}
a&b\\
c&d
\end{smallmatrix}
\right)}$ and then normalizing one obtains $(Z'\: D)$ where 
\(
Z' = (aZ+bD)(D^{-1}cZ + D^{-1}dD)^{-1}
\). This is the left action of $\Gamma_{D}$ on $\mathfrak{H}_{g}$ 
defined in \cite{LB}~p.219. The quotient $\Gamma_{D}\backslash 
\mathfrak{H}_{g}$ is the moduli space $\mathcal{A}_{D}$ for polarized 
abelian varieties of type $D$ (ibid). In conclusion one obtains a 
correspondence
\begin{equation}\label{es4.69}
(M,M_{\mathbb{C}}=U\oplus \overline{U}, Q)\quad \mapsto \quad Z(mod\: 
\Gamma_{D}) 
\in \mathcal{A}_{D}.
\end{equation}
Using the dual polarized Hodge structure 
$(\Lambda,\Lambda_{\mathbb{C}}=W\oplus \overline{W}, \hat{Q})$
 one obtains similarly another 
correspondence
\begin{equation}\label{es4.69a}
(M,M_{\mathbb{C}}=U\oplus \overline{U}, Q)\quad \mapsto \quad \hat{Z}(mod\: 
\Gamma_{\hat{D}}) \in \mathcal{A}_{\hat{D}}.
\end{equation}
Comparing with Proposition~\ref{s4.64} we see that \eqref{es4.69} and 
\eqref{es4.69a} associate to a polarized Hodge structure of weight one 
respectively the isomorphism classes $[A(U)]\in \mathcal{A}_{D}$ and 
$[P(U)]\in \mathcal{A}_{\hat{D}}$.
\end{block}
\begin{block}\label{s4.70}
Let $T$ be a connected complex manifold. A polarized variation of 
Hodge structure of weight one (VHS) over $T$ is given by the following 
data.
\renewcommand{\theenumi}{\roman{enumi}}
\begin{enumerate}
\item
A flat bundle of rank $2g$ lattices $\mathbb{M}\to T$.
\item 
A holomorphic rank $g$ subbundle $\mathbb{F}\subset 
\mathbb{M}_{\mathbb{C}}$ such that $\mathbb{F}\oplus 
\overline{\mathbb{F}}=\mathbb{M}_{\mathbb{C}}$.
\item
A flat skew-symmetric form $Q:\mathbb{M}\times \mathbb{M}\to \mathbb{Z}$
which satisfies fiberwise the Riemann relations of (\ref{s4.58}).
\end{enumerate}
Given a VHS one may consider the dual  VHS \quad
$\mathbb{L}=Hom_{\mathbb{Z}}(\mathbb{M},\mathbb{Z}),\; 
\mathbb{L}_{\mathbb{C}}=\mathbb{G}\oplus \overline{\mathbb{G}}$, 
where $\mathbb{G}=(\overline{\mathbb{F}})^{\bot},\; 
\overline{\mathbb{G}}=(\mathbb{F})^{\bot}$ 
and the polarization 
$\hat{Q}$ is obtained from $Q$ as described in the beginning of 
(\ref{s4.60}). Dividing by an appropriate negative integer $-c$ one obtains 
flat, integer valued, skew-symmetric forms $E=-\hat{Q},\quad 
\hat{E}=-\frac{1}{c}Q$ such that for each $s\in T$ the forms 
$E_{s},\: \hat{E}_{s}$ are 
respectively
polarizations of types $D = (1,d_{2},\ldots,d_{g}),\; 
\hat{D}=(1,\hat{d}_{2},\ldots,\hat{d}_{g})$
of the 
associated complex tori $A(\mathbb{F}_{s}),\: 
P(\mathbb{F}_{s})$  as in Proposition~\ref{s4.64}. Let $S\subset T$ be 
an open set in the Hausdorff topology
over which $\mathbb{M}$ and $\mathbb{F}$ are trivial. 
Choosing a frame $\lambda_{1},\ldots,\lambda_{2g}$ of 
$\mathbb{L}|_{S}$, a frame $\omega_{1},\ldots, \omega_{g}$ of 
$\mathbb{F}|_{S}$ and normalizing the associated period matrices as in 
(\ref{s4.67}) one obtains a holomorphic mapping $\tilde{\varPhi}_S: S\to 
\mathfrak{H}_{g}$. Passing to the quotient $\Gamma_{D}\backslash 
\mathfrak{H}_{g}=\mathcal{A}_{D}$ one obtains a holomorphic mapping 
$\varPhi _{S}:S\to \mathcal{A}_{D}$. Covering $T$ with such open sets 
$T=\cup S_{i}$ and gluing $\varPhi _{S_{i}}$ as evident from 
(\ref{s4.67}) one obtains the period mapping 
$\varPhi :T\to \mathcal{A}_{D}$. 
One has 
$\varPhi (s)=[A(\mathbb{F}_{s})]$. The same 
argument applied to the dual  VHS yields a dual period mapping 
$\hat{\varPhi }:T\to \mathcal{A}_{\hat{D}}$ such that $\hat{\varPhi 
}(s)=[P(\mathbb{F}_{s})]$. This proves part of the following 
statement. 
\end{block}
\begin{pro}\label{s4.72}
Let $T$ be a connected complex manifold and let $\mathbb{M},\; 
\mathbb{M}_{\mathbb{C}}=\mathbb{F}\oplus \overline{\mathbb{F}},\; 
Q:\mathbb{M}\times \mathbb{M}\to \mathbb{Z}$ be a polarized variation 
of Hodge structures of weight one. Let $D$ and $\hat{D}$ be the dual 
polarization types as defined in (\ref{s4.70}). Then one can define 
period mappings $\varPhi :T\to \mathcal{A}_{D}$ and 
$\hat{\varPhi }:T\to \mathcal{A}_{\hat{D}}$ which transform $s\in T$ 
respectively into the isomorphism classes of polarized abelian 
varieties $\varPhi (s)=[A(\mathbb{F}_{s})],\; 
\hat{\varPhi }(s) = [P(\mathbb{F}_{s})]$. If $T$ is algebraic, so are the 
period mappings $\varPhi $ and 
$\hat{\varPhi }$. The mapping $\varPhi $ is dominant if and only if 
$\hat{\varPhi }$ is dominant.
\end{pro}
{\bf Proof.}
That $\varPhi $ and $\hat{\varPhi }$ are algebraic if $T$ is algebraic 
follows from Borel's extension theorem \cite{Bo}~Theorem~3.10. The 
last statement follows from comparing the differentials of $\varPhi $ 
and $\hat{\varPhi }$ by means of \eqref{es4.63}.
\hfill $\Box$

\begin{block}\label{s4.73}
We now consider a family of coverings of curves and associate to it two 
dual VHS. Suppose we are given a commutative diagram of holomorphic 
mappings
\begin{equation}\label{es4.73}
\begin{diagram}
\mathcal{X}&    &\rTo^p &       &\mathcal{Y} \\
           &\rdTo_f &     & \ldTo_q  &       \\
           &      &  T  &        &           \\
\end{diagram}
\end{equation}
where $T$ is a connected complex manifold, $f$ and $q$ are smooth, 
proper of relative dimension one and $p$ is surjective. Then we have 
the standard VHS associated with $f$ and $q$: \quad 
$\mathbb{H}_{\mathcal{X}} = R^{1}f_{*}\mathbb{Z},\; 
\mathbb{F}_{\mathcal{X}}\subset \mathbb{H}_{\mathcal{X}}\otimes \mathbb{C}$ 
where 
$\mathcal{O}_{T}(\mathbb{F}_{\mathcal{X}}) = 
f_{*}\Omega^{1}_{\mathcal{X}/T}$ and similarly for $q$. One has a 
morphism of VHS $p^{*}:
(\mathbb{H}_{\mathcal{Y}},\mathbb{F}_{\mathcal{Y}})
\to 
(\mathbb{H}_{\mathcal{X}},\mathbb{F}_{\mathcal{X}})
$. Define $\mathbb{M},\mathbb{F}\subset \mathbb{M}_{\mathbb{C}}$ as 
follows (cf. (\ref{s4.66}))
\[
\mathbb{M}=R^{1}f_{*}\mathbb{Z}/
R^{1}f_{*}\mathbb{Z}\cap p^{*}R^{1}q_{*}\mathbb{R}, \quad
\mathcal{O}_{T}(\mathbb{F})=
f_{*}\Omega^{1}_{\mathcal{X}/T}/p^{*}q_{*}\Omega^{1}_{\mathcal{Y}/T}.
\]
We obtain a dual VHS letting $\mathbb{L} = \mathbb{M}^{*},\; 
\mathbb{G} = (\mathbb{F})^{\bot}$. Let 
$\tilde{Q}:R^{1}f_{*}\mathbb{Z}\times R^{1}f_{*}\mathbb{Z}\to 
R^{2}f_{*}\mathbb{Z} = \mathbb{Z}$ be the cup-product which is an
unimodular polarization of $(\mathbb{H}_{\mathcal{X}},
\mathbb{F}_{\mathcal{X}}
)$ (cf. Example~\ref{s4.59}). 
According to  (\ref{s4.63}) we 
let $\tilde{Q}\sphat :\mathbb{H}_{\mathcal{X}}^{*}\times 
\mathbb{H}_{\mathcal{X}}^{*} \to \mathbb{Z}
$ be the dual polarization. Since $\mathbb{M}$ is a quotient of $
\mathbb{H}_{\mathcal{X}}
$ the dual $\mathbb{L} = \mathbb{M}^{*}$ may be embedded in $
\mathbb{H}_{\mathcal{X}}^{*}
$. Restricting $\tilde{Q}\sphat$ on $\mathbb{L}$ and dividing the obtained 
flat 
skew-symmetric form by its least elementary divisor one obtains a polarization $\hat{Q}:\mathbb{L}\times \mathbb{L}\to \mathbb{Z}$. Its dual 
$Q:\mathbb{M}\times \mathbb{M}\to \mathbb{Z}$ polarizes the VHS
$(\mathbb{M},\mathbb{F})$. The flat forms $E=-\hat{Q}$ and 
$\hat{E}=-Q$ polarize fiberwise respectively the associated complex 
tori $Ker(Nm_{p_{s} })^{0} $ and 
$Pic^{0}\mathcal{X}_{s}/p_{s}^{*}Pic^{0}\mathcal{Y}_{s}$ for $\forall 
s\in T$. Applying Proposition~\ref{s4.72} we obtain the following 
result.
\end{block}
\begin{pro}\label{s4.74}
Let $p: \mathcal{X}\to \mathcal{Y}$ be a covering 
of smooth $T$-curves over a smooth connected algebraic base $T$ (cf. 
\eqref{es4.73}).  Fix $o\in T$  and suppose the restriction of the 
intersection form $(\; ,\; )_{\mathcal{X}_{o}}$ on 
$Ker((p_{o})_{*}:H_1(\mathcal{X}_{o},\mathbb{Z})\to 
H_1(\mathcal{Y}_{o},\mathbb{Z}))$ has elementary divisors 
$(m,md_{2},\ldots,md_{g})$,\; $d_{i}|d_{i+1}$. Let 
$D=(1,d_{2},\ldots,d_{g}),\; 
\hat{D}=(1,\hat{d}_{2},\ldots,\hat{d}_{g})$ where 
$\hat{d}_{i}=d_{g}/d_{g-i+1}$. 
Then the period mappings $\varPhi :T\to \mathcal{A}_{D}$ and 
$\hat{\varPhi }:T\to \mathcal{A}_{\hat{D}}$ given by $\varPhi (s)=
[Ker(Nm_{p_{s} })^{0}]$ and $\hat{\varPhi }(s)=
[Pic^{0}\mathcal{X}_{s}/p_{s}^{*}Pic^{0}\mathcal{Y}_{s}]$ are 
algebraic morphisms. The morphism $\varPhi$ is dominant if and only if 
the morphism $\hat{\varPhi}$ is dominant. If this is the case and if 
$T$ is unirational, then both $\mathcal{A}_{D}$ and 
$\mathcal{A}_{\hat{D}}$ are unirational varieties.
\end{pro}
\begin{rem}\label{s4.75one}
A recent result of Birkenhake and Lange \cite{BL3} shows that 
$\mathcal{A}_{D}$ and $ \mathcal{A}_{\hat{D}}$ are in fact isomorphic 
to each other.
\end{rem}
\begin{exa}\label{s4.75}
In \cite{BCV}~p.124 it is proved that $\mathcal{A}_{4}(1,2,2,2)$ is 
unirational considering the Prym mapping for a family as in 
\eqref{es4.73} where $T$ is a Zariski open set in 
$|\mathcal{O}_{\mathbb{P}^{2}}(4)(-p_{1}-p_{2}-p_{3}-p_{4}|$, where 
$p_{1},\ldots p_{4}$ are general points in $\mathbb{P}^{2}$, \quad 
$q:\mathcal{Y}\to T$ is the corresponding family of plane quartics and 
$p:\mathcal{X}\to \mathcal{Y}$ is a suitable double covering. From 
Proposition~\ref{s4.74} it follows that every general  
abelian variety of dimension four
with polarization of the type $(1,1,1,2)$ is isomorphic 
to $Pic^{0}\mathcal{X}_{s}/p_{s}^{*}Pic^{0}\mathcal{Y}_{s}$ for some 
$s\in T$ and $\mathcal{A}_{4}(1,1,1,2)$ is unirational. The latter 
follows of course from the result of Birkenhake and Lange cited above.
\end{exa}
\begin{que}\label{s4.75two}
Is it true that the moduli space $\mathcal{A}_{4}(1,1,2,2)$ is 
unirational?
\end{que}
\begin{block}\label{s4.75a}
We now wish to give  a formula for the differential of the period 
mapping $\varPhi $ of the VHS considered in  
Proposition~\ref{s4.74}. Let us first consider the general set-up of polarized VHS 
of weight one. Since the problem is local replacing $T$ by a smaller 
open set $S$ we may restrict ourselves to the case where 
$\mathbb{M}=M\times S$ is constant and the holomorphic subbundle 
$\mathbb{F}\subset \mathbb{M}_{\mathbb{C}}$ is trivial. Let us fix a 
basis $\lambda_{1}^{*},\ldots,\lambda_{2g}^{*}$ of $M$ and a frame 
$\omega_{1},\ldots,\omega_{g}$ of $\mathbb{F}$. Transposing 
\eqref{es4.67} we may write 
\[
(\omega_{1}(s),\ldots,\omega_{g}(s))\ =\ (
\lambda_{1}^{*},\ldots,\lambda_{2g}^{*}
)\; ^{t}\Pi(s).
\]
By the first Riemann relation (\ref{s4.58}(iii)) every $g$-plane 
$\mathbb{F}(s)\subset M_{\mathbb{C}}$ is isotropic with respect to 
$Q_{\mathbb{C}}$. Let us denote by $\check{\mathbb{D}}\subset 
Gr(g,M_{\mathbb{C}})$ the simplectic grassmanian of isotropic 
$g$-planes and let $\mathbb{D}\subset 
\check{\mathbb{D}}$ be the open subset of those $g$-planes satisfying 
the second Riemann relation \cite{GS}~p.54. One obtains a holomorphic 
mapping $\tilde{\varPhi }:S\to \mathbb{D}\subset \check{\mathbb{D}},\quad
\tilde{\varPhi }(s)= \mathbb{F}(s) $. 
Considering $g$-planes is equivalent to taking quotient modulo the 
equivalence relation $^{t}\Pi \sim \: ^{t}\Pi\: ^{t}A,\quad A\in 
GL(g,\mathbb{C})$. Thus $\tilde{\varPhi }$ is a coordinate free 
description of the mapping $\tilde{\varPhi }_{S}:S\to 
\mathfrak{H}_{g}$ considered in (\ref{s4.70}).

Let $\mathbb{U}\to \check{\mathbb{D}}$ be the tautological vector 
bundle. The cotangent bundle $(T\, \check{\mathbb{D}})^{*}$ is 
isomorphic to $Sym^{2}\mathbb{U}$. Fiberwise this isomorphism is 
explicitly described as follows. If $z=[U\subset 
M_{\mathbb{C}}]\in \check{\mathbb{D}}$ the vector $\phi \in 
Hom(U,M_{\mathbb{C}}/U) = T_{z}Gr(g,M_{\mathbb{C}})$ belongs to the 
tangent space to $\check{\mathbb{D}}$ if and only if 
$Q(\phi(u),v)+Q(u,\phi(v))=0$ for $\forall u,v \in U$. Since $Q$ is 
skew-symmetric this is equivalent to saying that the bilinear form 
$q_{\phi}(u,v)=Q(\phi(u),v)$ is symmetric. Considering the trilinear 
form
\begin{equation}\label{es4.75bone}
T\, \check{\mathbb{D}}\times \mathbb{U}\times \mathbb{U}\to 
\mathbb{C}, \qquad (\phi,u,v)\mapsto Q(\phi(u),v)
\end{equation}
yields the isomorphisms
\begin{equation}\label{es4.75btwo}
\begin{gathered}
T\, \check{\mathbb{D}} \overset{\sim}{\lto} 
Sym^{2}\mathbb{U}^{*},\quad \phi \mapsto q_{\phi}, \quad 
q_{\phi}(u,v)=Q(\phi(u),v),\\
Sym^{2}\mathbb{U} \overset{\sim}{\lto} (T\, \check{\mathbb{D}})^{*}, 
\qquad \langle \phi,u\otimes v\rangle = Q(\phi(u),v).
\end{gathered}
\end{equation}
 Consider the 
period mapping $\tilde{\varPhi }:S\to \check{\mathbb{D}}$. 
Let $s\in S$. The differential of $\tilde{\varPhi }$ at $s\in S$ yields 
a trilinear form 
\begin{equation}\label{es4.75bthree}
T_{s}S\times \mathbb{F}(s)\times \mathbb{F}(s)\lto \mathbb{C},\qquad
(\frac{\partial}{\partial \tau},\omega_{1},\omega_{2})=
Q(d\, \tilde{\varPhi }\left(\frac{\partial}{\partial 
\tau}\right)(\omega_{1}),\omega_{2}).
\end{equation}
From \eqref{es4.75btwo} one obtains a formula for 
$^{t}d\tilde{\varPhi }(s)$:
\begin{equation}\label{es4.75c}
^{t}d\tilde{\varPhi }(s): Sym^{2}\mathbb{F}(s)\lto (T_{s}S)^{*}, 
\;
\langle \frac{\partial}{\partial \tau},
{^{t}d}\tilde{\varPhi }(\omega_{1}\cdot \omega_{2})\rangle =
Q(d\, \tilde{\varPhi }\left(\frac{\partial}{\partial 
\tau}\right)(\omega_{1}),\omega_{2}),
\end{equation}
where $\cdot$ denotes the product in the symmetric algebra induced 
from $\otimes$.
\end{block}
\begin{block}\label{s4.75c}
Let us apply the above to the VHS associated with a surjective holomorphic 
mapping $p:\mathcal{X}\to \mathcal{Y}$ of smooth families of curves over $T$  
as in 
(\ref{s4.73}). Let $s_{0}\in T$ and let $X=\mathcal{X}_{s_{0}},\quad 
Y=\mathcal{Y}_{s_{0}},\quad \pi=p_{s_{0}}$. We may replace $T$ by a 
smaller neighborhood $S$ of $s_{0}$ such that $f:\mathcal{X}\to S$ and 
$q:\mathcal{Y}\to S$ are $C^{\infty}$-trivial fibrations and all 
bundles occurring in the constructions of (\ref{s4.73}) are trivial (in 
the corresponding category). Abusing notation let us denote by 
$\pi_{*}$ both the homomorphism $\pi_{*}:H_1(X,\mathbb{Z})\to 
H_1(Y,\mathbb{Z})$ and the Gysin homomorphism $D_{Y}\circ \pi_{*}\circ 
D_{X}^{-1}:H^1(X,\mathbb{Z})\to H^1(Y,\mathbb{Z})$, where $D_{X},\; 
D_{Y}$ are the Poincar\'{e} isomorphisms. Let $(\quad)^{-}$ denote the 
kernel of $\pi_{*}$ on the corresponding object. The Gysin 
homomorphism preserves the Hodge type and 
coincides with the 
trace map $Tr_{p_{s}}$ on 
$H^{1,0}(\mathcal{X}_{s}) = 
H^0(\mathcal{X}_{s},\Omega_{\mathcal{X}_{s}})$ . We thus obtain a VHS as 
follows: 
$M'=H_1(X,\mathbb{C})^{-},\; \mathbb{M}'=M'\times S,\; \mathbb{F}^{-} 
= \mathbb{M}'_{\mathbb{C}}\cap \mathbb{F}_{\mathcal{X}} = 
Ker(\pi_{*}:\mathbb{F}_{\mathcal{X}}\to \mathbb{F}_{\mathcal{Y}})$. 
Notice that for the Weyl operators one has $D_{X}\circ 
(-{^{t}C_{\mathcal{X}_{s}}}) = C_{\mathcal{X}_{s}}\circ D_{X}$, thus 
$D_{X}(H_1(X,\mathbb{Z})^{-}) = H^1(X,\mathbb{Z})^{-},\;
D_{X}(\mathbb{G})=\mathbb{F}^{-},\; 
D_{X}(\overline{\mathbb{G}}) = \overline{\mathbb{F}}^{-}$ (cf. (\ref{s4.73})). 
Let $Q^{-}$ be the restriction of $\tilde{Q}$ on $\mathbb{M}'$. Thus 
$Q^{-}$ is the constant skew-symmetric form induced from
\begin{equation}\label{es4.75d}
Q^{-}(\omega_{1},\omega_{2}) = \int_{X}\omega_{1}\wedge 
\omega_{2},\quad \text{where}\quad \omega_{1},\omega_{2}\in 
H^1(X,\mathbb{Z})^{-}.
\end{equation}
The canonical homomorphism $j:H^1(X,\mathbb{Z})^{-}\to 
H^1(X,\mathbb{Z})/\pi^{*}H^1(Y,\mathbb{Z})$ of lattices of equal rank 
has kernel 0 and induces a morphism of VHS\quad 
$(\mathbb{M}',\mathbb{F}')\to (\mathbb{M},\mathbb{F})$. We wish to 
compare the polarization $Q$ on $(\mathbb{M},\mathbb{F})$ as defined 
in (\ref{s4.73}) with $Q^{-}$ from \eqref{es4.75d}.
\end{block}
\begin{lem}\label{s4.75d}
Suppose the restriction of the intersection form $(\; ,\; )_{X}$ on 
$H_1(X,\mathbb{Z})^{-}$ has elementary divisors $(m,md_2,\ldots,md_g),\; 
d_{i}|d_{i+1}$. Then $j^{*}Q=md_{g}Q^{-}$
\end{lem}
{\bf Proof.}
Using the notation of (\ref{s4.66}) and (\ref{s4.73}) we have $\Lambda = 
H_1(X,\mathbb{Z})^{-},\; (\; ,\; )_{X}|_{\Lambda} = -mE$. The 
homomorphism $\varphi :\Lambda \to \Lambda^{*} = M$ is given by 
$\langle \delta ,\varphi(\lambda)\rangle = E(\lambda,\delta)$ and one 
has 
\[
\varphi ^{*}(Q) = \varphi^{*}(-\hat{E}) = \varphi^{*}(-d_{g}E^{*}) = 
-d_{g}E = \frac{d_{g}}{m}(\; ,\; )_{X}|_{\Lambda}.
\]
On the other hand the Poincar\'{e} isomorphism 
$D_{X}:H_1(X,\mathbb{Z})\to H^1(X,\mathbb{Z})$ transforms $\Lambda = 
H_1(X,\mathbb{Z})^{-}$ onto $H^1(X,\mathbb{Z})^{-}$ and one has 
$D_{X}^{*}(Q^{-}) = (\; ,\; )_{X}|_{\Lambda}$. Consider the 
composition mapping
\[
\varphi' :\Lambda \lto H_1(X,\mathbb{Z})
\overset{D_{X}}{\lto} H^1(X,\mathbb{Z}) 
\overset{j}{\lto} H^1(X,\mathbb{Z})/\pi^{*}H^1(Y,\mathbb{Z}).
\]
One has for every $\lambda,\delta\in\Lambda$
\[
\langle \delta ,\varphi'(\lambda)\rangle = \langle \delta 
,D_{X}(\lambda)\rangle = (\lambda,\delta)_{X} = -mE(\lambda,\delta) = 
\langle \delta, -m\varphi(\lambda)\rangle.
\]
Hence $\varphi'=-m\varphi$ and $(\varphi')^{*}(Q) = 
m^{2}\varphi^{*}(Q) = md_{g}(\; ,\; )_{X}|_{\Lambda}$. Comparing the 
equality $D_{X}^{*}\circ j^{*}(Q) = (\varphi')^{*}(Q) = 
md_{g}(\; ,\; )_{X}|_{\Lambda}
$ with $D_{X}^{*}(Q^{-}) = (\; ,\; )_{X}|_{\Lambda}$ we conclude that 
$j^{*}(Q) = md_{g}(Q^{-})$.
\hfill $\Box$

\bigskip
\noindent
We conclude from the lemma that via the isomorphism $j:M'_{\mathbb{C}} 
\overset{\sim}{\lto} M_{\mathbb{C}}$ the  mapping 
$\tilde{\varPhi}: S\to \mathbb{D}\subset \check{\mathbb{D}}$ 
associated with the polarized VHS\quad $(\mathbb{M},\mathbb{F},Q)$ may 
be identified with the  mapping associated with 
$(\mathbb{M}',\mathbb{F}',Q^{-})$. Via this identification 
$\check{\mathbb{D}}$ is the simplectic grassmanian of $g$-planes in 
$H^1(X,\mathbb{C})^{-}$ isotropic with respect to 
$Q^{-}_{\mathbb{C}}(\omega_{1},\omega_{2}) 
= \int_{X}\omega_{1}\wedge \omega_{2}$ and $\tilde{\varPhi}(s) = 
H^0(\mathcal{X}_{s},\Omega_{\mathcal{X}_{s}})^{-}$. 
\begin{pro}\label{s4.75f}
Let $p:\mathcal{X} \to \mathcal{Y}$ be a surjective holomorphic 
mapping of smooth families of curves over $T$ as in \eqref{es4.73}. Let 
$(\mathbb{M},\mathbb{F},Q)$ be the polarized variation of Hodge 
structures defined in (\ref{s4.73}). Let $s_{0}\in T,\; 
X=\mathcal{X}_{s_{0}}\; Y=\mathcal{Y}_{s_{0}},\; \pi = p_{s_{0}}$ and 
let $H^0(X,\Omega_{X})^{-} = Ker(Tr_{\pi}:H^0(X,\Omega_{X})\to 
H^0(Y,\Omega_{Y}))$. Let $S$ be  a small neighborhood of $s_{0}$ as in 
(\ref{s4.75a}) and let $\tilde{\varPhi}:S\to \mathbb{D}\subset 
\check{\mathbb{D}}$ be the local period mapping. Let 
$\rho:T_{s_{0}}S\to H^1(X,T_{X})$ be the Kodaira-Spencer mapping. Then 
the differential of $\tilde{\varPhi}$ at $s_{0}$ is given by the 
following formula, where $(\; ,\; )$ is the Serre duality pairing (in the next formula $\omega_{1}\cdot \omega_{2}$ denotes the product in the 
symmetric algebra, while $\omega_{1}\omega_{2}$ is a quadratic 
differential).

\begin{equation*}
^{t}d\tilde{\varPhi}(s_{0}): Sym^{2}H^0(X,\Omega_{X})^{-} \lto 
(T_{s_{0}}S)^{*}, 
\quad \langle \frac{\partial}{\partial \tau},\, 
^{t}d\tilde{\varPhi}(\omega_{1}\cdot \omega_{2})
\rangle = 
(\rho(\frac{\partial}{\partial \tau})\; ,\;
\omega_{1}\omega_{2})
\end{equation*}
\end{pro}
{\bf Proof.}
Let $\frac{\partial}{\partial \tau}\in T_{s_{0}}S$ and let $Z\subset 
S$ be  a smooth, complex analytic curve tangent to $\frac{\partial}{\partial 
\tau}$ at $s_{0}$. The restricted family $\mathcal{X}_{Z}\to Z$ yields 
the standard VHS\quad $H^{1,0}(\mathcal{X}_{u})\subset 
H^1(X,\mathbb{C}),\; u\in Z$. 
 By \cite{Gr}~p.816 the 
differential of the period mapping of this VHS transforms 
$\frac{\partial}{\partial 
\tau}$ into an element of $Hom(H^{1,0}(X),H^{0,1}(X))$ given by the 
cup-product $\omega\mapsto 
\rho(\frac{\partial}{\partial \tau})\circ \omega$ (i.e. $\cup$ 
composed with contraction). 
From the splitting
\[
H^{1,0}(\mathcal{X}_{u})\ =\ 
H^{1,0}(\mathcal{X}_{u})^{-}\oplus 
\pi^{*}H^{1,0}(\mathcal{Y}_{u})\subset H^{1}(X,\mathbb{C})^{-}\oplus 
\pi^{*}H^{1}(Y,\mathbb{C})
\]
one concludes that $d\tilde{\varPhi}(\frac{\partial}{\partial \tau}) = 
\phi \in Hom(H^{1,0}(X)^{-},H^{0,1}(X)^{-})$ is given as well by the cup 
product: $\phi(\omega)=\rho(\frac{\partial}{\partial \tau})\circ 
\omega$. 
Applying \eqref{es4.75c} we have that for $\forall 
\omega_{1},\omega_{2}\in H^0(X,\omega_{X})$ 
\[
\langle \frac{\partial}{\partial \tau},\, 
{^{t}d\tilde{\varPhi}}(\omega_{1}\cdot \omega_{2})\rangle \ =\
Q^{-}(\varphi(\omega_{1}),\omega_{2})\ =\ 
\int_{X}\left(\rho(\frac{\partial}{\partial \tau})\circ 
\omega_{1}\right)\wedge \omega_{2}.
\]
We may replace in this integral 
$\rho(\frac{\partial}{\partial \tau})\circ \omega_{1}$ by 
$v^{0,1}\wedge \omega_{1}$, where $v^{0,1}\in 
H_{\overline{\partial}}^{1}(T_{X})$ is a Dolbeault representative of 
$\rho(\frac{\partial}{\partial \tau})$.
The  new integrand $(v^{0,1}\wedge \omega_{1})\wedge \omega_{2}$
is a $(1,1)$ form which is given 
locally by 
\[
((a\frac{\partial}{\partial z}\otimes d\overline{z})\wedge b\, dz)\wedge 
c\, dz\ =\ ab\, d\overline{z} \wedge c\, dz\ =\ abc\, 
d\overline{z}\wedge dz.
\]
Locally $\omega_{1}\omega_{2} = bc(dz)^{2}$, so 
$(v^{0,1}\wedge \omega_{1})\wedge \omega_{2}$
 equals 
$v^{0,1} \wedge
(\omega_{1}\omega_{2})$. Integrating this $(1,1)$ form we obtain by 
the definition of Serre's duality 
 the required formula for $^{t}d\tilde{\varPhi}$.
\hfill $\Box$

\section{Local study of the Prym mapping}\label{s2}

\begin{block}\label{s2.0}
Let $Y$ be a smooth, projective curve of genus $\geq 0$ and let $\pi :X\to Y$ 
be a simple 
covering of degree $d\geq 2$ branched in $n$ points 
$B=\{b_1,\ldots,b_n\}$. Let $y_{0}\in Y-B$ and let $\Delta_{i}$ be 
a small open disk centered at $b_{i},\quad i=1,\ldots,n$. 
Fixing the monodromy 
$m:\pi_{1}(Y-\bigcup_{i=1}^{n}\Delta_{i},y_{0})\to S_{d}$ and varying 
the branch points in $\Delta_{i}$ one obtains a family of $d$-sheeted 
coverings $(\varPsi,f):\mathcal{X}\to Y\times H$ where $H = 
\Delta_{1}\times \cdots \times \Delta_{n}$. Let $w_{i}$ be local 
coordinates of $Y$ in $\Delta_{i}$ satisfying $w_{i}(b_{i})=0$. We 
define coordinates $\mathbf{t} = (t_{1},\ldots,t_{n})$ in $H$ by 
$t_{i}(y_{1},\ldots,y_{i},\ldots, y_{n}) = w_{i}(y_{i})$. At points of 
ramification which project to $\Delta_{i}$ the mapping $\varPsi$ is 
given by $w_{i} = \varPsi_{i}(z_{i},\mathbf{t}) = z_{i}^{2}+t_{i}$ and 
$f(z_{i},\mathbf{t}) = \mathbf{t}$.
\end{block}
\begin{block}\label{s2.0a}
Suppose first $g(Y)=1$. Choose a point $c\notin 
\{b_{1},\ldots,b_{n},y_{0}\}$ and a small disk $D$ centered at $c$. 
Let $\pi^{-1}(D) = D_{1}\cup \ldots \cup D_{d}$ be a disjoint union of 
disks biholomorphically equivalent to $D$. Let $x_{i} = 
\pi^{-1}(c)\cap D_{i}$. Let $v$ be a local coordinate in $D$ with 
$v(c)=0$. Let $u_{i} = \pi^{*}(v)|_{D_{i}}$ be the corresponding local 
coordinates in $D_{i}$. We consider a Schiffer variation of $Y$ (we 
follow here \cite{ACGH}~Vol.~II). Namely paste together $A\times \Delta := 
\{\zeta : |\zeta| \leq 
\epsilon \}\times \{s : |s| < r\}$ with $(Y-\{v : |v|\leq 
\frac{\epsilon}{2}\})\times \Delta$ by means of $\zeta = v + 
\frac{s}{v}$. Here $0< \epsilon \ll 1$ and $r$ depends on $\epsilon$. 
One obtains a family $q:\mathcal{Y}\to N=\Delta$ with 
$\mathcal{Y}_{o}\cong Y$. Performing the corresponding pastings by 
means of 
\begin{equation}\label{es2.0a1}
\zeta_{i} = u_{i} + \frac{s}{u_{i}}\qquad \text{at all}\qquad D_{i}
\end{equation}
and varying the branch points as in the preceding paragraph one 
obtains a smooth proper mapping $f:\mathcal{X}\to N \times H$ and a 
holomorphic mapping\; $p:\mathcal{X}\to \mathcal{Y}$ which fit 
into a commutative diagram 
\begin{equation}\label{es2.0a}
\begin{diagram}
\mathcal{X} & \rTo^p         & \mathcal{Y} \\
\dTo^f      &                &  \dTo_q     \\
N\times H   & \rTo^{\pi_{1}} & N           \\
\end{diagram}
\end{equation}
If $g(Y)\geq 2$ one chooses $m = 3g(Y)-3$ general points 
$c_{1},\ldots,c_{m}$ on $Y$ which impose independent conditions on 
$H^0(Y,\omega_{Y}^{\otimes 2})$, one repeats the same pasting 
procedure simultaneously at $c_{i},\; i=1,\ldots,m$ and obtains a 
 family $q:\mathcal{Y}\to N=\Delta^{m}$ with $\mathcal{Y}_{o}\cong Y$.\
 In both cases, $g(Y) = 1$ or $g(Y)\geq 2$,\quad $q:\mathcal{Y}\to 
N=\Delta^{m}$ is a minimal versal deformation of $Y$ since the Kodaira-Spencer 
mapping $\rho :T_oN\to H^1(Y,T_Y)$ is an isomorphism \cite{KS}.   Similarly to 
the 
above one obtains a holomorphic mapping\; $p:\mathcal{X}\to 
\mathcal{Y}$ and a smooth, proper mapping $f:\mathcal{X}\to N\times H$ 
which fit into a commutative diagram as in \eqref{es2.0a}.
The constructed family of coverings is versal for deformations of 
$\pi :X\to Y$ as shows the following proposition.
\end{block}

\begin{pro}\label{s2.0b}
Let $\pi :X\to Y$ be a covering of a smooth, projective curve of genus 
$g(Y)\geq 1$ simply branched in $B\subset Y$. Let 
\[
\begin{diagram}
\mathcal{X}'& \rTo^{p'} & \mathcal{Y}' \\
\dTo^{f'}       &      &   \dTo_{q'}\\
T          & \rTo^{h} &  Z      \\
\end{diagram}
\]
be a commutative diagram of holomorphic mappings of complex manifolds, 
where $f'$ and $q'$ are proper, smooth and surjective, $p'$ is finite 
and surjective, and there is a point $s_{0}\in T$ such that the covering 
$p'_{s_{0}}:\mathcal{X}'_{s_{0}}\to \mathcal{Y}'_{h(s_{0})}$ is 
isomorphic to $\pi :X\to Y$. Then  there exist neighborhoods $S$ and $V$ 
with $s_{0}\in S\subset T$,\; $h(s_{0})\in V\subset Z$, $h(S)\subset 
V$ and holomorphic mappings $\mu$ and $\nu$ which fit into the 
commutative diagram 

\begin{diagram}
S           & \rTo^h         & V         \\
\dTo^{\mu}  &                &  \dTo_{\nu}     \\
N\times H   & \rTo^{\pi_{1}} & N           \\
\end{diagram}

such that the restricted family of coverings $p':\mathcal{X}'_{S}\to 
\mathcal{Y}'_{V}$ is the pull-back via $\mu, \nu$ of the covering 
\eqref{es2.0a}
\end{pro}
{\bf Proof.}
As mentioned above the existence of $V\ni h(s_{0})$ and a holomorphic 
mapping $\nu:V\to N$ such that $\mathcal{Y}'_{V}\cong 
\mathcal{Y}\times_{N}V$ follows from Kodaira-Spencer's theorem of 
completeness \cite{KS}. We may replace $T$ by $h^{-1}(V)$ and 
$f':\mathcal{X}'\to T$ by the corresponding restriction. We then 
obtain a deformation into the family $q:\mathcal{Y}\to N$ as defined 
in \cite{Ho2}~\S~5:
\begin{equation}\label{es2.0d}
\begin{diagram}
\mathcal{X}'& \rTo     & \mathcal{Y} \\
\dTo^{f'}       &      &   \dTo_{q}\\
T          & \rTo^{\nu \circ h}  &  N   \\
\end{diagram}
\end{equation}
An easy calculation which uses \cite{Ho2}~Lemma~5.1 shows that the 
characteristic map \linebreak
$\tau: T_{(o,B)}N\times H \to D_{X/\mathcal{Y}}$ of 
the deformation \eqref{es2.0a} is an isomorphism. Therefore by 
Horikawa's theorem of completeness \cite{Ho2}~Theorem~5.2 there is a 
neighborhood $S$ of $s_{0}\in T$ and a holomorphic mapping $\mu:S\to 
N\times H$ such that the restriction of \eqref{es2.0d} on $S$ is 
isomorphic to the pull-back of \eqref{es2.0a} by $\mu$. This proves 
the proposition.
\hfill $\Box$

\bigskip
\noindent
Let $Y$ be a curve of genus $g\geq 0$ and let $\pi :X\to Y$ be a 
simple covering of degree $d\geq 2$ with $n\geq 2$ branch points. Let 
$\mathcal{C}\to M,\quad  \mathcal{C}_{o} \cong X$ be a
 minimal versal deformation of $X$.
Let $(\varPsi,f): \mathcal{X}\to Y\times H$ be the (local 
Hurwitz) family considered in (\ref{s2.0}). 
Shrinking $H$ if necessary,
there is a  holomorphic mapping $h:H\to M$ such 
that $f:\mathcal{X}\to H$ is the pull-back of $\mathcal{C}\to M$ via 
$h$.
\begin{pro}\label{s2.1}
Let $\pi :X\to Y$ be as above, let $B  = 
\{b_{1},\ldots,b_{n}\}\subset Y$ be the branch locus and let $R = 
\{a_{1},\ldots,a_{n}\}\subset X,\quad \pi(a_{i})=b_{i}$ be the 
ramification locus. Let $w_{i},\; w_{i}(b_{i})=0$ be local coordinates 
of $Y$ at $b_i,\; i=1,\ldots,n$ and let $\mathbf{t} = (t_{1},\ldots,t_{n})$ 
be the 
corresponding coordinates of $H$ as defined in (\ref{s2.0}). 
Then the transpose of the differential of \; $h$ 
\[
^{t}(d h)_{o} : H^0(X,\omega_{X}^{2})\cong (T_{o}M)^{*}\lto 
(T_{B}H)^{*} \cong \bigoplus_{i=1}^{n}\mathbb{C}(d\, t_{i})_{o}
\]
is given by the formula
\begin{equation}\label{es2.1}
\varphi \quad \mapsto \quad 
2\pi \sqrt{-1} 
\sum_{i=1}^{n}Res_{a_{i}}\frac{\varphi}{\pi^{*}dw_{i}}\; (dt_{i})_{o}
\end{equation}
\end{pro}
{\bf Proof.}
We are in the situation considered by Horikawa in \cite{Ho1}. At the 
ramification points of $(\varPsi,f):\mathcal{X}\to Y\times H$ the 
mapping $\varPsi :\mathcal{X}\to Y$ is given by (cf. (\ref{s2.0})) 
$w_i = \varPsi_{i}(z_{i},\mathbf{t}) = z_{i}^{2}+t_{i}$. Consider as in 
\cite{Ho1} the sequence
\[
0\lto T_{X} \overset{d\pi}{\lto} \pi^{*}T_{Y}\lto \mathcal{T} \lto 0.
\]
Here $\mathcal{T}$ is a skyscraper sheaf (non-canonically) isomorphic 
to $\oplus_{i=1}^{n}\mathbb{C}_{a_{i}}$. Then Horikawa's 
characteristic mapping $\tau :T_{B}H\to H^0(X,\mathcal{T})$ is given 
by 
\[
\tau (\frac{\partial}{\partial t_{i}})\quad = \quad 
\frac{\partial}{\partial w_{i}}\; mod\, d\pi(T_{X})_{a_{i}} \; \in\;
\mathbb{C}_{a_{i}}\; \subset \; \oplus_{k=1}^{n}\mathbb{C}_{a_{k}}.
\]
Consider the exact cohomology sequence
\[
0\lto H^0(X,T_{X})\lto H^0(X,\pi^{*}T_{Y})\lto H^0(X,\mathcal{T})
\overset{\delta}{\lto} H^{1}(X,T_{X}).
\]
The Kodaira-Spencer mapping equals the composition $\rho = \delta 
\circ \tau$\; \cite{Ho1}~p.376. In order to prove \eqref{es2.1} it 
suffices to verify that for every $\varphi\in 
H^0(X,\omega_{X}^{\otimes \,2})$ and every $i=1,\ldots,n$ for the 
Serre duality pairing one has
\begin{equation}\label{es2.3}
( \rho(\frac{\partial}{\partial t_{i}})\; ,\; \varphi ) \  = \
2\pi \sqrt{-1}\; Res_{a_{i}}
(\frac{\varphi}{\pi^{*}dw_{i}}
).
\end{equation}
For every $i=1,\ldots,n$ let $U_{i}$ be a small disc containing 
$a_{i}$ with local coordinate $z_{i},\; z_{i}(a_{i})=0$, such that 
$\pi(U_{i})\subset \Delta_{i}$ and such that $\pi $ is locally given 
by $w_{i}=z_{i}^{2}$. Let $U_{0}=X-R$. Consider the Stein covering 
$\mathfrak{U} = \{U_{0},U_{1},\ldots,U_{n}\}$. By the definition of 
$\tau(\frac{\partial}{\partial 
t_{i}})$ it is immediate that 
$ \rho(\frac{\partial}{\partial t_{i}}) = 
\delta( \tau(\frac{\partial}{\partial t_{i}}))$ is given by 
a 1-cochain $\{\xi_{\alpha \beta}\}\in C^{1}(\mathfrak{U},T_{X})$, 
such that $\xi_{i0} = \frac{1}{2z_{i}}\, \frac{\partial}{\partial 
z_{i}} = -\xi_{0i}$ while the other $\xi_{\alpha \beta} = 0$. Let 
$\varphi = f_{i}(z_{i})(d\, z_{i})^{2}$ in $U_{i}$. Then $\xi_{i0}\cdot 
\varphi = \frac{f_{i}(z_{i})}{2z_{i}}d\, z_{i} = 
\frac{\varphi}{\pi^{*}d\, w_{i}}$. A calculation similar to the one in 
\cite{ACGH}~pp.14,15 shows \eqref{es2.3}.
\hfill $\Box$

\begin{cor}\label{s2.5}
Let the hypothesis be as in Proposition~\ref{s2.1}. Then the 
annihilator of\; $d h (T_{B}H)$ in $(T_{o}M)^{*}$ equals 
$H^{0}(X,\omega_{X}^{\otimes 2}(-R))$
\end{cor}
\begin{rem}\label{s2.5a}
When $Y=\mathbb{P}^{1}$ this result was proved  by Donagi 
and Green by a different argument \cite{DS}~Appendix.
\end{rem}

\smallskip
\noindent
Let us consider now the case 
of a simple branched covering $\pi :X\to Y$ of degree $d$ where
 $g(Y)\geq 1$. Let $q:\mathcal{Y}\to N$ be the minimal versal deformation of 
$Y$ 
 constructed by Schiffer variations at 
non-branched points of $Y$ (cf. (\ref{s2.0a})) and let us consider the 
deformation \eqref{es2.0a}.
Let $\mathcal{C}\to 
M,\quad \mathcal{C}_{o}\cong X$ be a minimal versal deformation of $X$.
Shrinking $N$ and $H$ if necessary,
there exists a  holomorphic mapping $h:N\times H\to M$ such 
that $f:\mathcal{X}\to N\times H$ is the pull-back of $\mathcal{C}\to M$ via 
$h$.
\begin{pro}\label{s2.6}
Let $\pi :X\to Y$ be as above, let $B  = 
\{b_{1},\ldots,b_{n}\}\subset Y$ be the branch locus and let $R = 
\{a_{1},\ldots,a_{n}\}\subset X,\quad \pi(a_{i})=b_{i}$ be the 
ramification locus. Let $m=1$  if  $g(Y)=1$  and  let  $m=3g(Y)-3$  if 
$g(Y)\geq 2$. Let $c_{1},\ldots,c_{m}$ be general points  in  $Y$  and 
let $\pi^{-1}(c_{i}) = \{x_{i1},\ldots,x_{id}\}$. Let $v_{i}$
 be local coordinates at $c_{i},\quad i=1,\ldots,m$ and 
let $w_{j}$ be local coordinates at $b_{j},\quad j=1,\ldots,n$. Then 
the transpose of the differential of \; $h$
\[
{^{t}dh}\ :\ H^{0}(X,\omega_{X}^{\otimes 2})\lto 
(T_{(o,B)}N\times M)^*\; =\;
\left(\oplus_{i=1}^{m}\mathbb{C}(d\, s_{i})_{o}\right) \oplus 
\left(\oplus_{j=1}^{n}\mathbb{C}(d\, t_{j})_{o}\right)
\]
is determined by the following formulas 
\begin{align}
\langle \frac{\partial}{\partial s_{i}}\; ,\; {^{t}dh}(\varphi)\rangle 
\quad & = \quad 
2\pi \sqrt{-1} \:
\sum_{k=1}^{d}\frac{\varphi}{\pi^{*}(dv_{i})^{2}}(x_{ik})\\
\langle \frac{\partial}{\partial t_{j}}\; ,\; {^{t}dh}(\varphi)\rangle 
\quad & = \quad 
2\pi \sqrt{-1} \:
Res_{a_{j}}\frac{\varphi}{\pi^{*}dw_{j}}
\end{align}
\end{pro}
{\bf Proof.}
The second formula was proved in Proposition~\ref{s2.1}. We give the 
proof of the first one in the case $g(Y)=1$. The case $g(Y)\geq 2$ 
is similar. Using the notation of (\ref{s2.0a}) we consider the Stein 
covering of $X$,\quad $X=D_{0}\cup D_{1}\cup \ldots \cup D_{d}$, where 
$D_{0} = X-\pi^{-1}(a)$. The Kodaira-Spencer class $[\theta] = 
\rho(\frac{\partial}{\partial s})\in 
H^{1}(X,T_{X})$ of the deformation of $X$ defined by \eqref{es2.0a1}, the 
branch points remaining fixed, is 
given by the 1-cocycle 
$\theta_{i0}=-\theta_{0i}=\frac{1}{u_{i}}\frac{\partial}{\partial 
u_{i}} \in H^0(D_{i}\cap D_{0},T_{X})$ for $i\geq 1$ and 
$\theta_{ij} = 0$ for $i,j\geq 1$. Let $\varphi\in 
H^{0}(X,\omega_{X}^{\otimes 2})$ and let 
$\varphi_{i}=f_{i}(du_{i})^{2}$ in $D_{i}$. Then 
$\theta_{i0}\cdot \varphi = \frac{f_{i}(u_{i})}{u_{i}}du_{i}$. Thus 
$Res_{x_{i}}(\theta_{i0}\cdot \varphi) = f_{i}(0) = 
\frac{\varphi}{\pi^{*}(dv)^{2}}(x_{i})$. Using  a calculation 
similar to the one in \cite{ACGH}~pp.14,15 one obtains for the Serre 
duality pairing
\[
( [\theta]\; ,\; \varphi )\;   =\; 
2\pi \sqrt{-1} \:
\sum_{k=1}^{d}\frac{\varphi}{\pi^{*}(dv)^{2}}(x_{k})
\]
\hfill $\Box$

\begin{block}\label{s2.7}
Let $\pi :X\to Y$ be  a simple branched covering of degree $d\geq 2$ 
where $Y$ is a smooth projective curve of genus $g(Y)\geq 1$. Consider 
the deformation \eqref{es2.0a} of the covering $\pi$ as described at 
the end of (\ref{s2.0a}). Let $S=N\times H$ and let $q':\mathcal{Y}'\to 
S$ be the pullback of the  family $q:\mathcal{Y}\to N$ via 
$\pi_{1}:S\to N$. We obtain a commutative triangle as in 
\eqref{es4.73}\; :
\begin{equation*}
\begin{diagram}
\mathcal{X}&    &\rTo^{p'} &       &\mathcal{Y}' \\
           &\rdTo_f &     & \ldTo_{q'}  &       \\
           &      &  S  &        &           \\
\end{diagram}
\end{equation*}
By (\ref{s4.73}) one associates to it a polarized VHS\quad 
$(\mathbb{M},\mathbb{F},Q)$. Shrinking $N$ and $H$ if necessary one 
may define a lifting of the period mapping $\tilde{\varPhi}:S\to 
\mathbb{D}\subset \check{\mathbb{D}}$ as in (\ref{s4.75a}) and furthermore 
this 
mapping may be identified with the local period mapping associated with the 
polarized VHS \quad $(\mathbb{M}',\mathbb{F}',Q^{-})$ as in (\ref{s4.75c}) and 
Lemma~\ref{s4.75d}. Let $s_{0} = 
(o,B)$ be the reference point of $S=N\times H$.
\end{block}
\begin{pro}\label{s2.9}
Using the notation of Proposition~\ref{s2.6} one has the following 
formula for the differential of $\tilde{\varPhi}$ at $s_{0}$:
\[
{^{t}d\tilde{\varPhi}}(s_{0}): Sym^{2}H^0(X,\omega_{X})^{-} \lto 
\left(\oplus_{i=1}^{m}\mathbb{C}(ds_{i})_{o}\right) \oplus 
\left(\oplus_{j=1}^{n}\mathbb{C}(dt_{j})_{o}\right).
\]
For every $\omega_{1},\omega_{2}\in H^{0}(X,\omega_{X})^{-}$ it holds
\begin{align}
\langle \frac{\partial}{\partial s_{i}}\; ,\; 
{^{t}d\tilde{\varPhi}}(\omega_{1}\cdot \omega_{2})\rangle 
\quad & = \quad 
2\pi \sqrt{-1} \:
\sum_{k=1}^{d}\frac{\omega_{1}\omega_{2}}{\pi^{*}(dv_{i})^{2}}
(x_{ik})\label{es2.9a}\\
\langle \frac{\partial}{\partial t_{j}}\; ,\; 
{^{t}d\tilde{\varPhi}}(\omega_{1}\cdot \omega_{2})\rangle 
\quad & = \quad 
2\pi \sqrt{-1} \:
Res_{a_{j}}\frac{\omega_{1}\omega_{2}}{\pi^{*}dw_{j}}\label{es2.9b}
\end{align}
\end{pro}
{\bf Proof.}
One applies Proposition~\ref{s4.75f} with $\frac{\partial}{\partial 
\tau} = \frac{\partial}{\partial s_{i}}$ or $\frac{\partial}{\partial 
t_{j}}$. That 
$(\rho(\frac{\partial}{\partial \tau})\; ,\; 
\omega_{1}\omega_{2})$ equals the right-hand side of \eqref{es2.9a} 
and \eqref{es2.9b} in the respective cases follows from the proofs of 
Proposition~\ref{s2.6} and Proposition~\ref{s2.1}.
\hfill $\Box$

\begin{block}\label{s2.11}
We now restrict ourselves to the case of  coverings of elliptic curves $\pi 
:X\to Y$ with $g = g(X)\geq 3$ and first we consider 
deformations of $\pi :X\to Y$ with fixed $Y$. Consider the canonical 
map $\phi_{K}:X\to \mathbb{P}^{g-1} = |\omega_{X}|^{*}$. The 
space $H^0(\omega_{X})^{-}  =   = \{\omega \; |\; Tr_{\pi}(\omega) = 0\}$ is a 
hyperplane in $H^{0}(\omega_{X})$ since $g(Y)=1$. It defines a point 
$q^-\in \mathbb{P}^{g-1}$ with the property that the differentials 
$\omega\in H^{0}(\omega_{X})^{-}$ define hyperplanes in 
$\mathbb{P}^{g-1}$ containing $q^-$. Let $\alpha \in 
H^0(Y,\omega_{Y}),\quad \alpha \ne 0$. Let $H_{\alpha }\subset 
\mathbb{P}^{g-1}$ be the hyperplane in $\mathbb{P}^{g-1}$ which 
corresponds to $\pi^{*}\alpha \in H^{0}(\omega_{X})$. Clearly 
$div(\pi^{*}\alpha ) = R = \phi_{K}^{*}(H_{\alpha })$.
\end{block}
\begin{lem}\label{s2.11a}
Let $\pi :X\to Y$ be a covering of an elliptic curve of degree $d\geq 
2$. Then $q^-\notin H_{\alpha }$
\end{lem}
{\bf Proof.}
$Tr_{\pi}(\pi^{*}\alpha ) = d\alpha \ne 0$, hence by definition 
$q^-\notin H_{\alpha }$.
\hfill $\Box$

\bigskip
\noindent
Let $\pi :X\to Y$ be a simple covering of an elliptic curve of degree 
$d\geq 2$ branched at $B\subset Y$ where $g=g(X)\geq 3$. Consider the 
deformation $(\varPsi,f):\mathcal{X}\to Y\times H$ defined in (\ref{s2.0}). 
Let 
$\tilde{\varPhi}_{Y}:H\to \mathbb{D}\cong \mathfrak{H}_{g-1}$ be the 
period mapping corresponding to the variation 
$H^{1,0}(\mathcal{X}_{s})^{-}\subset H^{1}(X,\mathbb{C})^{-},\quad 
s\in H$ (cf. (\ref{s4.75c})). Clearly the mapping $\tilde{\varPhi}_{Y}$ is 
invariant with respect to the action of the germ of $Aut(Y)$ at $id$. 
So $\dim Ker\, d\tilde{\varPhi}_{Y}(s_{0})\geq 1$ where $s_{0}=B\in 
H$.
\begin{pro}\label{s2.12}
Let $\pi :X\to Y$ be a simple covering of an elliptic curve. Suppose $X$ is 
not 
hyperelliptic and $g(X)\geq 3$. Then 
$\dim Ker d\tilde{\varPhi}_Y(s_{0}) = 1$. In particular if $g(X)=3$, 
then $d\tilde{\varPhi}_Y(s_{0})$ is an epimorphism, and if $g(X)=4$ the 
image $d\tilde{\varPhi}_Y(T_{B}H)$ is a hypersurface in 
$T_{\Pi}\mathbb{D}$ where $\Pi=\tilde{\varPhi}(s_{0})$.
\end{pro}
{\bf Proof.}
We identify $X$ with its image  
$\phi_{K}(X)\subset \mathbb{P}^{g-1}$. Then $H_{\alpha 
}\cap X = R$ where $R$ is the ramification divisor of $\pi :X\to Y$ 
and the span $\langle R\rangle = H_{\alpha }$. By \eqref{es2.9b} the 
annihilator of $d\tilde{\varPhi}_Y(T_{B}H)$ equals the subspace of 
$S^{2}H^{0}(\omega_{X})^{-}$ consisting of elements which vanish in 
$R$. Identifying $H^{0}(\omega_{X})$ with 
$H^{0}(\mathbb{P}^{g-1},\mathcal{O}_{\mathbb{P}^{g-1}}(1))$ we 
consider the restriction map $r:S^{2}H^{0}(\omega_{X})^{-}\to 
H^0(H_{\alpha },\mathcal{O}_{H_{\alpha }}(2))$. This is an isomorphism 
since both spaces have the same dimension and the kernel is 0. Indeed 
any element of the kernel corresponds to a reducible quadric $H'\cup 
H_{\alpha }$ with center $H'\cap H_{\alpha }$. The quadrics of 
$|S^{2}H^{0}(\omega_{X})^{-}|$ are singular and contain $q^-$ in their 
center. By Lemma~\ref{s2.11a} $q^-\notin H_{\alpha }$, thus no quadric 
of the type $H'\cup H_{\alpha }$ could belong to 
$|S^{2}H^{0}(\omega_{X})^{-}|$. We conclude that $Ker\,
{^{t}d\tilde{\varPhi}_Y(s_{0})} \cong H^0(H_{\alpha 
},\mathcal{O}_{H_{\alpha }}(2)\otimes J_{R})$ where $J_{R}$ is the 
ideal sheaf of $R\subset H_{\alpha }$.

Let $g(X)=3$. Then $X\subset \mathbb{P}^{2}$ is a quartic and 
$H_{\alpha }$ is a line intersecting $X$ in 4 distinct points. Here 
$
H^0(H_{\alpha 
},\mathcal{O}_{H_{\alpha }}(2)\otimes J_{R}) \cong 
H^0(\mathbb{P}^{1},\mathcal{O}_{\mathbb{P}^{1}}(-2)) = 0$. Thus 
$d\tilde{\varPhi}_Y(s_{0})$ is epimorphic with kernel of dimension 1 
since $\dim H = 4,\; \dim \mathbb{D} = 3$. 

Suppose $g(X)\geq 4$. The restriction mapping 
\[
H^0(\mathbb{P}^{g-1},\mathcal{O}_{\mathbb{P}^{g-1}}(2)\otimes J_{X}) 
\lto
H^0(H_{\alpha },\mathcal{O}_{H_{\alpha }}(2)\otimes J_{R}) 
\]
is an isomorphism (see the proof of \cite{Is}~Lemma~2.10). Therefore 
$h^0(H_{\alpha },\mathcal{O}_{H_{\alpha }}(2)\otimes J_{R}) = 
\frac{(g-2)(g-3)}{2}$. 
We conclude that
\[
\dim d\tilde{\varPhi}_Y(T_{B}H)\quad =\quad \frac{(g-1)g}{2} - 
\frac{(g-2)(g-3)}{2}\quad = \quad 2g-3.
\]
By Hurwitz' formula $2g-2 = n + d(2g(Y)-2) = n = \dim H$. Hence $\dim 
Ker\, d\tilde{\varPhi}_Y(s_{0}) = 1$.
\hfill $\Box$

\begin{rem}\label{s2.13}
If $X$ is hyperelliptic, then similar calculation shows that $\dim 
d\tilde{\varPhi}_Y(T_{B}H) = g-1$. Hence $\dim Ker\, 
d\tilde{\varPhi}_Y(s_{0}) = g-1 \geq 2$ for $g\geq 3$. 
The next proposition shows that if $n\geq 4$  and if $[X\to Y]\in 
\mathcal{H}_{d,n}(Y)$ is sufficiently general, then  the curve $X$ is not 
hyperelliptic.
\end{rem}
\begin{pro}\label{s2.13a}
Let $d\geq 2, n=2e\geq 4$. Then the locus of points $[X\to Y]$ in 
$\mathcal{H}_{d,n}(Y)$ with hyperelliptic $X$ is contained in a closed 
subset of $\mathcal{H}_{d,n}(Y)$ of codimension $\geq \frac{n}{2}-1$
\end{pro}
{\bf Proof.}
Suppose $\pi :X\to Y$ is a covering of degree $d$ simply branched in 
$n$ points such that $X$ is hyperelliptic. Let $\mu:X\to 
\mathbb{P}^{1}$ be the corresponding covering of degree 2. Consider 
$\nu = (\mu,\pi):X\to \mathbb{P}^{1}\times Y$. Let $\nu(X)=X'$. The 
morphism $\nu:X\to X'$ is birational. Indeed, since $\mu=\pi_{1}\circ 
\nu$, the other possibility would be $\deg \nu = 2, X'\cong 
\mathbb{P}^{1}$. This is however absurd, since $\pi_{2} :X'\to Y$ is 
an epimorphism. The curve $X'$ belongs to 
$|\pi_{1}^{*}\mathcal{O}_{\mathbb{P}^{1}}(d)\otimes \pi_{2}^{*}L|$ for 
some invertible sheaf $L\in Pic^{2}Y$. Let $\varphi = \varphi_{L}:Y\to 
\mathbb{P}^{1}$ and let $C=(\mathbf{1}\times \varphi)(X')\subset 
\mathbb{P}^{1}\times \mathbb{P}^{1}$. The curve $C$ is irreducible and 
belongs to $|\mathcal{O}_{\mathbb{P}^{1}\times \mathbb{P}^{1}}(d,1)|$, 
hence it is isomorphic to $\mathbb{P}^{1}$. We obtain that the 
covering $\pi :X\to Y$ fits into a commutative diagram
\begin{equation}\label{es2.13b}
\begin{diagram}
X & \rTo^{\nu} & X'   & \rTo  &  Y \\
  & \rdTo_{\mu}& \dTo &       & \dTo_{\varphi} \\
  &            & C    & \rTo^{\psi}& \mathbb{P}^{1}\\
\end{diagram}
\end{equation}
where $\psi:C\to \mathbb{P}^{1}$ is a morphism of degree $d$. 
Conversely, given $\varphi:Y\to \mathbb{P}^{1}$ of degree $2$ and 
$\psi:C\to\mathbb{P}^{1}$ of degree $d$ one defines 
$X'=C\times_{\mathbb{P}^{1}}Y$. Normalizing $X'$ one obtains 
$\pi :X\to Y$ such that $X$ is hyperelliptic. Composing in 
\eqref{es2.13b} with an automorphism $\tau:\mathbb{P}^{1}\to 
\mathbb{P}^{1}$ does not change the equivalence class of $\pi :X\to 
Y$. Thus we may suppose the discriminant locus of $\varphi$ is a fixed 
set $B_{\varphi}=\{0,1,\lambda,\infty\}$ where $\lambda$ is determined 
by the $j$-invariant of $Y$ up to a finite number of choices 
\cite{Ha}~Ch.IV\S~4. If $c\in \mathbb{P}^{1}-B_{\varphi}$ is a 
branch point of $\psi$, then the two points $\varphi^{-1}(c)$ are 
branch points of $\pi$. One concludes firstly that all branch points 
of $\psi$ contained in $\mathbb{P}^{1}-B_{\varphi}$ are simple and 
secondly $\# (B_{\psi}-B_{\varphi})\leq \frac{n}{2}$. Families of 
coverings $\pi :X\to Y$ that fit into the diagram \eqref{es2.13b} are 
constructed by: (a) considering Hurwitz spaces of coverings $\psi:C\to 
\mathbb{P}^{1}$ where $\deg \psi =d,\; C\cong \mathbb{P}^{1}$, and the 
branch points of $\psi $ belonging to $\mathbb{P}^{1}-B_{\varphi}$ are simple, 
(b) composing the morphisms in the diagram \eqref{es2.13b} with an 
automorphism of $Y$. The dimensions of the possible Hurwitz spaces of 
(a) are $\leq \frac{n}{2}$. Hence the equivalence classes of coverings 
$[X\to Y]$ with hyperelliptic $X$ are contained in a closed subset of 
$\mathcal{H}_{d,n}(Y)$ of dimension $\leq \frac{n}{2}+1$.
\hfill $\Box$

\bigskip
\noindent
The next lemma  follows from  the proof of 
Proposition~\ref{s2.12}. Abusing notation we will not distinguish 
between $S^{m}H^{0}(X,\omega_{X})$ and 
$H^{0}(\mathbb{P}^{g-1},\mathcal{O}_{\mathbb{P}^{g-1}}(m))$. We denote 
by $S^{2}H^{0}(X,\omega_{X})^{^{-}}(-R)$ the space of elements of 
$S^{2}H^{0}(X,\omega_{X})^{-}$ which vanish on the ramification locus 
$R$.
\begin{lem}\label{s2.14}
Let the assumptions be  as in Proposition~\ref{s2.12}. 
Let $g(X)\geq 4$. Then for every $F\in 
H^{0}(\mathbb{P}^{g-1},\mathcal{O}_{\mathbb{P}^{g-1}}(2)\otimes 
J_{X})$ there is a unique element $F^{-}\in 
S^{2}H^{0}(X,\omega_{X})^{^{-}}(-R)$ such that 
\[
F^{-}\quad =\quad \omega \cdot \pi^{*}\alpha + F
\]
for some $\omega\in H^{0}(X,\omega_{X})$. The mapping $F\mapsto 
F^{-}$ yields a linear isomorphism between 
$H^{0}(\mathbb{P}^{g-1},\mathcal{O}_{\mathbb{P}^{g-1}}(2)\otimes 
J_{X})$ and $S^{2}H^{0}(X,\omega_{X})^{^{-}}(-R)$.
\end{lem}
We now consider the  deformation \eqref{es2.0a} of $\pi :X\to Y$. Let 
$\tilde{\varPhi}:N\times H\to \mathbb{D} $ be the period mapping 
defined in (\ref{s2.7}). We want to calculate $Ker\, 
d\tilde{\varPhi}(s_{0})$ where $s_{0}=(o,B)$. Since 
$\tilde{\varPhi}|_{\{0\}\times H} = \tilde{\varPhi}_{Y}$ this kernel 
has dimension $\geq 1$. 
\begin{pro}\label{s2.15}
Let $\pi :X\to Y$ be a simple  covering of an elliptic curve.
Suppose $X$ is not hyperelliptic and $g(X)\geq 4$.  Then 
$\dim Ker\, d\tilde{\varPhi}(s_{0}) = 1$ if and only if the point 
$q^-\in 
\mathbb{P}^{g-1}$ which corresponds to $H^{0}(\omega_{X})^{-}$ does not 
belong to the intersection of quadrics which contain $\phi_{K}(X)$.
\end{pro}
{\bf Proof.}
Identify $X$ with 
$\phi_{K}(X)\subset \mathbb{P}^{g-1}$.
Let $c\in Y$ and let $v$ be a local coordinate as in (\ref{s2.0a}). Let 
$\pi^{-1}(c)=\{x_{1},\ldots,x_{d}\}$. Let 
$m:S^{2}H^{0}(X,\omega_{X})\to H^{0}(X,\omega_{X}^{\otimes 2})$ be 
the multiplication map. By Proposition~\ref{s2.12} and 
Proposition~\ref{s2.9} one 
has $\dim Ker\, d\tilde{\varPhi}(s_{0}) = 1$ if and only if the linear 
functional $\lambda$  on $S^{2}H^{0}(X,\omega_{X})^{^{-}}(-R)$ defined 
by 
\begin{equation}\label{es2.16}
\lambda(F^{-})\quad = \quad 
\sum_{k=1}^{d}\frac{m(F^{-})}{\pi^{*}(dv)^{2}}(x_{k})
\end{equation}
is not identically zero. For every $F^{-}$ we have by 
Lemma~\ref{s2.14} the representation $F^{-}=\omega\cdot
\pi^{*}\alpha + F$ where $F\in 
H^{0}(\mathbb{P}^{g-1},\mathcal{O}_{\mathbb{P}^{g-1}}(2)\otimes 
J_{X})$. It is clear that 
\[
\lambda(F^{-})\quad =\quad \frac{Tr_{\pi}(\omega)}{dv}(c)\:
\frac{\alpha }{dv}(c).
\]
Since $Tr_{\pi}(\omega)=const\, \alpha $ we conclude that $\lambda 
(F^{-}) = 0$ if and only if $Tr_{\pi}(\omega) = 0$, i.e. iff $\omega\in 
H^{0}(X,\omega_{X})^{-}$, or equivalently iff the linear form of 
$H^{0}(\mathbb{P}^{g-1},\mathcal{O}_{\mathbb{P}^{g-1}}(1))$ 
corresponding to $\omega$ vanishes in $q^-$. 

Suppose there is a quadric $Q\subset \mathbb{P}^{g-1}$ given by an 
equation $F = 0$ such that $Q\supset X$ and $q^-\notin Q$. Let 
$F^{-}=\omega\cdot \pi^{*}\alpha + F$ be the element in 
$S^{2}H^{0}(X,\omega_{X})^{^{-}}(-R)$ corresponding to $F$ as in 
Lemma~\ref{s2.14}. Since $q^-$ belongs to the center of the quadric 
$\{F^{-}=0\}$ the form $\omega$ cannot belong to 
$H^{0}(X,\omega_{X})^{-}$, thus $\lambda(F^{-})\ne 0$.

Conversely, let $\lambda(F^{-})\ne 0$ for some $F^{-}\in 
S^{2}H^{0}(\omega_{X})^{^{-}}(-R)$. Let $F^{-}=\omega\cdot \pi^{*}\alpha 
+ F$ be the representation of Lemma~\ref{s2.14}. By Lemma~\ref{s2.11a} 
and the argument above the linear forms of 
$H^{0}(\mathbb{P}^{g-1},\mathcal{O}_{\mathbb{P}^{g-1}}(1))$ which 
correspond to $\omega$ and $\pi^{*}\alpha $ do not vanish in $q^-$. 
Since $q^-$ belongs to the center of $\{F^{-}=0\}$ we conclude that $F$ 
does not vanish in $q^-$.
\hfill $\Box$

\begin{cor}\label{s2.18}
Let $X\subset \mathbb{P}^{3}$ be a canonical curve which is a simple 
covering    of    an   elliptic 
curve $\pi  :X\to  Y$. Let 
$Q$ be the unique quadric in $\mathbb{P}^{3}$ which contains  $X$  and 
let $q^- \in \mathbb{P}^{3}$ be the point corresponding to the hyperplane 
of holomorphic differentials with trace  0.  Let  $\tilde{\varPhi}:N\times 
H\to \mathbb{D}\cong \mathfrak{H}_{3}$ be the period  mapping  defined 
in (\ref{s2.7}). Then the differential $d\tilde{\varPhi}(s_{0})$
at the point corresponding to $\pi :X\to Y$ is surjective if  and  only 
if $q^-\notin Q$.
\end{cor}
\begin{pro}\label{s2.18a}
Let $\pi :X\to Y$ be a covering of an elliptic curve of degree $d=2$. 
Suppose $X$ is not hyperelliptic and $g(X)\geq 4$. Let 
$\tilde{\varPhi}:N\times H\to \mathbb{D}\cong \mathfrak{H}_{g-1}$ be 
the mapping of (\ref{s2.7}) obtained by the periods of the Prym 
differentials. Then $\dim Ker\, d\tilde{\varPhi}(s_{0}) = 1$. In 
particular if $g(X) = 4$ then $d\tilde{\varPhi}(s_{0})$ is surjective.
\end{pro}
{\bf Proof.}
$X$ is a bi-elliptic curve and it is well-known that $\phi_{K}(X)$ 
lies on a normal elliptic cone $C$ of degree $g-1$ (cf. \cite{ACGH}~p.269). 
 The vertex of the cone is exactly the point $q^-$ 
defined in (\ref{s2.11}). One has $q^- \notin X$, the projection from the 
vertex  maps $X$ into $Y\subset \mathbb{P}^{g-2}$ and 
this projection coincides with $\pi$. Every quadric $Q$ which contains 
$X\cup \{q^-\}$ should contain $C$ as well since the generators of $C$ 
are secants of $X$. By Enriques-Babbage theorem $\cap_{Q\supset X}Q$ 
equals either $X$ or a surface of degree $g-2$. Since $\deg C = g-1$ we 
conclude there is a quadric $Q$ which contains $X$ and does not 
contain $q^-$ so that Proposition~\ref{s2.15} may be applied.
\hfill $\Box$

\begin{block}\label{s2.28}
We now consider the deformation \eqref{es2.0a} in the case of triple 
covers of genus 4 and we wish to prove that the differential of the 
period mapping 
$\tilde{\varPhi}:N\times H\to \mathbb{D}\cong \mathfrak{H}_{3}$ is 
surjective for $[X\to Y]$ belonging to 
 a dense open subset of $\mathcal{H}_{3,6}(Y)$. Applying Corollary~\ref{s2.18} 
and using the irreducibility of the Hurwitz spaces 
$\mathcal{H}_{3,6}(Y)$ proved in Theorem~\ref{s3.53} it suffices to 
construct a single $\pi :X\to Y$ with non-hyperelliptic $X$ such that 
$q^- \notin Q$. We do not know explicit examples of triple covers of genus 4 which 
satisfy this condition. For example our calculations show that in the 
case of cyclic triple covers one has in fact $q^- \in Q$. 
We resolve the problem by considering singular curves and smoothing.

Let $Y$ be a fixed elliptic curve and let 
$B'=\{b_{1},b_{2},b_{3}\}\subset Y$. Let $C_{1}=Y$, let
$p_{1}:C_{1}\to Y$ be the identity mapping and let
$x_{i}=p_{1}^{-1}(b_{i})$. Let $p_{2}:C_{2}\to Y$ be an isogeny of 
degree 2. Let $p_{2}^{-1}(b_{i})=\{y_{i},y'_{i}\}$. We consider the 
following triple covering of $Y$:
\begin{equation}\label{es2.28a}
X'\ =\ C_{1}\sqcup C_{2}/\{x_{i}\sim y_{i}\}_{i=1}^{3},
\qquad \pi'=p_{1}\cup p_{2}:X'\to Y.
\end{equation}
The regular differentials of $X'$ are of the form 
$\omega=(\omega_{1}, \omega_{2})$ where $\omega_{1}\in 
H^0(C_{1},\Omega_{C_{1}}(\sum x_{i}))$, $\omega_{2}\in 
H^0(C_{2},\Omega_{C_{2}}(\sum y_{i}))$ and 
$Res_{x_{i}}\omega_{1}+Res_{y_{i}}\omega_{2}=0$ for $\forall i = 
1,2,3$. In order to describe the canonical image $\phi_{K}(X')$ we 
proceed as follows. We have embeddings 
$\phi_{1}=\phi_{|\sum x_{i}|}:C_{1}\eto \mathbb{P}^{2}_{1}$ and 
$\phi_{2}=\phi_{|\sum y_{i}|}:C_{2}\eto \mathbb{P}^{2}_{2}$ such that 
$\sum x_{i}$ and $\sum y_{i}$ are respectively pull-backs of the lines 
$\ell_{1}$ and $\ell_{2}$. There is a unique projective linear mapping 
of $\ell_{1}$ into $\ell_{2}$ that transforms $x_{i}$ into $y_{i}$. 
Identifying $\ell_{1}$ with $\ell_{2}$ along this transformation we 
obtain a reducible quadric
\[
Q\ =\ \mathbb{P}^{2}_{1}\sqcup \mathbb{P}^{2}_{2}/
\ell_{1}\sim \ell_{2}\ =\ H_{1}\cup H_{2}\; \subset \; \mathbb{P}^{3}.
\]
By dimension count one easily verifies that 
\begin{equation}\label{es2.28b}
\phi_{1}(C_{1})\cup \phi_{2}(C_{2})\quad =\quad Q\cap F
\end{equation}
where $F$ is an irreducible cubic surface in $\mathbb{P}^{3}$. Since 
the degrees of both sides in \eqref{es2.28b} are equal to 6 one 
obtains that $\phi_{1}(C_{1})\cup \phi_{2}(C_{2})$ is a canonical 
curve of arithmetic genus 4. This shows 

\smallskip
\noindent
{\bf Claim.} \emph{The canonical map $\phi_{K}:X'\to 
\mathbb{P}^{3}$ is a regular embedding and $\phi_{K}(X')=Q\cap F$ 
where $Q$ is a reducible quadric $H_{1}\cup H_{2}$ and $F$ is an 
irreducible cubic surface.}
\end{block}
\begin{pro}\label{s2.29}
Let $H^0(X',\omega_{X'})^{^{-}}$ be the hyperplane of regular 
differentials with trace 0 and let $q^- \in 
|\omega_{X'}|^* = \mathbb{P}^{3}$ be the corresponding point. Then 
$q^- \notin H_{1}\cup H_{2} = Q$. 
\end{pro}
{\bf Proof.}
Let $A=\{y_{1},y_{2},y_{3}\}$, Let $\sigma: C_{2}\to C_{2}$ be the 
involution interchanging the branches of $p_{2}:C_{2}\to Y$ and let 
$G=\{id,\sigma\}$. The mapping $p_{2}^{*}\circ Tr_{p_{2}}$ 
transforms injectively $H^0(C_{2},\Omega_{C_{2}}(A))$ into 
$H^0(C_{2},\Omega_{C_{2}}(A+\sigma A))^{G}$. Since these spaces as 
well as $H^0(Y,\Omega_{Y}(B'))$ all have dimension 3 we conclude that

\begin{equation}\label{es30.1}
Tr_{p_{2}}\ :\ H^0(C_{2},\Omega_{C_{2}}(A))\ \overset{\sim}{\lto}
H^0(Y,\Omega_{Y}(B'))
\end{equation}
is an isomorphism. Let $(\omega_{1},\omega_{2})\in 
H^0(X',\omega_{X'})$ be a regular differential with trace 
0,\quad $Tr_{\pi'}(\omega_{1},\omega_{2}) = 
Tr_{p_{1}}(\omega_{1}) + Tr_{p_{2}}(\omega_{2}) = 0$. Then $\omega 
= Tr_{p_{1}}(\omega_{1})\in H^0(Y,\Omega_{Y}(B'))$ and 
$Tr_{p_{2}}(\omega_{2}) = - \omega$. Conversely by the isomorphism 
\eqref{es30.1} the regular differentials with trace 0 are the pairs 
$(\omega_{1},\omega_{2})$ that satisfy the latter property. The fact that the 
restrictions $(\omega_{1},\omega_{2})\mapsto \omega_{1}$ and 
$(\omega_{1},\omega_{2})\mapsto \omega_{2}$ are isomorphisms of 
$H^0(X',\omega_{X'})^{^{-}} \cong 
H^0(\mathbb{P}^{3},\mathcal{O}_{\mathbb{P}^{3}}(1)\otimes J_{q})$ with 
respectively $H^0(C_{1},\Omega_{C_{1}}(\sum x_{i})) \cong 
H^0(H_{1},\mathcal{O}_{H_{1}}(1))$ and 
$H^0(C_{2},\Omega_{C_{2}}(\sum y_{i}))
\cong H^0(H_{2},\mathcal{O}_{H_{2}}(1))$ implies that $q^- \notin H_{1}\cup 
H_{2}$.
\hfill $\Box$

\begin{pro}\label{s2.36a}
Let $\pi' :X'\to Y$ be the singular triple covering of 
(\ref{s2.28}) with sufficiently general branch points 
$\{b_{1},b_{2},b_{3}\}\subset Y$. Then there exists a  
finite, flat, surjective morphism  $p=(\varPsi,f):\mathcal{X}\to Y\times 
T$ where $\mathcal{X}$ and $T$ are 
smooth and irreducible varieties such that the following properties hold: the 
morphism 
$f:\mathcal{X}\to 
T$ is proper and flat; the scheme-theoretic fibers of $f$ are  reduced and 
have dimension 1 and 
arithmetic 
genus 4; for every sufficiently general $s\in T$ the covering $p_{s} 
:\mathcal{X}_{s}\to Y\times \{s\}$ is simple; there is a point 
$s_{0}\in T$ such that $p_{s_{0}}:\mathcal{X}_{s_{0}}\to Y$ is equivalent 
to  $\pi' :X'\to Y$.
\end{pro}
{\bf Proof.}
Consider the curve $X_{0} = C_{1}\sqcup C_{2}$ and the covering 
$\pi_{0} = p_{1}\sqcup p_{2}:X_{0}\to Y$. One has 
$(p_{2})_{*}\mathcal{O}_{C_{2}} = \mathcal{O}_{Y}\oplus \eta$ where 
$\eta = Ker (Tr_{p_{2}}:(p_{2})_{*}\mathcal{O}_{C_{2}} \to 
\mathcal{O}_{Y})$ is an invertible sheaf of $Y$ with $\eta^{\otimes 
2}\cong \mathcal{O}_{Y}$. The Tschirnhausen module of $\pi_{0} 
:X_{0}\to Y$ is the quotient sheaf $E^{\vee}_{0}$ defined by 
\[
0 \lto \mathcal{O}_{Y} \overset{\pi_{0}^{\#}}{\lto} 
(\pi_{0})_{*}\mathcal{O}_{X_{0}} \lto E^{\vee}_{0} \lto 0.
\]
It is isomorphic to $Ker(
Tr_{\pi_{0}}:(\pi_{0})_{*}\mathcal{O}_{X_{0}} \to 
\mathcal{O}_{Y})
)$. There is a canonical embedding 
$X_{0}\overset{i}{\eto}\mathbf{P}(E_{0})$ such that 
$i^{*}\mathcal{O}_{\mathbf{P}(E_{0})}(1) \cong \omega_{X_{0}/Y} \cong 
\omega_{X_{0}}\otimes \pi_{0}^{*}(\omega_{Y})^{-1}$ (cf. \cite{CE}~p.448). 
We claim $E^{\vee}_{0}\cong \mathcal{O}_{Y}\oplus \eta$. Indeed, 
$(\pi_{0})_{*}\mathcal{O}_{X_{0}}\cong \mathcal{O}_{Y}\oplus 
(\mathcal{O}_{Y}\oplus \eta)$. The trace map is given on local 
sections by $Tr_{\pi_{0}}(a_{1},a_{2}+b) = a_{1}+2a_{2}$. Thus the 
homomorphism $\mathcal{O}_{Y}\oplus \eta\to Ker\, Tr_{\pi_{0}}$ given 
by $(a,b)\mapsto (-2a,a+b)$ is an isomorphism. Let 
$W_{0}=\mathbf{P}(E_{0}) = \mathbf{P}(\mathcal{O}_{Y}\oplus 
\eta)$. This is an elliptically ruled surface of invariant $e_{0}=0$ 
and with two sectional curves $Y_{0}$ and $Y_{\infty}$ with minimal 
self-intersection $-e_{0}=0$. Clearly  $C_{2}\ne Y_{0}$ and  
$C_{2}\ne Y_{\infty}$, thus one may choose $b_{1}$ in such a way that if 
$\pi_{0}^{-1}(b_{1})=\{x_{1},y_{1},y'_{1}\}$ with  $x_{1}\in C_{1},\; 
y_{1},y'_{1}\in C_{2}$ none of the points $y_{1}$ or $y'_{1}$ belongs 
to $Y_{0}\cup Y_{\infty}$. Performing an elementary transformation 
with center at the point $y'_{1}$, which consists of blowing-up 
$y'_{1}$ and blowing-down the strict transform of the fiber of 
$W_{0}\to Y$ passing through $x_{1}$ and $y_{1}$ we obtain a new 
ruled surface $W_{1}\to Y$ with invariant $e_{1}=-1$ 
\cite{Se}~p.210 
and an embedding of $X_{1} = X_{0}/x_{1}\sim y_{1}$ in $W_{1}$. Let 
$\pi_{1}:X_{1}\to Y$ be the covering induced by $\pi_{0}$. Choose a 
second non-branch point $b_{2}\in Y,\; b_{2}\ne b_{1}$ and let 
$\pi_{1}^{-1}(b_{2})=\{x_{2},y_{2},y'_{2}\}$ with $x_{2}\in C_{1}$,\; 
$y_{2},y'_{2}\in C_{2}$. There is a sectional curve with minimal 
self-intersection $-e_{1}=1$ passing through $y'_{2}$ 
\cite{Ha}~Ex.V.2.7. Performing an elementary transformation with center in 
$y'_{2}$ we obtain a new ruled surface $W_{2}\to Y$ with invariant 
$e_{2}=0$ \cite{Se}~p.210 and an embedding of $X_{2}=X_{1}/x_{2}\sim 
y_{2}$ in $W_{2}$. We claim that $W_{2}$ may have only a finite number 
of sectional curves with minimal self-intersection $-e_{2}=0$. Indeed, 
otherwise there would be an $\infty^{1}$ family of such curves and if 
$z_{2}\in W_{2}$ is the image of the blown-down curve of 
$\tilde{W}_{1}$ there would be a sectional curve with minimal 
self-intersection passing through $z_{2}$. Performing an elementary 
transformation of $W_{2}$ with center in $z_{2}$ we would obtain a 
ruled surface with invariant $e=1$. This surface is however 
isomorphic to $W_{1}$ with invariant $e_{1}=-1$ which is absurd. Let 
$\pi_{2}:X_{2}=X_{0}/\{x_{i}\sim y_{i}\}_{i=1}^{2}\to Y$ be the 
mapping induced by $\pi_{0}$. Since $C_{2}\subset X_{2}$ cannot be a 
sectional curve we may choose $b_{3}$ in such a way that if 
$\pi_{2}^{-1}(b_{3}) = \{x_{3},y_{3},y'_{3}\}$ with $x_{3}\in C_{1},\; 
y_{3},y'_{3}\in C_{2}$ then none of $y_{3}$ or $y'_{3}$ belongs to the 
union of the sectional curves with minimal self-intersection. 
Performing an elementary transformation of $W_{2}$ with center in 
$y'_{3}$ we obtain $W=W_{3}$ with invariant $e=-1$ and an embedding of 
$X' = X_{0}/\{x_{i}\sim y_{i}\}_{i=1}^{3}$ in $W$. Let $Y_{0}$ be a 
sectional curve with minimal self-intersection ${Y_{0}}^{2} = -e =1$ 
and let $F$ be a fiber of $W\to Y$. We have numerical equivalence $X' 
\equiv 3Y_{0}+bF$ since $X'\cdot F = \deg (\pi':X'\to Y) = 3$. For the 
canonical class $K_{W}$ of $W$ one has $K_{W}\equiv -2Y_{0}+F$ 
\cite{Ha}~Cor.V.2.11. Hence 
\[
2p_{a}(X')-2 = 6 = X'\cdot (X'+K_{W}) = (3Y_{0}+bF)(Y_{0}+(b+1)F) = 
3+3(b+1)+b = 6 + 4b.
\]
Therefore $b=0$. We now apply a criterion of very ampleness due to 
Biancofiore and Livorni \cite{BLi}~p.183. Namely, $b=0>3e +2 = -1$, 
hence $\mathcal{O}_{W}(X')$ is a very ample invertible sheaf. The 
required family in the proposition is constructed as follows. Let 
$\mathbb{P}^{n}=|\mathcal{O}_{W}(X')|\spcheck$ and let $\varphi : W\to 
\mathbb{P}^{n}$ be the closed embedding. We let $Z=\check{P}^{n}$, 
$\mathcal{X}' \subset W\times Z$ be the closed subset $\mathcal{X}'= 
\{(x,\gamma)|\varphi(x)\in H_{\gamma}\}$. It is clear that 
$\mathcal{X}'$ is smooth and projective. The projection 
$\pi_{1}:\mathcal{X}'\to W$ has fibers $\cong \mathbb{P}^{n-1}$, so 
$\mathcal{X}'$ is irreducible. The projection $W\to Y$ induces 
$p':\mathcal{X}'\to Y\times Z$ which is a finite surjective morphism to a 
smooth 
variety, so it is flat. Let $f'=\pi_{2}:\mathcal{X}'\to Z$. This is a flat 
morphism since it is a composition of two flat morphisms: $p':\mathcal{X}'\to 
Y\times Z$ and $Y\times Z\to Z$. Let 
$T=\{s\in Z| \mathcal{X}'_{s}\quad \text{is reduced}\}$. We prove 
below $T$ is open in $Z$. Assuming this we let 
$\mathcal{X}=(p')^{-1}(Y\times T)$ and $p:\mathcal{X}\to Y\times T$ be the 
restriction of $p'$ on $\mathcal{X}$. 

For proving the openness of $T$ we use that a 
locally noetherian scheme 
is reduced if and only if it satisfies the conditions $R_{0}$ and 
$S_{1}$ \cite{AK}~p.132.  Every $\mathcal{X}'_{s}$ is a hyperplane 
sections of the smooth surface $W$, so it is Cohen-Macaulay and 
satisfies the conditions $S_{k}$ for all $k$. The condition $R_{0}$ 
means the scheme is smooth at every generic point. Let $R'\subset 
\mathcal{X}'$ be the closed subset where the differential of $f'$ 
has not a maximal rank.
 Since $f'$ is proper $B=f'(R')$ is closed in $Z$ and 
$B\ne Z$ by Bertini's theorem. Let $B_{1}=\{s\in B | \dim 
(f'|_{R'})^{-1}(s)\geq 1\}$. This is a closed subset of B and the 
condition $R_{0}$ holds for $s\in Z$ iff $s\in Z-B_{1}$. Thus 
$T=Z-B_{1}$ is open. 

It remains to  prove 
that for sufficiently general 
$s\in T$ the covering $p_{s}:\mathcal{X}_{s}\to Y$ is simple.
 The ruled 
surface $W$ has invariant $e=-1$, so $W=\mathbf{P}(F)$ where $F$ is a 
normalized indecomposable locally free sheaf \cite{Ha}~Ch.V\S~2. If 
$E$ is the dual of the Tschirnhausen module of $\pi' :X'\to Y$, then by 
\cite{Mi}~p.1150 one has $E\cong F\otimes L$ for a certain 
invertible sheaf $L$ on $Y$, hence $E$ is indecomposable of degree 3 
(see the proof of Lemma~\ref{s3.41c}). Let $A=\det E$. By 
Proposition~\ref{s3.49} and Lemma~\ref{s3.51a} every sufficiently 
general simple covering $[X\to Y]\in \mathcal{H}_{3,A}(Y)$ has 
Tschirnhausen module isomorphic to $E^{\vee}$ and $X\in 
|\mathcal{O}_{\mathbf{P}(E)}(3)\otimes \pi^{*}(\det E)^{-1}|$. The 
curve $X'$ belongs to the same linear system. Since simpleness is an 
open condition we obtain 
that for sufficiently general 
$s\in T$ the covering $p_{s}:\mathcal{X}_{s}\to Y$ is simple.
\hfill $\Box$

\begin{pro}\label{s6516}
Let $Y$ be an elliptic curve and let $A\in Pic^{3}Y$. 
Then every sufficiently general simple triple covering $\pi :X\to Y$ 
having a Tschirnhausen module $E^{\vee}$ with $\det E\cong A$ 
satisfies the condition of Corollary~\ref{s2.18}: $X$ is not 
hyperelliptic and $q^- \notin Q$.
\end{pro}
{\bf Proof.}
Composing $\pi$ by an automorphism of $Y$ does not change neither 
$q^- \in \mathbb{P}^{3}$ nor the quadric $Q\supset \phi_{K}(X)$. By the 
argument of Lemma~\ref{s3.51a} it suffices to prove the claim only for 
one line bundle in $Pic^{3}(Y)$. Let $E^{\vee}$ be the Tschirnhausen 
module of the reducible covering $\pi' :X'\to Y$ as at the end of the 
preceding proof. Let $A = \det E$. By the proofs of 
Theorem~\ref{s3.53} (cf. Corollary~\ref{s3.55a}) and 
Proposition~\ref{s2.36a} a 
Zariski open nonempty subset of $\mathcal{H}_{3,A}(Y)$ consists of 
equivalence classes of coverings $[X\to Y]$ such that $X\in 
|\mathcal{O}_{W}(X')|$. 

We claim that for every $C\in 
|\mathcal{O}_{W}(X')|$ the dualizing sheaf $\omega_{C}$ is spanned, 
i.e. the canonical homomorphism $H^0(C,\omega_{C})\otimes 
\mathcal{O}_{C}\to \omega_{C}$ is surjective. Indeed, using the 
notation of the proof of Proposition~\ref{s2.36a} one has 
$K_{W}+X'\equiv Y_{0}+F$. Hence $\mathcal{O}_{W}(K_{W}+X')$ is 
base-point-free by \cite{BLi}~p.183. Thus by the adjunction formula 
$\omega_{C}\cong \mathcal{O}_{W}(K_{W}+X')\otimes \mathcal{O}_{C}$ is 
spanned. Let us consider the family $f:\mathcal{X}\to T$ constructed 
in Proposition~\ref{s2.36a}. Let $\omega_{_{\mathcal{X}/T}} := 
\Omega_{\mathcal{X}}^{n}\otimes f^{*}(\Omega_{T}^{n-1})^{-1}$. By 
flatness $f_{*}\omega_{_{\mathcal{X}/T}}$ is locally free of rank 4. By 
the above claim the canonical morphism 
$f^{*}(f_{*}\omega_{_{\mathcal{X}/T}})\to \omega_{_{\mathcal{X}/T}}$ is an 
epimorphism. This yields a well-defined relative canonical map 
\begin{equation*}
\begin{diagram}
\mathcal{X}&    &\rTo^{\varphi} &       &\mathbf{P}(f_{*}\omega_{_{\mathcal{X}/T}}) = 
\mathbb{P} \\
           &\rdTo_f &     & \ldTo_g  &               \\
           &      &  T  &        &               \\
\end{diagram}
\end{equation*}
The trace mapping $Tr_{p}:p_{*}(\Omega^{n}_{\mathcal{X}})\to 
\Omega^{n}_{Y\times T} 
$ is an epimorphism since for every $\Phi \in \Gamma(V,
\Omega^{n}_{Y\times T})$ one has $Tr_{p}(\frac{1}{\deg p}p^{*}\Phi) = 
\Phi$. This yields an epimorphism $Tr : 
f_{*}\omega_{_{\mathcal{X}/T}}\to H^0(Y,\omega_{Y})\otimes 
\mathcal{O}_{T} = \mathcal{O}_{T}$, hence a section $\mu$ of 
$g:\mathbb{P}\to T$. According to Proposition~\ref{prel2} the 
restriction of $Tr$ at the fiber over any $s\in T$ equals 
$Tr_{p_{s}}:H^0(\mathcal{X}_{s},\omega_{\mathcal{X}_{s}})\to 
H^0(Y,\omega_{Y})$. Hence $\mu(s)$ is the point $q^-$ for 
$p_{s}:\mathcal{X}_{s}\to Y$ as defined in (\ref{s2.11}).

Consider the homomorphism 
$\varphi^{*}:g_{*}\mathcal{O}_{\mathbb{P}}(2)\to f_{*}\omega^{\otimes 
2}_{\mathcal{X}/T}$. If $s_{0}\in T$ is the point which corresponds to 
the reducible covering $\pi' :X'=C_{1}\cup C_{2}\to Y$, then 
$\varphi^{*}\otimes k(s_{0}):S^{2}H^0(\mathcal{X}_{s_{0}},\omega_{
\mathcal{X}_{s_{0}}
})\to H^0(\mathcal{X}_{s_{0}},\omega_{\mathcal{X}_{s_{0}}}^{\otimes 
2})$ is surjective since $\phi_{K}(\mathcal{X}_{s_{0}})$ is contained in a 
unique quadric. Hence replacing $T$  by an affine 
neighborhood of $s_{0}$ if necessary we may assume that for every $s\in T$ the 
fiber 
homomorphism $\varphi^{*}\otimes k(s)$ is surjective with one 
dimensional kernel. In particular if $\mathcal{X}_{s}$ is smooth it is 
non-hyperelliptic. We obtain a relative quadric 
$\mathcal{Q}\subset \mathbb{P}$ which contains $\varphi(\mathcal{X})$. 
The set $V = \{s\in T | \mu(s)\in \mathcal{Q}(s)\}$ is closed in $T$ and 
$s_{0}\notin V$. The proposition is proved.
\hfill $\Box$

\section{Unirationality results}\label{s5}
We first give an alternative proof of a known result \cite{BL1}.
\begin{thm}\label{s4.76}
The moduli spaces of polarized abelian surfaces $\mathcal{A}_{2}(1,2)$ 
and $\mathcal{A}_{2}(1,3)$ are unirational.
\end{thm}
{\bf Proof.}
Let $d=2$ or $3$. We fix an elliptic curve $Y$  
and let $T$ be the Hurwitz space $\mathcal{H}_{d,4}(Y)$. We 
consider a commutative diagram as in \eqref{es4.73}
\begin{equation*}
\begin{diagram}
\mathcal{X}&    &\rTo^p &       &Y\times T \\
           &\rdTo_f &     & \ldTo_{q=\pi_{2}}  &       \\
           &      &  T  &        &           \\
\end{diagram}
\end{equation*}
where\; $p$\; is the universal family of simple 
coverings  of  degree\:  $d$\:  branched  in  4   points.   According   to 
Theorem~\ref{s3.53}\; $T$ is irreducible. Applying 
Lemma~\ref{s1.1} and Proposition~\ref{s4.74} we obtain a period 
mapping $\varPhi:T\to \mathcal{A}_{2}(1,d)$. We claim
 $\varPhi$ is dominant. Let $s_{0}\in T$ and let  $H\subset T$  be  a 
polydisk centered in $s_{0}$ as in (\ref{s2.0}). Then we have a  lifting 
of the period mapping and a commutative diagram
\begin{equation*}
\begin{diagram}
H      & \rTo^{\tilde{\varPhi}} & \mathfrak{H}_{2} \\
\dTo_i  &                      &   \dTo           \\
T          & \rTo^{\varPhi} &  \mathcal{A}_{2}(1,3) \\
\end{diagram}
\end{equation*}
By Propositions~\ref{s2.13a} and (\ref{s2.12}) if $s_0$ is chosen 
general enough 
the differential $d\tilde{\varPhi}(s_{0})$ is surjective. Hence by the 
implicit 
function theorem a neighborhood of $\tilde{\varPhi}(s_{0})$ is contained in 
$\tilde{\varPhi}(H)$. Therefore 
$\varPhi$ is dominant. Let us choose $A\in Pic^{2}Y$. Then the images
by $\varPhi$ of $\mathcal{H}_{d,4}(Y)$ and $\mathcal{H}_{d,A}(Y)$ are 
the same (cf. Proposition~\ref{s4.74} and the proof of 
Lemma~\ref{s3.51a}). Thus the restriction of 
$\varPhi:\mathcal{H}_{d,A}(Y)\to \mathcal{A}_{2}(1,d)$ is dominant as 
well. By Theorem~\ref{s3.53}\quad $\mathcal{H}_{d,A}(Y)$ is 
unirational. Therefore $\mathcal{A}_{2}(1,d)$ is unirational.
\hfill $\Box$

\bigskip
\noindent
Gritsenko proved in \cite{Gri1} that if 
$\tilde{\mathcal{A}}_{2}(1,d)$ is a nonsingular, projective model of 
$\mathcal{A}_{2}(1,d)$, then the geometric genus $p_{g}(
\tilde{\mathcal{A}}_{2}(1,d)
)\geq 1$ when $d\geq 13$ and $d\neq 14,15,16,18,20,24$, $30, 36$. Sankaran 
proved in \cite{Sa}
that if $d$ is prime and $d\geq 173$ then $\tilde{\mathcal{A}}_{2}(1,d)$ is 
of general type. If $g(Y)=1$ none of the  Hurwitz spaces
$\mathcal{H}_{d,n}(Y)$ is uniruled because of the 
epimorphism $h: \mathcal{H}_{d,n}(Y)\to Pic^{n/2}Y$ (cf. (\ref{s3.51})). 
We denote by $\mathcal{H}^{0}_{d,n}(Y)$ the subset of 
$\mathcal{H}_{d,n}(Y)$ whose points correspond to coverings $\pi :X\to 
Y$  with the property that $\pi_*:H_1(X,\mathbb{Z})\to 
H_1(Y,\mathbb{Z})$ is surjective, or equivalently that 
$ \pi^*:J(Y)\to J(X)$ is injective. This property is preserved under 
deformation, so $\mathcal{H}^{0}_{d,n}(Y)$ is a union of connected 
components of $\mathcal{H}_{d,n}(Y)$. We denote by $\mathcal{H}^{0}_{d,A}(Y)$ 
the intersection 
$\mathcal{H}^{0}_{d,A}(Y)\cap \mathcal{H}_{d,A}(Y)$. 
We notice that $\mathcal{H}^{0}_{d,n}(Y)$ and $\mathcal{H}^{0}_{d,A}(Y)$ are 
non-empty for every $n\geq 2$ and every $A\in Pic^{n/2}Y$ as follows from 
Lemma~\ref{s3.41}, Lemma~\ref{s1.3} and Lemma~\ref{s3.51a}.
In a direction opposite to the one of Theorem~\ref{s4.76} we have the 
following result.

\begin{thm}\label{s4.77a}Let 
$Y$ be an elliptic curve and let $A\in Pic^{2}Y$.
Let  $d\geq 13$ and let $d\ne 
14,15,16,18,20,24,30,36$. Then every connected component of 
$\mathcal{H}^{0}_{d,A}(Y)$ has the property that each of its 
nonsingular, projective models has geometric genus $\geq 1$. In 
particular none of the connected components of 
$\mathcal{H}^{0}_{d,A}(Y)$ is uniruled. If 
furthermore $d$ is prime and $d\geq 173$, then every projective 
nonsingular model of any connected component of 
$\mathcal{H}_{d,A}(Y)$ is of 
general type.
\end{thm}
{\bf Proof.}
Using Lemma~\ref{s1.3} and repeating the argument of Theorem~\ref{s4.76} 
we obtain 
a morphism $\varPhi :\mathcal{H}^0_{d,A}(Y)\to 
\mathcal{A}_{2}(1,d)$ whose restriction on every connected 
 component is dominant. The theorem follows thus from the 
results of Gritsenko and Sankaran cited above.
\hfill $\Box$

\begin{thm}\label{s4.77}
The moduli spaces of polarized abelian threefolds 
$\mathcal{A}_{3}(1,1,2),\; \mathcal{A}_{3}(1,2,2)$,\linebreak 
$\mathcal{A}_{3}(1,1,3)$ and $\mathcal{A}_{3}(1,3,3)$ are unirational.
\end{thm}
{\bf Proof.}
Let $d=2$ or 3. We consider the family of coverings from (\ref{s3.55}) 
and Proposition~\ref{s3.57}. According to Lemma~\ref{s1.1} and 
Proposition~\ref{s4.74} the period mapping $\varPhi:T\to 
\mathcal{A}_{3}(1,1,d)$ defined by $\varPhi(s)=[Ker(Nm_{p_{s}})]$ is 
an algebraic morphism. We wish to prove $\varPhi$ is dominant. Let 
$s_{0}\in T$. Let us denote by $\pi :X\to Y$ the covering 
corresponding to $s_{0}$. By the definition of $\varPhi$ (cf. 
(\ref{s4.70})) one may choose a neighborhood $S$ of $s_{0}$ (in the 
Hausdorff topology) such that $\varPhi|_{S}$ may be lifted to a 
holomorphic mapping in $\mathfrak{H}_{3}$
\[
\begin{diagram}
S      & \rTo^{\tilde{\varPhi}'} & \mathfrak{H}_{3} \\
\dTo_i  &                      &   \dTo           \\
T          & \rTo^{\varPhi} &  \mathcal{A}_{3}(1,1,3) \\
\end{diagram}
\]
Consider the restriction of the family \eqref{es3.56} on $S$. 
Applying Proposition~\ref{s2.0b} to the deformation \eqref{es3.56a} we 
conclude that,
 shrinking $S$ if necessary, there is a holomorphic 
mapping $\mu:S\to N\times H$ such that the family of coverings over $S$ 
is the pull-back of the family induced from \eqref{es2.0a}. Considering 
the period 
mappings one obtains a commutative diagram
\begin{diagram}
S          & \rTo^{\tilde{\varPhi}'} & \mathfrak{H}_{3} \\
\dTo^{\mu} & \ruTo_{\tilde{\varPhi}} &                  \\
N\times H  &                         &                   \\
\end{diagram}
By Proposition~\ref{s6516} if $s_{0}\in T$ is chosen general enough 
the differential $d\tilde{\varPhi}(s_{0})$ is epimorphic. Hence by the 
implicit function theorem a neighborhood of $\tilde{\varPhi}(s_{0})$ 
is contained in $\tilde{\varPhi}(N\times H)$. Consider the 
Tschirnhausen module $\mathcal{E}^{\vee}$ of the covering 
$\mathcal{X}\to \mathcal{Y}\times _{N}(N\times H)$ induced from 
\eqref{es2.0a}. The set of $u\in N\times H$ such that 
$\mathcal{E}_{u}$ is stable is open in $N\times H$ (cf. Proposition~\ref{B1}). 
For vector bundles of rank 2 and degree 3 over an elliptic curve 
being 
stable and being indecomposable are equivalent conditions  
\cite{Tu}~p.20. Hence by Proposition~\ref{s3.49} if one chooses 
$s_{0}\in T$ general enough a neighborhood of $\mu(s_{0})$ in $N\times 
H$ corresponds to triple coverings with indecomposable Tschirnhausen 
modules of degree $-3$. Given a triple covering of an elliptic curve, 
composing the covering by a translation results in translation of the 
determinant of the Tschirnhausen module (cf. the proof of 
Lemma~\ref{s3.51a}) while the kernel of the norm map of the Jacobians 
remains the same. This shows that $\tilde{\varPhi}'(S)$
 contains the image by 
$\tilde{\varPhi}$ of a certain neighborhood of $\mu(s_{0})$ in $N\times 
H$. Therefore a neighborhood of $\tilde{\varPhi}'(s_{0})$ in 
$\mathfrak{H}_{3}$ is contained in $\tilde{\varPhi}'(S)$. This shows 
that $\varPhi:T\to \mathcal{A}_{3}(1,1,d)$ is dominant.

By construction $T$ is rational. Therefore $\mathcal{A}_{3}(1,1,d)$ is 
unirational. The unirationality of $\mathcal{A}_{3}(1,d,d)$ then 
follows either from the result of Birkenhake and Lange \cite{BL3} or 
from the weaker statement in Proposition~\ref{s4.72}.
\hfill $\Box$

\bigskip
\noindent
The proofs of  Theorem~\ref{s4.76}, Theorem~\ref{s4.77a} and 
Theorem~\ref{s4.77} together with Proposition~\ref{s4.74} yield 
the following corollary.

\begin{cor}\label{s4.80}
Let $Y$ be an elliptic curve.
Every sufficiently general abelian surface with polarization of type $(1,d)$ 
is isomorphic to the Prym variety of a simple degree $d$ covering of  $Y$ 
branched in 4 points.
For every sufficiently general abelian threefold $A$ with polarization of 
type $(1,1,d)$, $d=2$ or $3$ there exists an elliptic curve $Y$ (depending on 
$A$) such that $A$ is isomorphic to the Prym variety of a simple covering 
$\pi:X\to Y$ of degree $d$ branched in 6 points. 
 For every sufficiently 
general abelian threefold $B$ with polarization of type $(1,d,d)$, $d=2$ or 
$3$ 
there exists an elliptic curve $Y$ and a simple covering $\pi:X\to Y$ of 
degree $d$ branched in 6 points such that $B$
is 
isomorphic to $Pic^{0}X/\pi^{*}Pic^{0}Y$. 
\end{cor}

\appendix
\section{Traces of differential forms}\label{prel}
\begin{block}\label{prel1}
In this appendix we assume the base field $k$ is algebraically closed of 
characteristic 0.
Let $X$ and $Y$ be smooth, irreducible varieties of dimension $n$, let 
$p:X\to Y$ be a finite, surjective morphism. One defines the trace 
mapping 
$Tr_{p}:\Omega^{n}(X)\to \Omega^{n}(Y)$ for rational differential forms 
$n$-forms as follows. Let $Y_0\subset Y$ be an open subset such that if 
$X_{0}=p^{-1}(Y_{0})$ the restriction  $p:X_{0}\to Y_{0}$ 
is \'{e}tale. If $v_{1},\ldots,v_{n}$ are 
local parameters at some point $y\in Y_{0}$, then 
$Tr_{p}(a\,dp^{*}v_{1}\wedge \cdots \wedge dp^{*}v_{n}) := Tr_{p}(a)
dv_{1}\wedge \cdots \wedge dv_{n}
$. It is well-known that the trace mapping transforms regular 
differentials of $X$ into regular differentials of $Y$ 
(cf. \cite{Li}~Example 2.1.2). More 
generally one may define in the same way $Tr_{p}:\Omega^{n}(X)\to 
\Omega^{n}(Y)$ for every finite, surjective morphism $p:X\to Y$ between 
reduced, equidimensional schemes of dimension $n$ where $\Omega^{n}(X)$ and 
$\Omega^{n}(Y)$ are the  sheaves of rational $n$-forms.
\end{block}

\begin{pro}\label{prel2}
Suppose we have a commutative diagram 
\begin{equation*}
\begin{diagram}
\mathcal{X}&    &\rTo^p &       &\mathcal{Y} \\
           &\rdTo_f &     & \ldTo_g  &       \\
           &      &  T  &        &           \\
\end{diagram}
\end{equation*}
where $\mathcal{X},\mathcal{Y}$ and $T$ are smooth and irreducible, 
$T$ is affine, $\dim \mathcal{X} =\dim \mathcal{Y} = n$, $p$ is 
finite and surjective, the morphisms $f$ and $g$ are proper of 
relative dimension 1 with reduced fibers. Suppose the discriminant 
subscheme of $p$ does not contain fibers of $g$. Then for every 
$s\in T$ one has canonical isomorphisms 
given by Poincar\'{e} residues $(\Omega^{n}_{\mathcal{X}})_{s}\cong 
\omega_{\mathcal{X}_{s}},\; (\Omega_{\mathcal{Y}}^{n})_{s}\cong 
\omega_{\mathcal{Y}_{s}}$ and the following commutative diagram holds

\begin{equation}\label{eprel2}
\begin{diagram}
H^0(\mathcal{X},\Omega^{n}_{\mathcal{X}})\otimes k(s) & 
\rTo^{Tr_{p}\otimes k(s)} &
H^0(\mathcal{Y},\Omega_{\mathcal{Y}}^{n})\otimes k(s) \\
\dTo^{\cong} &          & \dTo_{\cong}     \\
H^0(\mathcal{X}_{s},\omega_{\mathcal{X}_{s}}) &
\rTo^{Tr_{p_{s}}}  &
H^0(\mathcal{Y}_{s},\omega_{\mathcal{Y}_{s}}) \\
\end{diagram}
\end{equation}
In particular for every $s\in T$ the trace mapping $Tr_{p_{s}}$ 
 transforms 
the regular differentials of $\mathcal{X}_{s}$ into regular differentials 
of $\mathcal{Y}_{s}$.
\end{pro}
{\bf Proof.}
Let $t_{1},\ldots,t_{n-1}$ be local parameters of $T$ at $s$. If $x\in 
\mathcal{X}_{s}$ and if $u_{1},\ldots,u_{n}$ are local parameters of 
$\mathcal{X}$ at $x$ then a generator of 
$(\omega_{_{\mathcal{X}_{s}}})_{x}$ is given by the Poincar\'{e} residue 
$\frac{du_{1}\wedge\cdots \wedge du_{n}}{f^{*}(dt_{1}\wedge \cdots 
\wedge dt_{n-1})}|_{\mathcal{X}_{s}}$. The vertical homomorphisms in 
\eqref{eprel2} are induced by the Poincar\'{e} residue and it is a standard 
fact that they are isomorphisms when $T$ is affine. Let 
$\mathcal{Y}_{0}\subset \mathcal{Y}$,\; $\mathcal{X}_{0}=p^{-1}(
\mathcal{Y}_{0}
)$ be open subsets such that $p:\mathcal{X}_{0}\to 
\mathcal{Y}_{0}$ is \'{e}tale and $\mathcal{Y}_{0}\cap 
\mathcal{Y}_{s}$ and $\mathcal{X}_{0}\cap \mathcal{X}_{s}$ are smooth 
and dense. 
In order to verify \eqref{eprel2} it suffices to prove the 
commutativity of an analogous diagram with $\mathcal{X},\mathcal{Y}$ 
replaced by $\mathcal{X}_{0},\mathcal{Y}_{0}$  and 
$\mathcal{X}_{s},\mathcal{Y}_{s}$ replaced by 
$\mathcal{X}_{s}\cap \mathcal{X}_{0},\mathcal{Y}_{s}\cap \mathcal{Y}_{0}$ 
respectively. Let $\varphi\in H^0(\mathcal{X}_{0},
\Omega^{n}_{\mathcal{X}})$.
Let $y\in \mathcal{Y}_{s}\cap 
\mathcal{Y}_{0}$ and let 
$v_{1},\ldots,v_{n}$ be local parameters of $\mathcal{Y}$ at $y$ 
defined in an affine neighborhood $V\subset \mathcal{Y}_{0}$. Let 
$U=p^{-1}(V)$ and let $\varphi|_{U}=a\, p^{*}(dv_{1}\wedge \cdots 
\wedge dv_{n})$ with $a\in \Gamma(U,\mathcal{O}_{\mathcal{X}})$. Then 
$Tr_{p}(\varphi)=Tr_{p}(a)dv_{1}\wedge \cdots \wedge dv_{n}$. The 
Poincar\'{e} residue of this form equals 
$Tr_{p}(a)|_{\mathcal{Y}_{s}}\frac{
dv_{1}\wedge \cdots \wedge dv_{n}
}{g^{*}(dt_{1}\wedge \cdots \wedge dt_{n-1})}|_{\mathcal{Y}_{s}}$. The 
Poincar\'{e} residue of $\varphi$ on $\mathcal{X}_{s}$ equals 
\[
a\bigr|_{\mathcal{X}_{s}}
\frac{p^{*}dv_{1}\wedge \cdots \wedge dv_{n}}
{f^{*}dt_{1}\wedge \cdots \wedge dt_{n-1}}\bigr|_{\mathcal{X}_{s}}\ =\ 
a\bigr|_{\mathcal{X}_{s}}p^{*}_{s}\left(
\frac{dv_{1}\wedge \cdots \wedge dv_{n}}
{g^{*}dt_{1}\wedge \cdots \wedge dt_{n-1}}
\bigr|_{\mathcal{X}_{s}}\right).
\]
Applying $Tr_{p_{s}}$ to this differential one obtains
$Tr_{p_{s}}(a|_{\mathcal{X}_{s}})\frac{
dv_{1}\wedge \cdots \wedge dv_{n}
}{g^{*}(dt_{1}\wedge \cdots \wedge dt_{n-1})}|_{\mathcal{Y}_{s}}$.
So the commutativity of \eqref{eprel2} follows from the equality 
$Tr_{p}(a)|_{\mathcal{Y}_{s}} = Tr_{p_{s}}(a|_{\mathcal{X}_{s}})$ (cf. 
\cite{AK}~p.123).
\hfill $\Box$

\section{Openness conditions for families of vector bundles over families of elliptic curves}\label{a2}
In this appendix we will make the customary identification between 
vector bundles and locally free sheaves. In the next proposition we 
prove a result which should be well-known, but we could not find an 
appropriate reference. The openness of the stability and the semistability conditions for families of vector bundles over a 
fixed curve of genus $\geq 1$ is proved in \cite{NS}~Theorem~2 and in 
fact their proof is easily adapted to the case we consider in the 
proposition below. Families of 
vector bundles varying over a family of curves of genus 
$\geq 2$ are treated in \cite{Pa}~p.459. 

Given a vector bundle over a smooth, 
projective curve we denote by $r(E),\, d(E)$ and $\mu(E)=d(E)/r(E)$ 
the rank, degree and the slope of $E$.
The following facts about vector 
bundles over an elliptic curve are known (see e.g. \cite{Oe},\cite{Br},\cite{Tu})
 A vector bundle is 
semistable if and only if it is a direct sum of indecomposable vector 
bundles of the same slope. A vector bundle  
is stable if and only if it is indecomposable and its rank and degree are coprime. If one fixes one indecomposable 
vector bundle $E(r,d)$ of rank $r$ and degree $d$ then all others  are isomorphic to $E(r,d)\otimes L$ 
where $L\in Pic^{0}Y$ \cite{At}~Theorem~10. 
\begin{dfn}
Let $E$ be a semistable vector bundle  over an elliptic curve.\quad $E$ is called 
\emph{polystable} if $E\cong E_1\oplus \cdots \oplus E_n$, where $E_i$ is stable 
for every $i=1,\ldots,n$. Let $r(E)=rh$, $d(E)=dh$, where $(r,d)=1$.\quad $E$ is called \emph{regular} if $h^0(\End E)=h$ (cf. \cite{FMW}).
\end{dfn}
\begin{lem}\label{B1b}
Let $E$ be a semistable vector bundle over an elliptic curve of rank $rh$ and degree $dh$, where $(r,d)=1$. Then $h^0(\End E)\geq h$. The vector bundle $E$ is regular, polystable if and only if $E\cong E_1\oplus \cdots \oplus E_h$, where every $E_i$ is indecomposable of rank $r$ and degree $d$ and $E_i\ncong E_j$ for $\forall i\ne j$.
\end{lem}
{\bf Proof.}
Let $E\cong E_1\oplus \cdots \oplus E_n$ where every $E_i$ is indecomposable. We have $\mu(E_i)=\mu(E)=d/r$, so $r(E_i)=rh_i,\; d(E_i)=dh_i$. By \cite{At}~Lemma~24 and the proof of Lemma~23 one has for every $i=1,\ldots,n$ that  $E_i\cong E'_i\otimes F_{h_i}$, where $E'_i$ is indecomposable of rank $r$ and degree $d$ and furthermore $h^0(\End E_i)=h_i$. Thus $h^0(\End E)\geq \sum h_i = h$. The vector bundle $E$ is polystable if and only $h_i=1$ for every $i$ and therefore $n=h$. It is regular polystable if moreover $h^0(E_i\otimes E_j^{\vee})=0$, or equivalently $E_i\ncong E_j$ for $\forall i\ne j$,  
\hfill $\Box$
\begin{lem}\label{B1c}
Let $E$ be a vector bundle over an elliptic curve. Then $E$ 
is not stable (resp. 
not semistable) if and only if there exists a stable vector bundle $F$
of rank $<r(E)$ and slope $\mu(F)\geq \mu(E)$ (resp. $\mu(F)> \mu(E)$) 
and a nonzero homomorphism of $F$ into $E$.\quad $E$ is semistable, but is not regular polystable if and only if there exists an indecomposable vector bundle $G$ with $\mu(G)=\mu(E)$ and $r(G)\leq r(E)$ such that $\dim Hom(G,E)\geq 2$.
\end{lem}
{\bf Proof.} The statement about non-stability resp. non-semistability is proved in \cite{NS}~Proposition~4.6 for vector bundles over smooth, projective curves of arbitrary genus $\geq 1$. Let $E$ be semistable and let $E\cong E_1\oplus \cdots \oplus E_n$ be a direct sum of indecomposable vector bundles, which by the semistability have  slopes equal to $\mu(E)$. Suppose $E$ is not regular polystable. Then either  one of $E_i$ has $(r(E_i),d(E_i))>1$ and in this case we let $G=E_i$, or  $E_j\cong E_k$ for some pair $j\ne k$ and in this case we let $G=E_j$. Then clearly $\dim Hom(G,E)\geq 2$. Conversely, suppose $E$ is regular polystable. Let $r(E)=rh,\; d(E)=dh$ where $(r,d)=1$ and let $E\cong E_1\oplus \cdots \oplus E_h$ be the decomposition of Lemma~\ref{B1b}. If $G$ is indecomposable with $\mu(G)=d/r$, then $G\cong G'\otimes F_{\ell}$, where $G'$ is indecomposable of rank $r$ and degree $d$ \cite{At}~Lemma~26. For every $i$ we have $G^{\vee}\otimes E_i\cong (G')^{\vee}\otimes E_i \otimes F_{\ell}^{\vee} \cong (\sum_{j=1}^{r^2}L_j)\otimes F_{\ell}$ where $L_j$ are line bundles of degree 0  (cf. \cite{At}~pp.434,437,439). Only one of  $L_j$ might be isomorphic to the trivial line bundle and this happens if and only if $G'\cong E_i$. By \cite{At}~Theorem~5 we conclude $\dim Hom(G,E) = h^0(G^{\vee}\otimes E)\leq 1$.
\hfill $\Box$

\begin{pro}\label{B1}
Let $q:\mathcal{Y}\to B$ be a smooth family of elliptic curves. Here in the algebraic setting $B$ 
is a scheme and $q$ is a smooth morphism, in the complex analytic 
setting $B$ is an analytic space and $q$ is a smooth holomorphic 
mapping. Let $\mathcal{E}$ be a vector bundle over $\mathcal{Y}$. Then 
the sets $B_{s}=\{b\in B|\:
\mathcal{E}_{b}\; \text{is stable over}\; \mathcal{Y}_{b}\}$,
$B_{ss}=\{b\in B|\:
\mathcal{E}_{b}\; \text{is semistable over}\; \mathcal{Y}_{b}\}$,
$B_{rss}=\{b\in B|\:
\mathcal{E}_{b}\; \text{is regular semistable over}\; \mathcal{Y}_{b}\}$ and
$B_{rps}=\{b\in B|\:
\mathcal{E}_{b}\; \text{is regular polystable over}\; \mathcal{Y}_{b}\}$
are all open in $B$.
\end{pro}
{\bf Proof.}
Replacing $B$ by $B_{red}$ we may assume $B$ is reduced. The statement 
is local so we may furter assume $B$ is connected and is either affine 
scheme (in the algebraic setting) or is a Stein space (in the complex 
analytic setting). Furthermore it is obvious it suffices to prove the statement of the proposition for elliptic fibrations which have a section 
$\sigma : B\to \mathcal{Y}$. 
We need a lemma.

\begin{lem}[Atiyah]\label{B1a}
Let $r,d$ be a pair of integers, $r\geq 1$. Then there exists a vector 
bundle $E(r,d)$ over $\mathcal{Y}$ such that for every $b\in B$ the 
fiber $E(r,d)_{b}$ is indecomposable of rank $r$ and degree $d$ over $\mathcal{Y}_{b}$ and if 
$d=0$ it holds $h^0(\mathcal{Y}_{b},E(r,d)_{b})= 1$.
\end{lem}
{\bf Proof.}
We only indicate the modifications one needs to make in order to apply 
the arguments of \cite{At}. Let $D=\sigma(B)$. Tensoring by powers of 
$\mathcal{L}(D)$ we see it suffices to construct $E(r,d)$ for $0\leq 
d<r$. One proceeds by induction on $r$. If $r=1$ one lets 
$E(1,0)=\mathcal{O}_{\mathcal{Y}}$. Suppose $E(\ell,d)$ is constructed 
for all $\ell <r$. Consider $E(r-d,d)$. By \cite{At}~Lemma~15 one has
\linebreak 
$h^1(\mathcal{Y}_{b},E(r-d,d)_{b}^{\vee})=h$ where $h=\max\{1,d\}$. By 
Grauert's theorem the sheaf $G=R^{1}q_{*}E(r-d,d)^{\vee}$ is locally 
free. The extensions of $q^{*}G$ by $E(r-d,d)^{\vee}$ are classified 
by $Ext^{1}(q^{*}G,E(r-d,d)^{\vee})\cong 
H^1(\mathcal{Y},q^{*}G^{\vee}\otimes E(r-d,d)^{\vee})\cong 
H^0(B,R^{1}q_{*}(q^{*}G^{\vee}\otimes E(r-d,d)^{\vee}))\cong 
Hom_{B}(G,R^{1}q_{*}E(r-d,d)^{\vee})$. To the identity homomorphism it 
corresponds an extension
\[
0\lto  E(r-d,d)^{\vee} \lto \mathcal{F} \lto \mathcal{G} \lto 0.
\]
One lets $E(r,d)=\mathcal{F}^{\vee}$.
\hfill $\Box$

\medskip
\noindent
We continue the proof of Proposition~\ref{B1}. Let us consider 
$\mathcal{Y}\times_{B} \mathcal{Y}$. Let $\Delta\subset 
\mathcal{Y}\times_{B} \mathcal{Y}$ be the diagonal. Let 
$\mathcal{E}(r,d)=p_{1}^{*}E(r,d)\otimes 
\mathcal{L}(\Delta-p_{2}^{*}D)$ where 
$p_{i}: \mathcal{Y}\times_{B} \mathcal{Y}\to \mathcal{Y}$ 
are the two 
projections. Looking at 
$p_{2}: \mathcal{Y}\times_{B} \mathcal{Y}\to \mathcal{Y}$ as a base 
change of the fibration $q:\mathcal{Y}\to B$ we have that 
$\mathcal{E}(r,d)$ is a Poincar\'{e} vector bundle for the elliptic 
fibration $q:\mathcal{Y}\to B$. Namely, every fiber 
$\mathcal{E}(r,d)_{z},\quad z\in \mathcal{Y}$ is an indecomposable 
vector bundle over $\mathcal{Y}_{b}$, where $q(z)=b$, and for every 
indecomposable vector bundle $E$ of rank $r$ and degree $d$ over 
$\mathcal{Y}_{b}$ there exists a $z\in \mathcal{Y}_{b}$ such that 
$E\cong \mathcal{E}(r,d)_{z}$ (cf. \cite{At}~Theorem~10). 

 Let 
$\mu=d(\mathcal{E}_{b})/r(\mathcal{E}_{b})$. Let us fix integers $r\geq 1$ 
and $d\in \mathbb{Z}$ such that $(r,d)=1$ and $d/r\geq \mu$ (resp. 
$>\mu$). Let $A\subset B$ be the subset of $b\in B$ such that there 
exists a stable vector bundle $F$ over $\mathcal{Y}_{b}$ of rank $r$ 
and degree $d$ and a nonzero homomorphism $F\to \mathcal{E}_{b}$. We 
claim $A$ is closed in $B$. Indeed, consider 
$\mathcal{H}=\mathcal{E}(r,d)^{\vee}\otimes 
p_{1}^{*}\mathcal{E}$ over $\mathcal{Y}\times_{B} \mathcal{Y}$. From the upper semi-continuity theorem for 
cohomology applied to $\mathcal{H}$ and 
$p_{2}: \mathcal{Y}\times_{B} \mathcal{Y}\to \mathcal{Y}$
it follows that the set $S=\{z\in 
\mathcal{Y}|h^0(\mathcal{Y}_{b},\mathcal{E}(r,d)^{\vee}_{z}\otimes 
\mathcal{E}_{b})\geq 1\quad \text{for}\quad b=q(z)\}$ is closed in 
$\mathcal{Y}$. Since $q:\mathcal{Y}\to B$ is proper and $A=q(S)$ we 
conclude $A$ is closed in $B$.

If a vector bundle $E$ over a smooth, projective curve $X$ is not 
stable (resp. not semistable) there exists a proper stable vector 
subbundle $F$ of $E$ such that $\mu(F)\geq \mu(E)$ (resp. $\mu(F)> 
\mu(E)$) (cf. \cite{NS}~Proposition~4.5). Furthermore the slopes of 
such vector bundles are bounded from above by a constant. We recall 
how one obtains a majorant for $\mu(F)$ (see e.g. \cite{Br}~p.82). Let 
$L$ be a line bundle on $X$ of positive degree and let 
$E^{\vee}\otimes L^{m},\quad m\gg 0$ be generated by global sections. 
Then $E$ is isomorphic to a vector subbundle of the semistable bundle 
$(L^{m})^{\oplus N}$, so for each vector subbundle $F\subset E$ one 
has $\mu(F)\leq m\deg(L)$. 

If $q:\mathcal{Y}\to B$ is an algebraic family then there exists an 
$m\gg 0$ such that $q^{*}q_{*}
\mathcal{E}^{\vee}\otimes \mathcal{L}(mD)\to 
\mathcal{E}^{\vee}\otimes \mathcal{L}(mD)$ is surjective. If 
$\mathcal{E}_{b}$ is not stable (resp. not semistable) then as we saw 
in the preceeding paragraph there exists a proper, stable vector 
subbundle $F\subset \mathcal{E}_{b}$ such that 
$\mu=\mu(\mathcal{E}_{b)}\leq \mu(F)\leq m$ (resp.
$\mu=\mu(\mathcal{E}_{b)}< \mu(F)\leq m$). The rank and degree of such 
vector subbundles belong to a finite set of integers. So, from what 
we proved above and from Lemma~\ref{B1c} the set $B-B_{s}$ (resp. $B-B_{ss}$) consisting
of $b\in B$ such that $\mathcal{E}_{b}$ is not stable (resp. not 
semistable) is a finite union of closed subsets of $B$ and hence it is 
closed in $B$.

In the complex analytic setting a small change of the above argument 
is necessary. It suffices to check that $U\cap (B-B_{s})$ (resp. 
$U\cap (B-B_{ss})$) is closed for every relatively compact open subset 
$U\subset B$. For each such $U$ there exists an $m\gg 0$ such that 
$q^{*}(q_{*}
\mathcal{E}^{\vee}\otimes \mathcal{L}(mD)|_{U})\to 
\mathcal{E}^{\vee}\otimes \mathcal{L}(mD)|_{q^{-1}U}$ is surjective. 
Repeating the argument of the preceeding paragraph we see that both 
$B-B_{s}$ and $B-B_{ss}$ are closed in $B$.

That $B_{rss}$ is open in $B_{ss}$ follows from Lemma~\ref{B1b} and the upper semi-continuity theorem for cohomology. 

 Let $r(\mathcal{E}) = rh,\; d(\mathcal{E}) = dh$ where $(r,d)=1$. Let $ 1\leq \ell \leq h$. Consider the vector bundle 
$\mathcal{H}=\mathcal{E}(r\ell,d\ell)^{\vee}\otimes 
p_{1}^{*}\mathcal{E}$  over $\mathcal{Y}\times_{B} \mathcal{Y}$.
From the upper semi-continuity theorem for 
cohomology applied to $\mathcal{H}$ and 
$p_{2}: \mathcal{Y}\times_{B} \mathcal{Y}\to \mathcal{Y}$
it follows that the set $S_{\ell}=\{z\in 
\mathcal{Y}|h^0(\mathcal{Y}_{b},\mathcal{E}(r\ell,d\ell)^{\vee}_{z}\otimes 
\mathcal{E}_{b})\geq 2\quad \text{for}\quad b=q(z)\}$ is closed in 
$\mathcal{Y}$. Since $q:\mathcal{Y}\to B$ is proper  the image $A_{\ell}=q(S_{\ell})$  is closed in $B$. According to Lemma~\ref{B1c}  the set $B_{rps}$ is the complement in $B_{ss}$ of $\cup_{\ell = 1}^h A_{\ell}$. Therefore $B_{rps}$ is open in $B$.
\hfill $\Box$

\bibliographystyle{amsalpha}
\providecommand{\bysame}{\leavevmode\hbox to3em{\hrulefill}\thinspace}

\bigskip

\noindent
{\sc
Department of Mathematics, University of Palermo\\
Via Archirafi n. 34, 90123 Palermo, Italy}\qquad and 
\par \smallskip \noindent
{\sc Institute of Mathematics, Bulgarian Academy of Sciences}

\medskip \noindent
{\it E-mail address:} {\bf kanev@math.unipa.it}

\end{document}